\documentclass[11pt]{article}

\usepackage[left=2.7cm,top=2cm,right=2.7cm,nohead]{geometry}
\usepackage{amsmath,amssymb,amsthm}
\usepackage{enumerate}

\usepackage{pifont} 

\newtheorem{theorem}{Theorem}[section]
\newtheorem{lemma}[theorem]{Lemma}
\newtheorem{corollary}[theorem]{Corollary}
\newtheorem{remark}[theorem]{Remark}
\newtheorem{assumption}[theorem]{Assumption}

\newcommand{\spgl}{{\sc{spgl}\footnotesize{1}}}

\newcommand{\spg}{{\sc{spg}}}
\newcommand{\abbrv}[1]{{\small{#1}}}

 \usepackage{graphicx}
 \usepackage{color}
 \usepackage[colorlinks,citecolor={black},urlcolor={black},linkcolor={black}]{hyperref}
 \usepackage{url}
 \input{macros}

\usepackage[nofillcomment]{algorithm3e}
\SetAlgoSkip{}\SetAlgoInsideSkip{smallskip}
\renewcommand{\AlTitleFnt}[1]{\textbf{#1}\unskip}
\newenvironment{algo}[1]
{\def\theAlgoLine{\arabic{AlgoLine}}
 \begin{algorithm}[#1]%
    \small%
    \dontprintsemicolon%
    \SetArgSty{texttsf}%
    \SetTitleSty{textsf}{}%
    \SetKwInput{Inputs}{Input}%
    \SetKwInput{Outputs}{Output}%
    \SetKwData{Converged}{converged}%
    \SetKwComment{tcc}{[}{]}
    \SetKwFunction{BuildHeap}{BuildHeap}
    \SetKwFunction{DeleteMax}{DeleteMax}
    \SetKwFunction{SoftThreshold}{SoftThreshold}
  }
  {\end{algorithm}}

\definecolor{softblue}{rgb}{0.95,0.95,0.95}
\definecolor{gray}{rgb}{0.5,0.5,0.5}
\definecolor{comment}{rgb}{0,0,0}

\usepackage{boxedminipage}
\usepackage{array}

\makeatletter\newenvironment{btheorem}{%
\begin{lrbox}{\@tempboxa}\begin{minipage}{0.97\textwidth}\begin{theorem}}%
{\end{theorem}\end{minipage}\end{lrbox}%
\par\hbox{}\noindent%
{\setlength{\fboxsep}{0pt}\colorbox{softblue}{\setlength{\fboxsep}{4pt}\begin{boxedminipage}{\textwidth}\usebox{\@tempboxa}\end{boxedminipage}}}%
\vspace{0.5\baselineskip}}
\makeatother

\makeatletter\newenvironment{blemma}{%
\begin{lrbox}{\@tempboxa}\begin{minipage}{0.97\textwidth}\begin{lemma}}%
{\end{lemma}\end{minipage}\end{lrbox}%
\par\hbox{}\noindent%
{\setlength{\fboxsep}{0pt}\colorbox{softblue}{\setlength{\fboxsep}{4pt}\begin{boxedminipage}{\textwidth}\usebox{\@tempboxa}\end{boxedminipage}}}%
\vspace{0.5\baselineskip}}
\makeatother

\makeatletter%
{\end{corollary}\end{minipage}\end{lrbox}%
\par\hbox{}\noindent%
{\setlength{\fboxsep}{0pt}\colorbox{softblue}{\setlength{\fboxsep}{4pt}\begin{boxedminipage}{\textwidth}\usebox{\@tempboxa}\end{boxedminipage}}}%
\vspace{0.5\baselineskip}}
\makeatother

\makeatletter%
{\end{proposition}\end{minipage}\end{lrbox}%
\par\hbox{}\noindent%
{\setlength{\fboxsep}{0pt}\colorbox{softblue}{\setlength{\fboxsep}{4pt}\begin{boxedminipage}{\textwidth}\usebox{\@tempboxa}\end{boxedminipage}}}%
\vspace{0.5\baselineskip}}
\makeatother

\makeatletter\newenvironment{bassumption}{%
\begin{lrbox}{\@tempboxa}\begin{minipage}{0.97\textwidth}\begin{assumption}}%
{\end{assumption}\end{minipage}\end{lrbox}%
\par\hbox{}\noindent%
{\setlength{\fboxsep}{0pt}\colorbox{softblue}{\setlength{\fboxsep}{4pt}\begin{boxedminipage}{\textwidth}\usebox{\@tempboxa}\end{boxedminipage}}}%
\vspace{0.5\baselineskip}}
\makeatother

\makeatletter%
{\end{remark}\end{minipage}\end{lrbox}%
\par\hbox{}\noindent%
{\setlength{\fboxsep}{0pt}\colorbox{softblue}{\setlength{\fboxsep}{4pt}\begin{boxedminipage}{\textwidth}\usebox{\@tempboxa}\end{boxedminipage}}}%
\vspace{0.5\baselineskip}}
\makeatother

\newcommand{\sfrac}[2]{{\textstyle\frac{#1}{#2}}}

\newcommand{\prox}{\mathrm{prox}}
\newcommand{\dhull}[1]{\mathrm{diff\ hull({#1})}}
\newcommand{\selfproj}[1]{\mathrm{self\ proj}({#1})}

\let\OrigCaption\caption
\renewcommand{\caption}[1]{\OrigCaption{\small\it{#1}}}

\title{A Hybrid Quasi-Newton Projected-Gradient Method \\with
  Application to Lasso and Basis-Pursuit Denoise} \author{E. van den
  Berg\\[3pt]\small IBM Watson, Yorktown Heights, NY, USA}

\begin{document}

\maketitle

\begin{abstract}
  We propose a new algorithm for the optimization of convex functions
  over a polyhedral set in $\mathbb{R}^n$. The algorithm extends the
  spectral projected-gradient method with limited-memory BFGS iterates
  restricted to the present face whenever possible. We prove
  convergence of the algorithm under suitable conditions and apply the
  algorithm to solve the Lasso problem, and consequently, the
  basis-pursuit denoise problem through the root-finding framework
  proposed by van den Berg and Friedlander [SIAM Journal on Scientific
  Computing, 31(2), 2008]. The algorithm is especially well suited to
  simple domains and could also be used to solve bound-constrained
  problems as well as problems restricted to the simplex.
\end{abstract}

\section{Introduction}\label{Sec:Introduction}
In this paper we propose an algorithm for optimization problems of the
form
\begin{equation}\label{Eq:GeneralProblem}
\minimize{x}\quad f(x)\quad \st\quad x \in \mathcal{C},
\end{equation}
where $f : \mathbb{R}^n \to \mathbb{R}$ is a convex, twice
continuously differentiable function, and $\mathcal{C}$ is a
polyhedral set in $\mathbb{R}^n$. The main focus of the paper is the
specialization and application of the framework to the Lasso problem
\cite{TIB1996a}:
\begin{equation}\label{Eq:Lasso}
\minimize{x}\quad \half\norm{Ax-b}_2^2\quad \st\quad
\norm{x}_1\leq\tau,\tag{{\small LS}$_\tau$}
\end{equation}
where $\mathcal{C}$ is a, possibly weighted, one-norm ball and $f(x) =
\half\norm{Ax-b}_2^2$. The work in this paper was motivated by the
need for an efficient and accurate solver for the Lasso subproblems
appearing in the \spgl~\cite{BER2008Fb} solver for basis-pursuit denoise
\cite{CHE1998DSa} problems of the form
\begin{equation}\label{Eq:BPDN}
\minimize{x}\quad \norm{x}_1\quad \st\quad \half \norm{Ax-b}_2^2 \leq
\sigma.\tag{{\small BP}$_\sigma$}
\end{equation}
Both formulations are central to compressed sensing
\cite{CAN2006RTa,DON2006c} as a means of recovering exactly or
approximately sparse vectors $x_0$ from linearly compressed and often
noisy observations $b = Ax_0 + z$. In practice we may have a better
idea about about the noise level $\norm{z}_2$, appearing as
$\sigma$ in the basis-pursuit denoise formulation, rather than the
one-norm of the unknown signal $x_0$, appearing as $\tau$ in Lasso.
The \eqref{Eq:BPDN} formulation is therefore often a more natural
choice.

\begin{figure}[t]
\centering
\includegraphics[width=0.645\textwidth]{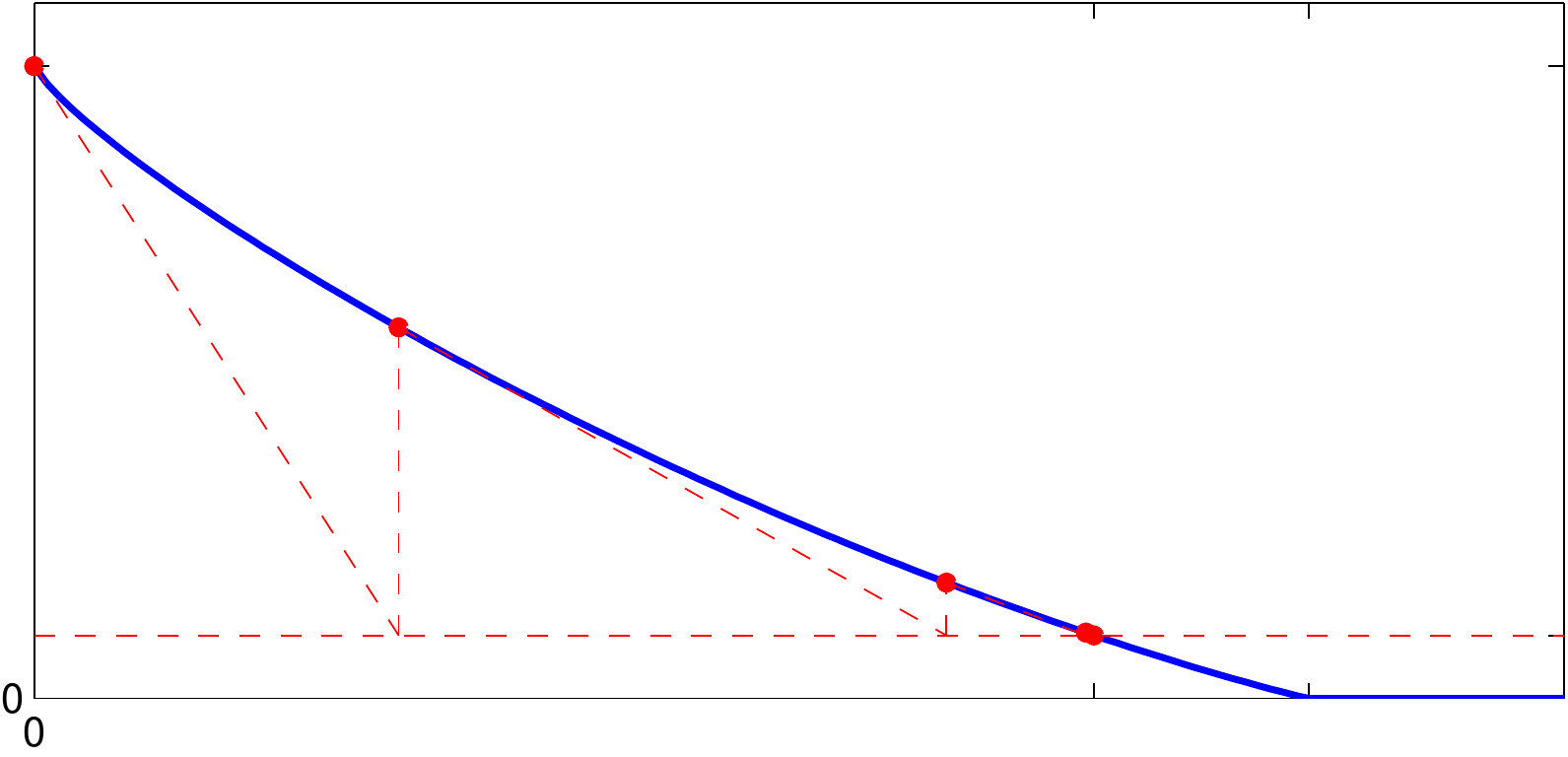}%
\begin{picture}(0,0)(0,0)
  \put(-168,59){\small$f(\tau)$}
  \put(-93,2.5){\small$\tau_{\sigma}$}
  \put(-52,2.5){\small$\tau_{0}$}
  \put(-300,20.5){\small$\sigma$}
\end{picture}
\vspace*{-8pt}
\caption{Example of root finding on the Pareto curve $f(\tau)$.}\label{Fig:RootFinding}
\end{figure}
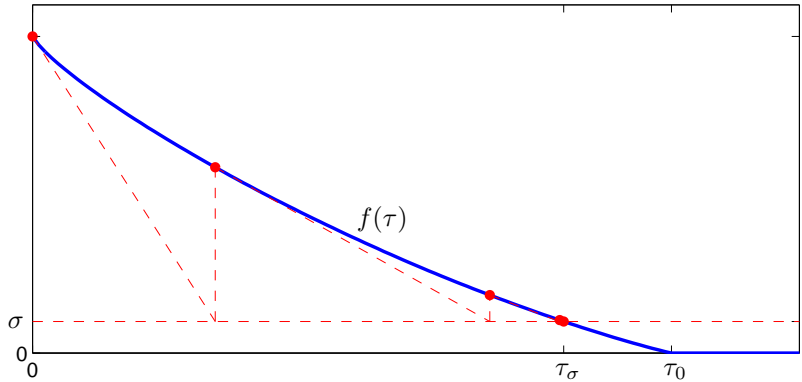

It was shown in \cite{BER2008Fb} that basis-pursuit denoise and Lasso are
connected through the Pareto curve
\[
f(\tau) = \min_{x}\quad \norm{Ax-b}_2 \quad \st\quad \norm{x}_1
\leq \tau,
\]
and that solving \eqref{Eq:BPDN} can be reduced to finding the
smallest $\tau$ for which the Lasso solution $x_{\tau}^*$ satisfies
$\norm{Ax_{\tau}^*-b}\leq\sigma$. Denoting by $\tau_{\sigma}$ this
critical value of $\tau$ and assuming that $b$ lies in the range space
of $A$ it was shown in \cite{BER2008Fb} that the Pareto curve is
convex and differentiable at all $\tau \in [0,\tau_0)$ with gradient
$\norm{A^Tr}_{\infty} / \norm{r}_2$ where $r$ denotes the misfit
$Ax_{\tau}^* - b$. Evaluation of both $f(\tau)$ and $f'(\tau)$ relies
on the misfit $r$, which can be obtained by solving \eqref{Eq:Lasso}.
The \spgl~solver proposed in \cite{BER2008Fb} applies root finding on
the Pareto curve, as illustrated in Figure~\ref{Fig:RootFinding}, to
solve $f(\tau) = \sigma$ and thereby reduce basis-pursuit denoise to a
series of Lasso problems. In \spgl\ these subproblems are solved using
the spectral projected-gradient (\abbrv{SPG})
algorithm~\cite{BIR2000MRa}, which we discuss in more detail in
Section~\ref{Sec:Background}. For certain problems it was found that
\abbrv{SPG} generates long sequences of iterates that all lie on the
same face of the feasible set.  This suggests an active-set type of
method in which a quasi-Newton method, such as the limited-memory
\abbrv{BFGS} (\abbrv{L-BFGS}) method \cite{LIU1989Na}, is used to
minimize the problem restricted to the current face. Rather than
carefully deciding when an active set has stabilized and then
accurately solving over the active set before switching back to the
global mode, we propose a hybrid algorithm that uses seamless and
lightweight switching between the two methods. By doing so, we are
able to take full advantage of the strengths of both methods, while
avoiding possibly costly subproblem solves, or complicated heuristics
that determine when to switch between the solvers.

\subsection{Paper outline}

In Section~\ref{Sec:Background} we provide a concise background on the
\abbrv{SPG} and \abbrv{L-BFGS} methods along with some of their theoretical
properties. We then describe the proposed algorithm for the general
problem formulation \eqref{Eq:GeneralProblem} in
Section~\ref{Sec:ProposedAlgorithm}. In Section~\ref{Sec:Applications}
we study the geometry of the constraints in the Lasso problem, and
develop the tools needed for an efficient implementation of the
framework for Lasso. Numerical experiments are provided in
Section~\ref{Sec:NumericalExperiments}, followed by the conclusions in
Section~\ref{Sec:Conclusions}.

\subsection{Notations and definitions}

We use caligraphic capital letters for sets. Given any two set
$\mathcal{S}_1$ and $\mathcal{S}_2$, we write $\mathcal{S}_1 +
\mathcal{S}_2$ for $\{x_1 + x_2 \mid x_1 \in \mathcal{S}_1,\
x_2\in\mathcal{S}_2\}$, and likewise for $\mathcal{S}_1 -
\mathcal{S}_2$. For a seeming lack of established terminology, we
define the \emph{difference hull} of a set $\mathcal{S}$,
$\dhull{\mathcal{S}}$, as the linear hull of differences $\{u_1 - u_2
\mid u_1,u_2\in\mathcal{S}\}$. The difference hull can be seen as the
linear subspace corresponding to the affine hull of $\mathcal{S}$
translated to contain the origin. For any $x$ in a polyhedral set
$\mathcal{C}$, we define $\mathcal{F}(x)$ to be the unique face
$\mathcal{F}$ of $\mathcal{C}$ for which $x \in
\mathrm{relint}(\mathcal{F})$; this may be $\mathcal{C}$ itself. The
normal cone of $\mathcal{C}$ at $x$ is given by $ \mathcal{N}(x) :=
\{d \in \mathbb{R}^n \mid \mathcal{P}(x + d) = x\}$. The normal cone
of a face $\mathcal{F}$ is understood to be $\mathcal{N}(x)$ for any
$x \in\mathrm{relint}(\mathcal{F})$.  Orthogonal projection of any
vector $v \in\mathbb{R}^n$ onto $\mathcal{C}$ is defined as
\[
\mathcal{P}(v) := \argmin_{x}\quad \norm{x-v}_2\quad \st\quad x\in\mathcal{C}.
\]
 We define
the \emph{self-projection cone} of a face $\mathcal{F} =
\mathcal{F}(x)$ as the closed and convex cone of directions
$d\in\mathbb{R}^n$ such that there exists an $\epsilon > 0$ for which
the projection of $x + \epsilon d$ lies on $\mathcal{F}$:
\[
\selfproj{\mathcal{F}} = \mathcal{S}(\mathcal{F}(x)) := \{ d \in
\mathbb{R}^n \mid \exists\ \epsilon > 0: \mathcal{F}[\mathcal{P}(x +
\epsilon d)] = \mathcal{F}(x)\} = \mathcal{N}(x) +
\dhull{\mathcal{F}(x)}.
\]
Note that $\mathcal{N}(x) \perp \dhull{\mathcal{F}(x)}$; in fact, the
difference hull of $\mathcal{F}$ is the orthogonal complement of the
linear hull of $\mathcal{N}(\mathcal{F})$. For any $k$-face
$\mathcal{F}$ of $\mathcal{C}$, $k\geq 1$, we denote by
$\Phi(\mathcal{F}) \in \mathbb{R}^{n \times k}$ an arbitrary but fixed
orthonormal basis for $\dhull{\mathcal{F}}$. We will never use
$\Phi(\mathcal{F})$ when $\mathcal{F}$ is a vertex and therefore leave
it undefined. We denote by $e_i$ the $i$-th column of an identity
matrix whose size is clear from the context. The proximation of a
function $f$ is defined as $\prox_f(u) := \argmin_x f(x) + \half\norm{x-y}_2^2$.


\section{Background}\label{Sec:Background}
\subsection{The nonmonotone spectral projected-gradient method}\label{Sec:SPG}

The nonmonotone spectral projected-gradient method (\abbrv{SPG}) was
introduced by Birgin, Mart\'{i}nez, and Raydan \cite{BIR2000MRa} for
problems of the form \eqref{Eq:GeneralProblem}, with $\mathcal{C}$ a
closed convex set in $\mathbb{R}^n$, and $f: \mathbb{R}^n \to
\mathbb{R}$ a function with continuous partial derivatives on an
open set that contains $\mathcal{C}$. The algorithm is outlined in
Algorithm~\ref{Alg:SPG}, and it can be seen that the main step in each
iteration is a line search along the curvilinear trajectory given by (see also
\cite{BER1976a}):
\begin{equation}\label{Eq:Curvilinear}
x(\alpha) = \mathcal{P}(x_i - \alpha \nabla f(x_i)),\quad \alpha \geq 0.
\end{equation}

\begin{algo}{bt}
\caption{The nonmonotone spectral projected-gradient method (\abbrv{SPG}).}\label{Alg:SPG}
\Inputs{$A$,\ $b$,\ $\tau$,\ $x_0$}
Initialize\ \,$i \gets 0$,\ choose\ \,$\alpha_0 \in [\alpha_{\min},\alpha_{\max}]$\;
\While{not done}{\smallskip
Compute\ \,$g_i = \nabla f(x_i)$\;

\medskip
{\color{comment}\textit{\# Compute the Barzilai-Borwein scaling parameter}}\;
\If{$i > 0$}{\smallskip
  Set\ \,$s \gets x_i - x_{i-1}$,\ \,$y \gets g_i - g_{i-1}$\;\smallskip
   \eIf{$s^Ty > 0$}
   {  $\alpha_i \gets \max\left\{\min\left\{\frac{s^Ts}{s^Ty},\;\alpha_{\max}\right\},\;\alpha_{\min}\right\}$\;
   }
   {  $\alpha_i \gets \alpha_{\max}$\;
   }
}

\bigskip
{\color{comment}\textit{\# Non-monotone curvilinear Armijo line-search}}\;
Initialize\ \,$k \gets 0$\;
\While{condition \eqref{Eq:NonmonotoneArmijo} is not satisfied}{\smallskip
$k \gets k + 1$\;
}
\medskip
{\color{comment}\textit{\# Proceed to the next iteration}}\;
Set\ \,$x_{i+1} = x(\beta^k\alpha_i)$\;
Update\ \,$i \gets i + 1$\;
}
\end{algo}

Two important modifications to the curvilinear projected-gradient
method, made to help speed up convergence, were introduced in
\cite{BIR2000MRa}. The first modification allows a limited level of
nonmonotonicity in the objective value. Given {$\mu,\gamma \in
  (0,1)$}, the Armijo-type line search starts with an initial step
length $\alpha_i$, and then finds the first nonnegative integer $k$
such that
\begin{equation}\label{Eq:NonmonotoneArmijo}
f(x(\mu^k \alpha_i)) \leq
 \max_{0 \leq j \leq\min\{i,M-1\} }\{f(x_{i-j})\} +
\gamma (\nabla f(x_i))^T( x(\mu^k\alpha_i) - x_i).
\end{equation}
The right-hand side of this condition ensures sufficient descend, but
only with respect to the maximum of up to $M$ of the most recent
objective values. In case $M=1$ this reduces to the standard Armijo
line-search condition. The second modification is the use of the
spectral step length, as proposed by Barzilai and Borwein
\cite{BAR1988Ba}. Given $s = x_i - x_{i-1}$ and $y = \nabla f(x_i) -
\nabla f(x_{i-1})$, the initial step length at iteration $i$ is defined as
\[
\alpha_{i} = \begin{cases}
\max\left\{\min\left\{\frac{s^Ts}{s^Ty},\ \alpha_{\max}\right\},\ \alpha_{\min}\right\} &
\mathrm{if}\ s^Ty > 0,
\\
\alpha_{\max} & \mathrm{otherwise},
\end{cases}
\]
where $0 < \alpha_{\min} < \alpha_{\max}$ are fixed parameters. More
information on the motivation behind this particular choice of step
length can be found in \cite{BIR2000MRa,FLE2012a}. Under the
conditions stated at the beginning of the section, it holds that any
accumulation point $x^*$ of the sequence $\{x_i\}$ is a constrained
stationary point \cite{BIR2000MRa}; that is a point $x^*
\in\mathcal{C}$ such that
\begin{equation}\label{Eq:StationaryPoint}
\norm{\mathcal{P}(x^* - \nabla f(x^*)) - x^*}_2 = 0.
\end{equation}
In practice a relaxed version of \eqref{Eq:StationaryPoint} is used as
a stopping criterion in Algorithm~\ref{Alg:SPG}, along with other
conditions.

\subsection{Limited-memory BFGS}

The \abbrv{L-BFGS} algorithm by Liu and Nocedal~\cite{LIU1989Na} is a popular
quasi-Newton method for unconstrained minimization of smooth functions
$f: \mathbb{R}^n \to \mathbb{R}$:
\[
\minimize{x}\quad f(x).
\]
At each iteration, the algorithm constructs a positive definite
approximation $H_i$ of the inverse Hessian of $f$ at $x_i$. This
construction is based on an initial positive definite matrix $H$ and
$\hat{n} = \min\{i,N\}$ of the most recent vector pairs $\{s_{i-j}$,
$y_{i-j}\}_{j=0}^{\hat{n}-1}$, with
\begin{equation}\label{Eq:StandardUpdate}
s_j = x_j - x_{j-1},\quad \mathrm{and}\quad y_j = \nabla f(x_j) -
\nabla f(x_{j-1}).
\end{equation}
The iterates are of the form $x_{i+1} = x_{i} + \alpha_i d_i$, where the search
direction $d_i$ is given by
\begin{equation}\label{Eq:QNSearchDirection}
d_{i} = -H_i \cdot\nabla f(x_i),
\end{equation}
and the step size $\alpha_i$ is chosen to satisfy the Wolfe conditions:
\begin{subequations}\label{Eq:Wolfe}
\begin{equation}\label{Eq:Wolfe1}
f(x_i + \alpha_i d_i) \leq f(x_i) + \gamma_1 \alpha_i (\nabla
f(x_i))^T d_i,
\end{equation}
\begin{equation}\label{Eq:Wolfe2}
(\nabla f(x_i + \alpha_i d_i))^Td_i \geq \gamma_2 (\nabla f(x_i))^Td_i.
\end{equation}
\end{subequations}
Parameters $\gamma_1$ and $\gamma_2$ are chosen such that $0 <
\gamma_1 < \half$, and $\gamma_1 < \gamma_2 < 1$. For details on the
structure of the inverse approximation $H_i$ and efficient ways of
evaluating the matrix-vector product in \eqref{Eq:QNSearchDirection},
see \cite{LIU1989Na,NOC2006Wa}.

\subsubsection{Convergence results}
For the analysis of the \abbrv{L-BFGS} algorithm, Liu and Nocedal
make the following assumptions \cite{LIU1989Na}:

\begin{bassumption}\label{Ass:LIU1989Na}
For a given starting point $x_0$, we have that:\\[2pt]
(1) The objective function $f$ is twice continuously differentiable;\\[2pt]
(2) The level set $\mathcal{D} := \{ x \in\mathbb{R}^n \mid f(x) \leq
f(x_0)\}$ is convex; and\\[2pt]
(3) There exist positive constants $\mu_1$ and $\mu_2$ such that
\begin{equation}\label{Eq:LIU1989Na-Eq.7.4}
\mu_1 \norm{v}^2 \leq v^T [\nabla^2 f(x)]\, v \leq \mu_2\norm{v}^2,\quad
\end{equation}
\phantom{(3)} for all $x \in \mathcal{D}$ and $v \in \mathbb{R}^n$.
\end{bassumption}
\noindent Under these conditions, and with some simplifications, they prove that
\begin{btheorem}[Liu and Nocedal \cite{LIU1989Na}]
  The \abbrv{L-BFGS} algorithm generates a sequence $\{x_i\}$ that converges
  to the unique minimizer $x_*$ in $\mathcal{D}$. Moreover, there
  exists a constant $c > 0$ such that
\begin{equation}\label{Eq:LIU1989Na-Eq.7.13b}
f(x_{i+1}) - f(x^*) \leq (1 - c)(f(x_i) - f(x^*)).
\end{equation}
\end{btheorem}


\section{Proposed algorithm}\label{Sec:ProposedAlgorithm} 
The proposed algorithm can be seen as a modification of the
\abbrv{SPG} method that allows the use of quasi-Newton steps over a
currently active face.  The basic idea is that whenever two successive
iterates $x_i$ and $x_{i-1}$ lie on the same face, we can form or
update a quadratic model of the objective function restricted to the
face. To avoid unnecessary updates to the model and, indeed, to ensure
convergence to a global optimizer, we require that $-g_i := -\nabla
f(x_i)$ lies in the self-projection cone of $\mathcal{F}_i :=
\mathcal{F}(x_i)$. Whenever a model for the current face is available,
the algorithm will attempt a quasi-Newton step that is restricted to
the face and satisfies the Wolfe conditions \eqref{Eq:Wolfe}. If the
quasi-Newton step fails, or is otherwise abandoned, the algorithm
simply falls back and takes a regular \abbrv{SPG} step. After each
step---regardless of the type---we again check the conditions required
to update the quadratic model and initiate the quasi-Newton step:
\begin{equation}\label{Eq:UpdateCriterion}
\mathcal{F}(x_{i}) = \mathcal{F}(x_{i-1})\quad\mbox{and}\quad -g_i \in \selfproj{\mathcal{F}_i}.
\end{equation}
If these conditions are not met, we discard the Hessian approximation
used in the quadratic model, for example, by setting it to the empty
set. Note that omitting the self-projection criterion from
\eqref{Eq:UpdateCriterion} could cause the algorithm to take repeated
quasi-Newton iterations that converge to a minimum on the relative
interior of the face that is not the global minimum.

\subsection{Quasi-Newton over a face}

One way of performing quasi-Newton over a face is by maintaining an
inverse Hessian approximation using the update vectors in
\eqref{Eq:StandardUpdate}, and computing the search direction $d_i$
using \eqref{Eq:QNSearchDirection}. However, this approach has some
major disadvantages. First, we may have that $d_i \not\in
\dhull{\mathcal{F}_i}$, which means that $x_i + \alpha d_i
\not\in\mathcal{F}_i$ for all nonzero $\alpha$. This could be
partially solved by projection onto the face, but such a projected
direction is no longer guaranteed to be a descent
direction~\cite{BER2003a}. This too could be addressed by modifying
the Hessian, but doing so would further complicate the algorithm.  A
second disadvantage is that we maintain the inverse Hessian
approximation for the ambient space, which typically has a much higher
dimension than the current face and may therefore not be very
accurate along $\mathrm{aff}(\mathcal{F}_i)$.

\begin{figure*}[!t]
\begin{algo}{H}
\caption{Outline of the proposed hybrid quasi-Newton projected-gradient method.}\label{Alg:HybridMethod}
\Inputs{$A$,\ $b$,\ $\tau$,\ $x_0$}

Initialize \quad$H_0 = \emptyset$,\quad $\mathcal{F}_0 = \mathcal{F}(x_0)$,\quad $g_0 = \nabla f(x_0)$\;
Set $t \gets 0$,\ \,choose $\alpha_0 \in
[\alpha_{\min},\alpha_{\max}]$\;
\While{not done}{

\smallskip
{\color{comment}\textit{\# Quasi-Newton step on current face}}\;
flagUpdated $\gets$ false\;
\If{($H_t \neq \emptyset$)}
{  $d_t = -\Phi_t H_k\Phi_t^T g_t$\;
   Perform Wolfe-type line-search on $x(\gamma) := x + \gamma d$ over
   $\{\gamma \mid x(\gamma) \in \mathcal{F}_t\}$\;
   \If{(line search successful)}
   {  Reset objective function history\;
      $x_{t+1} = x_t + \gamma_t d_t$\;
      flagUpdated $\gets$ true\;\smallskip
   }
} 

\medskip
{\color{comment}\textit{\# Projected-gradient step}}\;
\If{(flagUpdated $=$ false)}
{  
   Compute Barzilai-Borwein scaling parameter\;
   Nonmonotone curvilinear Armijo line-search along $x(\gamma) := \mathcal{P}(x_t
- \gamma g_t)$\;
   $x_{t+1} = x(\gamma_t)$\;\smallskip
}

\medskip
$t \gets t + 1$,\quad $g_t = \nabla f(x_t)$, \quad $\mathcal{F}_{t} =
\mathcal{F}(x_{t})$\;

\medskip
{\color{comment}\textit{\# Update the quadratic model of the current face}}\;
\eIf{($\mathcal{F}_{t} = \mathcal{F}_{t-1}$ and $-g_t \in \selfproj{\mathcal{F}_t}$)}
{ \eIf{($H_{t-1} = \emptyset$)}
  { Initialize $H_{t-1} \gets \mu I$\;
    Determine orthonormal basis $\Phi_t = \Phi(\mathcal{F}_t)$\;
  }
  {
      Set $\Phi_t = \Phi_{t-1}$\;
  }
  Set $H_t$ as L-BFGS update to $H_{t-1}$ using $s_t = \Phi_t^T(x_{t} - x_{t-1})$
  and $y_t = \Phi_t^T(g_{t} - g_{t-1})$\;
}
{  $H_t = \emptyset$\;\smallskip
}

} 
\end{algo}\hfill
\end{figure*}

The solution of the above problem is straightforward: we simply work
with a representation for the function restricted to
$\mathrm{aff}(\mathcal{F}_i)$. Let $\mathcal{F}_i$ be a
$k$-dimensional face with $k > 0$. Then we can find an orthonormal
basis $\Phi := \Phi(\mathcal{F}_i) \in \mathbb{R}^{n\times k}$ whose
span coincides with $\dhull{\mathcal{F}_i}$. Using $\Phi$ we can write
any point $x\in\mathcal{F}_i$ as $x = v + \Phi c$, where $v
\in\mathbb{R}^n$ is an arbitrary but fixed point in $\mathcal{F}_i$,
and $c \in\mathbb{R}^k$ is a coefficient in the lower-dimensional
space. The function $\hat{f}:\mathbb{R}^k\to\mathbb{R}$, which
restricts $f$ to the current face, is then given by $\hat{f}(c) = f(v
+ \Phi c)$. The idea then is to form the inverse Hessian approximation
over $\hat{f}$, and use it to obtain a search direction $\hat{d}
\in\mathbb{R}^k$, which can then be mapped back to the ambient space
for the actual line search. In particular, we can form the approximate
inverse Hessian $H_i \in \mathbb{R}^{k\times k}$ by updating an
initial positive definite $H$ using
\begin{equation}\label{Eq:LBFGSUpdate}
s_i = \Phi^T(x_i - x_{i-1})\quad \mbox{and}\quad
y_i = \Phi^T(g_i - g_{i-1}).
\end{equation}
In order to obtain the search direction we first compute $\nabla
\hat{f}(\Phi^T(x_i - v)) = \Phi^T g_i$ by projecting the gradient
$g_i$ onto the lower-dimensional space. We then apply the inverse
Hessian followed by back projection, giving:
\[
d_i = \Phi H_i\Phi^T(-g_i).
\]
The resulting vector clearly lies in $\dhull{\mathcal{F}_i}$ and we
therefore have that $x_i + \alpha d_i \in \mathcal{F}_i$ for all step
sizes $\alpha \in [0,\alpha_{\mathrm{bnd}}]$, for some
$\alpha_{\mathrm{bnd}} > 0$, possibly with $\alpha_{\mathrm{bnd}} =
+\infty$. Since the line search is done in the ambient dimension, we
try to find a step size $\alpha$ within the above range, such that the
original Wolfe conditions \eqref{Eq:Wolfe} are met. For the line
search we could start with a unit step length, whenever
$\alpha_{\mathrm{bnd}} \geq 1$, or we could try $\alpha =
\alpha_{\mathrm{bnd}}$ first to encourage exploring lower-dimensional
faces, provided of course that $\alpha_{\mathrm{bnd}} < \infty$. If no
suitable step length can be found, or a certain maximum number of
trial steps is taken, we abandon the quasi-Newton step and take a
spectral projected-gradient step instead. As a final remark, note that
condition \eqref{Eq:UpdateCriterion} should never be met for vertices,
since that would imply not only that $x_i = x_{i-1}$, but also that
$-g_{i-1} \in\selfproj{\mathcal{F}_{i-1}} =
\mathcal{N}(\mathcal{F}_{i-1})$, which means that the optimality
condition given in \eqref{Eq:StationaryPoint} would already have been
satisfied at $x_{k-1}$.

\subsection{Convergence}

For the convergence analysis of Algorithm~\ref{Alg:HybridMethod} we
rely on the results in \cite{BIR2000MRa} and \cite{LIU1989Na}, and add
a step in the algorithm that resets the objective-value history used
by \abbrv{SPG} after each series of successful quasi-Newton iterations
to ensure that any subsequent iteration has a lower objective
value. We use the following assumptions, which are somewhat more
restrictive than those in the aforementioned two papers (see, for
example, Assumption~\ref{Ass:LIU1989Na}):

\begin{bassumption}\label{Ass:General} We assume that \newline
(1) The objective function $f$ is convex, twice continuously
differentiable, and bounded below;\newline
(2) There exist constants $0 < \mu_1 \leq \mu_2 < \infty$ such that
for all $x,v \in \mathbb{R}^n$:
\begin{equation}\label{Eq:LIU1989Na-7.4}
\mu_1\norm{v}_2^2 \leq v^T H(x) v \leq \mu_2\norm{v}_2^2.
\end{equation}
\end{bassumption}

\noindent Under these assumptions, we have the following result:

\begin{btheorem}
  Let $f(x)$ satisfy Assumptions~\ref{Ass:General} and let $x_0 \in
  \mathcal{C}$. Then the sequence $\{x_t\}$ generated by
  Algorithm~\ref{Alg:HybridMethod} converges to the unique minimizer
  of \eqref{Eq:GeneralProblem}.
\end{btheorem}
\begin{proof}
Assumption~\ref{Ass:General} ensures the existence of a unique
minimizer $x^*$ to \eqref{Eq:GeneralProblem}, which satisfies
\[
-g(x^*) \in \mathcal{N}(x^*).
\]
If there exists a finite $t$ for which $x_t = x^*$, we are
done. Suppose, therefore that $x_t \neq x^*$ for all $t$. We consider
two cases. First, if there are finitely many quasi-Newton steps, there
must a $\bar{t}$ such that all iterations $t > \bar{t}$ are of the
projected gradient type. In this case the result follows directly from
the analysis in \cite{BIR2000MRa}. Second, consider the case where
there are infinitely many quasi-Newton steps. It can be seen that each
quasi-Newton step works towards minimizing the objective over the
affine hull of the current face $\mathcal{F}$:
\begin{equation}\label{Eq:MinAffineFace}
\minimize{x}\quad f(x)\quad \st\quad x \in \mathrm{aff}(\mathcal{F}).
\end{equation}
For any such step it follows from the analysis in Liu and Nocedal
\cite{LIU1989Na} (with minor modification to the number of update
vectors available, and working in the lower-dimensional and
unconstrained representation based on $\Phi$) that there exists a
constant $c > 0$ such that the quasi-Newton step satisfies
\begin{equation}\label{Eq:LIU1989Na-7.13b2}
  f(x_{t+1}) - f(x_{\mathcal{F}}^*) \leq (1 - c)(f(x_t) - f(x_{\mathcal{F}}^*)),
\end{equation}
where $x_{\mathcal{F}}^*$ denotes the minimizer of
\eqref{Eq:MinAffineFace}. Because the history of the $M$ most recent
objective values is reset after each successful quasi-Newton step, any
intermediate projected-gradient step will not increase the objective.
Based on this, Lemma~\ref{Lemma:FiniteFace} below, shows that the
number of quasi-Newton iterates on any $\mathcal{F}$ that does not
contain $x^*$ is finite. By polyhedrality of the domain, the number of
faces itself is bounded, and we must therefore take infinitely many
iterations on at least one face that contains $x^*$. Repeated
application of \eqref{Eq:LIU1989Na-7.13b2} then shows that the
objective value converges to $f(x_{\mathcal{F}}^*)$. Finally, it
follows from Assumption~\ref{Ass:General} that $\{x_t\}$ converges to
$x^*$.
\end{proof}

\begin{blemma}\label{Lemma:FiniteFace}
  Let $\mathcal{F}$ be a face of $\mathcal{C}$ such that $x^*
  \not\in\mathcal{F}$. Then the number of quasi-Newton steps on
  $\mathcal{F}$ taken by Algorithm~\ref{Alg:HybridMethod} is finite.
\end{blemma}
\begin{proof}
  Let $x_{\mathcal{F}}^*$ be the solution to \eqref{Eq:MinAffineFace},
  and denote by $x_{[j]}$ and $x_{[j]+1}$ the starting, respectively
  ending, point for the $j$-th quasi-Newton step on $\mathcal{F}$. It
  can be seen that
\begin{equation}\label{Eq:QNDescend}
f(x_{[j]}) - f(x_{\mathcal{F}}^*) \leq f(x_{[j-1]+1}) -
f(x_{\mathcal{F}}^*) \leq (1-c)(f(x_{[j-1]}) - f(x_{\mathcal{F}}^*)),
\end{equation}
for $j \geq 2$. This holds since any intermediate quasi-Newton
iteration can only reduce the objective, and likewise for
projected-gradient steps, as a consequence of resetting the
function-value history.

We consider two cases. First assume that $x_{\mathcal{F}}^*
\not\in\mathcal{F}$. Let $\bar{f}$ be the minimum of $f(x)$ over
$x\in\mathcal{F}$. Repeated application of \eqref{Eq:QNDescend} gives
\begin{equation}\label{Eq:QNkSteps}
  f(x_{[j]+1}) - f(x_{\mathcal{F}}^*) \leq (1-c)^j(f(x_{[1]}) - f(x_{\mathcal{F}}^*)).
\end{equation}
For sufficiently large, but finite $j$, the right-hand side in
\eqref{Eq:QNkSteps} must fall below $\bar{f} - f(x_{\mathcal{F}}^*)$,
which is strictly positive. Since every successful quasi-Newton step
results in a vector $x_{[j]+1} \in\mathcal{F}$ by construction, it
follows that the number of quasi-Newton iterates on $\mathcal{F}$ must
be bounded.

For the second case, assume that $x_{\mathcal{F}}^* \in\mathcal{F}$.
Because optimization is done over $\mathrm{aff}(\mathcal{F})$, it
holds that $-g(x_{\mathcal{F}}^*) \perp \dhull{\mathcal{F}}$.
For $-g(x_{\mathcal{F}}^*) \in \selfproj{\mathcal{F}}$, we must
therefore have $-g(x_{\mathcal{F}}^*)
\in\mathcal{N}(x_{\mathcal{F}}^*)$, but this cannot be the case since
it would imply that $x_{\mathcal{F}}^*$ is a global minimizer. (The
same holds when $x_{\mathcal{F}}^*$ lies on a lower-dimensional
subface on the boundary of $\mathcal{F}$.) Since $f$ is continuously
differentiable by assumption, it follows that $-g(x)
\not\in\selfproj{\mathcal{F}}$ for all points $x \in \mathcal{F}$
sufficiently close to $x_{\mathcal{F}}^*$.  Assumption
\ref{Ass:General} then allows us to define a sufficiently close
neighborhood as the level set $f(x) \leq \bar{f}$ over
$x\in\mathcal{F}$, where $\bar{f} > f(x_{\mathcal{F}}^*)$. Applying
the same argument we used above shows that the right-hand side of
\eqref{Eq:LIU1989Na-7.4} again falls below $\bar{f} -
f(x_{\mathcal{F}}^*)$ for sufficiently large $j$. Once this happens
all following iterates $x_t \in\mathcal{F}$ must have $f(x_t) \leq
\bar{f}$. Since the self-projection cone condition $-g(x)
\in\selfproj{\mathcal{F}}$ does not hold at these points, no more
quasi-Newton steps are taken on $\mathcal{F}$.
\end{proof}

A similar analysis holds when the spectral projected-gradient method
is replaced by another convergent algorithm, provided that the
iterates do not exceed the initial objective value.


\section{Application to Lasso}\label{Sec:Applications}
The proposed algorithm depends on a number of operations on the
constraint set. In particular, it has to determine in which face the
current iterate lies, check membership of the self-projection cone,
and determine an orthonormal basis for the current face. For the
algorithm to be of practical use, the constraint set therefore needs
to be simple enough to allow efficient evaluation of these operations.
As this work was motivated by improving the Lasso problem, we focus on
the weighted one-norm ball (which for unit weights is also known as
the cross-polytope or $n$-octahedron~\cite{GRU2003a}):
\[
\mathcal{C}_{w,1} = \{x\in\mathbb{R}^n \mid \norm{x}_{w,1} \leq\tau\},
\]
where $\norm{x}_{w,1} := \sum_i w_i\abs{x_i}$ positive $w_i$. The
proposed framework also applies naturally to bound or simplex
constrained problems, but these are outside the scope of this
paper.

The objective function we consider throughout this section is
\begin{equation}\label{Eq:ObjectiveF}
f(x) = \half\norm{Ax-b}_2^2 + {\textstyle\frac{\mu}{2}}\norm{x}_2^2
+ c^Tx,
\end{equation}
which can also be written in the form $\half\norm{Ax-b}_2^2 + c^Tx$,
with $A$ and $b$ appropriately redefined.  The benefit of having an
objective function of the form \eqref{Eq:ObjectiveF} is that it
permits closed-form expressions for step lengths satisfying certain
conditions.  In the remainder of this section we discuss practical
considerations for the line-search conditions and look at the specific
structure and properties of the set $\mathcal{C}_{w,1}$.

\subsection{Line search}

For most objective functions the line search is done by evaluating
$f(\mathcal{P}(x + \alpha d))$ or $f(x + \alpha d)$ for a series of
$\alpha$ values until all required conditions, such as Armijo and
Wolfe, are satisfied. The objective function in~\eqref{Eq:ObjectiveF}
has closed-form solutions for some of the problems arising in the line
search, thereby allowing us to simplify the algorithms and improve their
performance.

\subsubsection{Optimal unconstrained step size}

As a start we look at the step length that minimizes the objective
along $f(x+\alpha d)$:
\[
\alpha_{\mathrm{opt}} = \argmin_{\alpha}\quad f(x + \alpha d)
\]
Differentiating $f$ with respect to $\alpha$ and equating to zero
leads to the following expression:
\[
\alpha_{\mathrm{opt}} = -\frac{(A^Tr + \mu x + c)^Td}{\norm{Ad}_2^2 + \mu\norm{d}_2^2},
\]
with $r = Ax-b$. When $\mu = 0$ and $c = 0$ this reduces to
$\alpha_{\mathrm{opt}} = -r^TAd / \norm{Ad}_2^2$.

\subsubsection{Wolfe line search conditions}

The maximum step length for which the Armijo condition
\eqref{Eq:Wolfe1} is satisfied can be found by expanding the terms and
simplifying. Doing so gives the following bound:
\[
\alpha \leq \alpha_{\max} = 2(1-\gamma_1) \alpha_{\mathrm{opt}}.\label{Eq:AlphaMax}
\]
Likewise, the gradient condition \eqref{Eq:Wolfe2} reduces to
\[
\alpha \geq \alpha_{\min} = (1-\gamma_2)\alpha_{\mathrm{opt}}.
\]
The derivations of these quantities are given in
Sections~\ref{AppSec:Wolfe1} and \ref{AppSec:Wolfe2} for completeness.

\subsubsection{Projection arc}\label{Sec:PiecewiseProjectionArc}

Line search in gradient projection methods is often done by
backtracking from a single projection
\[
p(\alpha) = x + \alpha (\mathcal{P}_{\mathcal{C}}(x + d)),
\]
with $\alpha \in [0,1]$, or by search over the projection arc
\[
p(\alpha) = \mathcal{P}_{\mathcal{C}}(x + \alpha d),
\]
with $\alpha \geq 0$. The trajectory of the first method depends
strongly on the scaling of $d$ and is more likely to generate points
on the interior of the feasible set. The second method is invariant to
the scaling of $d$ and better captures the structure of the domain,
but can also be more expensive computationally. The theorem below
shows that the projection arc for polyhedral sets is piecewise linear
and continuous, with a finite number of segments. In
Sections~\ref{Sec:ProjectionArc} we provide an efficient algorithm for
generating the successive line segments of the projection arc over
$\mathcal{C}_{w,1}$. Combined with the closed-form solution of the
optimum along each segment this gives an efficient and reliable
algorithm for doing the line search along the entire projection arc.

\begin{figure}
\centering
\includegraphics[width=0.49\textwidth]{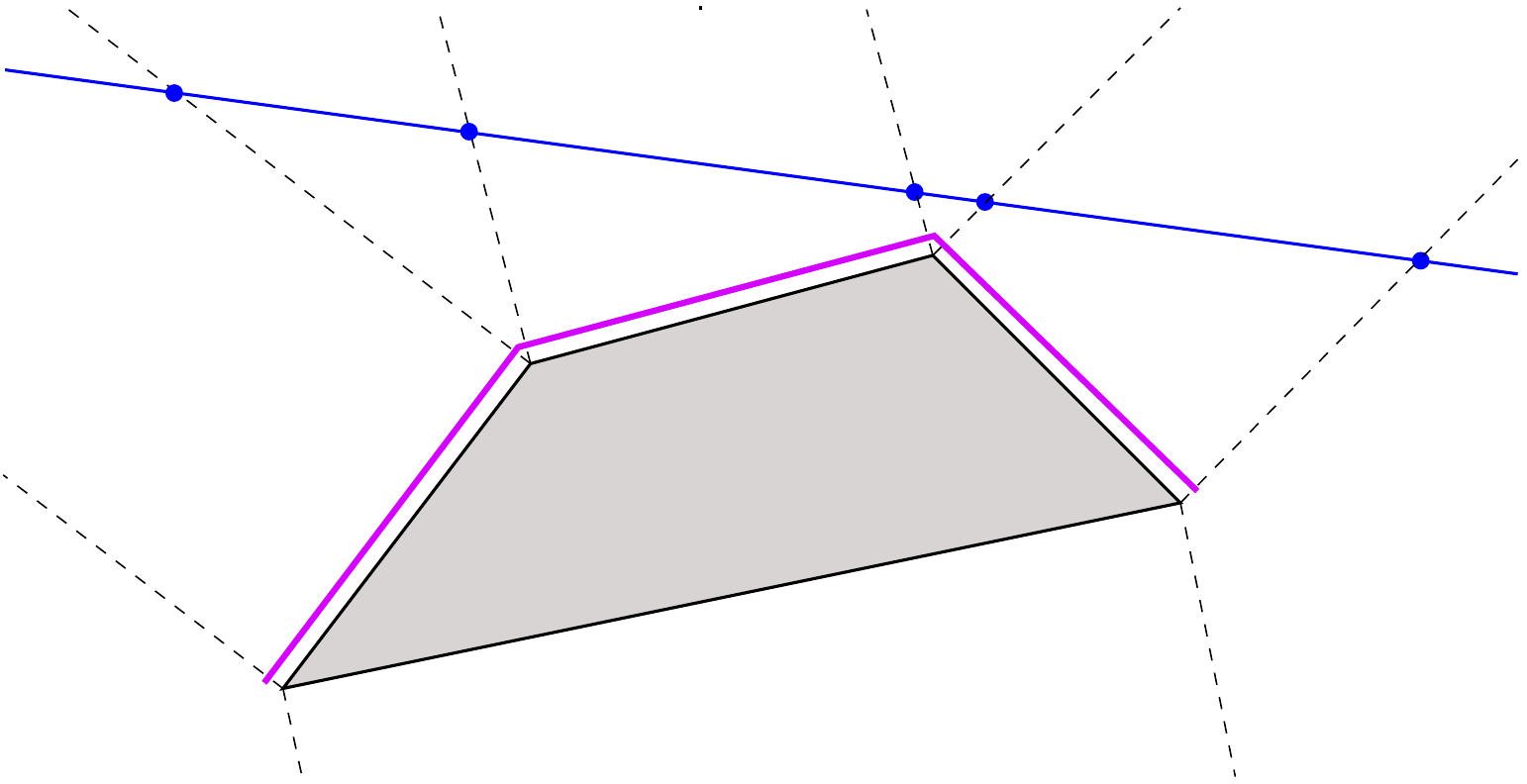}%
\begin{picture}(0,0)(0,0)
\put(-129,43){\footnotesize$\mathcal{R}_1=\mathcal{C}$}
\put(-217,10){\footnotesize$\mathcal{R}_2$}
\put(-200,57){\footnotesize$\mathcal{R}_3$}
\put(-185,108){\footnotesize$\mathcal{R}_4$}
\put(-137,108){\footnotesize$\mathcal{R}_5$}
\put(-81,108){\footnotesize$\mathcal{R}_6$}
\put(-40,95){\footnotesize$\mathcal{R}_7$}
\put(-25,38){\footnotesize$\mathcal{R}_{8}$}
\put(-115,8){\footnotesize$\mathcal{R}_{9}$}
\end{picture}
\caption{{Illustration of a polytope $\mathcal{C}$ and the associated
  partition of $\mathbb{R}^2$ into nine regions. The solid blue line
  (extended to infinity) intersects regions
  $\mathcal{R}_2$ through $\mathcal{R}_{8}$ and projects onto the
  corresponding faces, shown by the solid purple line, which is offset
  slightly from $\mathcal{C}$ for clarify.}}\label{Fig:Partition}
\end{figure}

\begin{btheorem}\label{Thm:PiecewiseLinear}
  Let $\mathcal{C}$ be a convex polyhedron in $\mathbb{R}^n$. Then for
  any $x,d \in\mathbb{R}^n$, the projection trajectory $p(\alpha) :=
  \mathcal{P}_{\mathcal{C}}(x + \alpha d)$ is piecewise linear and
  continuous with a finite number of segments.
\end{btheorem}
\begin{proof}
  The result is trivial for $d = 0$ so we assume $d \neq 0$.  Let
  $\mathcal{F}_i$, $i = 1,\ldots,N$ be the faces of
  $\mathcal{C}$. Then we partition the ambient space $\mathbb{R}^n$
  into convex polyhedral regions $\mathcal{R}_i = \mathcal{F}_i +
  \mathcal{N}(\mathcal{F}_i)$, where $\mathcal{N}(\mathcal{F}_i)$ is
  the normal cone to face $\mathcal{F}_i$ and $+$ denotes the
  Minkowski sum (see Figure~\ref{Fig:Partition} for an illustration).
  The line $x + \alpha d$ intersects the boundary of each region for
  at most two values of $\alpha$. Because of the one-to-one
  correspondence of regions and faces, it follows that the number of
  such breakpoints must be finite. Given two successive breakpoints
  $\alpha_k$ and $\alpha_{k+1}$ (possibly $\pm\infty$) corresponding
  to the intersection with some region $\mathcal{R}_i$. We first show
  that the projection trajectory $p(\alpha)$ for $\alpha_k \leq \alpha
  \leq \alpha_{k+1}$ is linear.  For regions generated by a vertex it
  holds that all points on the line segment project to the vertex,
  thus giving a constant trajectory, which is trivially linear, and it
  remains to show linearity for higher dimensional faces. Let $v =
  x+\alpha d$ with $\alpha \in [\alpha_1, \alpha_2]$. Then $v =
  \mathcal{P}_{\mathcal{C}}(v) + (v - \mathcal{P}_{\mathcal{C}}(v)) $,
  where the first and second terms lie respectively in $\mathcal{F}_i$
  and $\mathcal{N}(\mathcal{F}_i)$.  Let $Q = \Phi(\mathcal{F})$ be an
  orthonormal basis for the difference hull of $\mathcal{F}_i$ and
  choose any $u \in\mathcal{F}$ then
\begin{eqnarray*}
v - u & = & \mathcal{P}_{\mathcal{C}}(v) - u + (v -
\mathcal{P}_{\mathcal{C}}(v))\\
QQ^T(v - u) & = & QQ^T(\mathcal{P}_{\mathcal{C}}(v) - u) + QQ^T(v -
\mathcal{P}_{\mathcal{C}}(v))\\
QQ^T(v - u) & = & \mathcal{P}_{\mathcal{C}}(v) - u\\
\mathcal{P}_{\mathcal{C}}(v) & = & u + QQ^T(v - u)
\end{eqnarray*}
For the projection trajectory of the entire segment we therefore find
\begin{eqnarray*}
  p(\alpha) & = & u + QQ^T((x + \alpha d) - u)\\
  & = & (u + QQ^T(x - u)) + \alpha QQ^Td,
\end{eqnarray*}
which is linear in $\alpha$, as desired. Because the space is
partitioned, each finite breakpoint is simultaneously the end point of
the line segment in one region and the starting point of the line
segment in another region. Each value of $\alpha$ corresponds to one
point $x+\alpha d$, which has a unique projection, thereby showing the
continuity of the trajectory and completing the proof.
\end{proof}

\subsection{Properties of the weighted one-norm ball}

\subsubsection{Facial structure}

The weighted one-norm ball of radius $\tau$ is the convex hull of
vertices $\left\{\pm\tau/w_i\cdot e_i\right\}_i$. Every proper
$k$-face $\mathcal{F}$ of the weighted one-norm ball $\mathcal{C}_{w,1}$
can be written as the convex hull of $\{\sigma_i/w_i\cdot
e_i\}_{i\in\mathcal{I}}$, where $\mathcal{I}$ is a subset of
$\{1,\ldots,n\}$ with cardinality $k+1$, and $\sigma_i \in
\{-1,+1\}$. Throughout this section we assume that $\tau > 0$.

Given an $x\in\mathcal{C}$ we can determine $\mathcal{F}(x)$ as
follows. First, we need to check whether $\norm{x}_{w,1} < \tau$, in
which case $\mathcal{F}(x) = \mathcal{C}$.  Otherwise, $x$ lies on a
proper face, which can be uniquely characterized by the sign
vector $\sgn(x)$ whose $i$-th entry is given by
$\sgn(x_i)$. Determining $\mathcal{F}(x)$ and checking equality of
faces can therefore be done in $\mathcal{O}(n)$ time.

\subsubsection{Projection onto the feasible set}\label{Sec:LassoProjection}

Projection onto the weighted one-norm ball is discussed in
\cite{BER2010Fa} and is based on the solution of the prox function
\begin{equation}\label{Eq:ProxWeightedL1}
x_{\lambda}(u) = \prox_{\lambda\norm{\cdot}_{w,1}}(u) := \argmin_x\ \half\norm{x-u}_2^2
+ \lambda\norm{x}_{w,1} = \sign(u)\cdot\left[\abs{u} - \lambda w\right]_+,
\end{equation}
where $[\cdot]_+ = \max(0,\cdot)$, and the absolute value, sign
function, and multiplication in the right-hand side are evaluated
elementwise. Projection onto $\mathcal{C}_{w,1}$ then amounts to
finding the smallest $\lambda \geq 0$ for which
$\norm{x_{\lambda}(u)}_{w,1}\leq\tau$. The entries in
$x_{\lambda}(u)$, and therefore $\norm{x_{\lambda}(u)}_{w,1}$, are
continuous and piecewise linear in $\lambda$ with break points
occurring at $\lambda = \abs{u_i} / w_i$. We can obtain an
$\mathcal{O}(n\log n)$ algorithm that finds the optimal $\lambda$ and
subsequent projection by sorting the break points
\cite{BER2010Fa}. This can be reduced to an expected $\mathcal{O}(n)$
algorithm \cite{DUC2008SSCa} by avoiding the explicit sorting step.

\subsubsection{Self-projection cone of a face}\label{Sec:LassoSelfProjCone}

Given $x \in \mathcal{C}_{w,1}$ and search direction $d
\in\mathbb{R}^n$, we want to know if $d
\in\selfproj{\mathcal{F}(x)}$. When $\norm{x}_{w,1} < \tau$ it follows
that $x$ lies in the interior of $\mathcal{C}_{w,1}$ meaning that
$\mathcal{F}(x) = \mathcal{C}_{w,1}$ and $d \in
\selfproj{\mathcal{C}_{w,1}} = \mathbb{R}^n$, trivially.  For
$\norm{x}_{w,1} = \tau$, consider the support $\mathcal{I} = \{i \mid
x_i \neq 0\}$. Because the entries on the support are bounded away
from zero by definition and the soft-thresholding parameter $\lambda$
is initially linear in $\alpha$ it follows that the support of
$x(\alpha) = \mathcal{P}(x+\alpha D)$ includes $\mathcal{I}$ for all
sufficiently small $\alpha$. For $d$ to be in the self-projection cone
we therefore need to show that (1) $x+\alpha d$ does not move into the
polytope, and (2) that the support does not increase. It can be
verified that the first condition is satisfied if and only if
\begin{equation}\label{Eq:NondecreasingNorm}
\sum_{i\in\mathcal{I}} \sign(x_i)d_iw_i +
\sum_{i\not\in\mathcal{I}}\abs{d_i}w_i \geq 0.
\end{equation}
For the second condition to be satisfied we require for all
$i\not\in\mathcal{I}$ and sufficiently small $\alpha$ that the
absolute value of entry remains less than or equal to the threshold
value, namely $\alpha \abs{d_i} \leq w_i\lambda(\alpha) $. When the
support remains the same we find $\lambda(\alpha)$ by solving
\[
\sum_{i\in\mathcal{I}} w_i(\abs{x_i + \alpha d_i} -
w_i\lambda(\alpha)) =\tau,\quad\mbox{which gives}\quad
\lambda(\alpha) = \alpha\cdot\frac{ \sum_{i\in\mathcal{I}} 
w_i\sign(x_i)d_i}{\sum_{j\in\mathcal{I}}w_i^2},
\]
after writing $\abs{x_i + \alpha d_i} = \abs{x_i} + \sign(x_i)d_i$ and
recalling that $\norm{x}_{w,1} = \tau$.  A necessary (and sufficient)
condition for the support to remain the same is therefore that
\begin{equation}\label{Eq:IdenticalSupport}
\max_{i\not\in\mathcal{I}} \abs{d_i}/w_i \leq \frac{ \sum_{i\in\mathcal{I}} 
w_i\sign(x_i)d_i}{\sum_{i\in\mathcal{I}}w_i^2}.
\end{equation}
Summarizing the above we have:
\begin{btheorem}
  Given $x\in\mathcal{C}_{w,1}$ with support $\mathcal{I} = \{i \mid
  x_i \neq 0\}$, then $d \in\selfproj{\mathcal{F}(x)}$ if and only if
  $\norm{x}_{w,1}<\tau$, or $\norm{x}_{w,1} = \tau$ and
\[
\sum_{i\in\mathcal{I}} \sign(x_i)d_iw_i +
\sum_{i\not\in\mathcal{I}}\abs{d_i}w_i \geq 0,\quad\mathrm{and}\quad
\max_{i\not\in\mathcal{I}} \abs{d_i}/w_i \leq \frac{ \sum_{i\in\mathcal{I}} 
w_i\sign(x_i)d_i}{\sum_{i\in\mathcal{I}}w_i^2}.
\]
\end{btheorem}

\subsubsection{Orthogonal basis for a face}

For the construction of a quadratic approximation of objective
function $f$ over a face $\mathcal{F}$, we require an orthogonal basis
$\Phi$ for $\dhull{\mathcal{F}}$. For simplicity, consider the facet
of the unit cross polytope lying in the positive orthant in
$\mathbb{R}^3$. In other words, consider the unit simplex given by
$\mathrm{conv}\{e_1, e_2, e_3\}$. A first vector for the basis can
then be obtained by normalizing $e_2 - e_1$ to have unit norm. A
second vector orthogonal to the first can be obtained by connecting
the point halfway on the line segment $e_1$--$e_2$ to $e_3$, that is,
$e_3 - (e_1+e_2)/2$, followed again by normalization. This can be
generalized, and for a general $k$-simplex we find $(k+1)\times k$
basis $Q = [q_1,q_2,\ldots,q_k]$ with
\[
q_j = (e_{j+1} - \frac{1}{j}\sum_{i=1}^{j}e_i) / \sqrt{1 + 1/j}.
\]
In other words
\[
Q_{i,j} = \begin{cases}
-\sqrt{1 / (j^2 + j)} & i \leq j \\
\sqrt{j / (j+1)} & i = j+1 \\
0 & \mbox{otherwise}.
\end{cases}
\]
It can be seen that the above procedure amounts to taking a QR
factorization of the $k+1\times k$ matrix $[-e,\ I]^T$ of differences
between the first vertex and all others, and discarding the last
column in $Q$, whose entries are all equal to $1/\sqrt{n}$. The
special structure of $Q$ allows us to evaluate matrix-vector products
with $Q$ itself and its transpose in $\mathcal{O}(k)$ time, without
having to form the matrix explicitly.  For the general case, let
$\mathcal{F} = \mathcal{F}(x)$. For the case where $\mathcal{F} =
\mathcal{C}$ no projection is needed and we can simply choose $\Phi =
I$. Otherwise, let $\mathcal{I} = \{i \mid x_i \neq 0\}$ denote the
support of $x$. The desired basis can then be obtained by first
restricting the vector to its support and then normalizing the sign
pattern, thus giving:
\[
\Phi =
I_{\mathcal{I}}\cdot\mathrm{diag}(\mathrm{sgn}(x_{\mathcal{I}}))\cdot
Q.
\]
Matrix-vector products with $\Phi$ can be evaluated in
$\mathcal{O}(n)$ time, again without forming the matrix

\paragraph{Generic weighted one-norm ball}

\begin{figure}[t!]
\centering
\begin{tabular}{cccc}
\includegraphics[height=4.35cm]{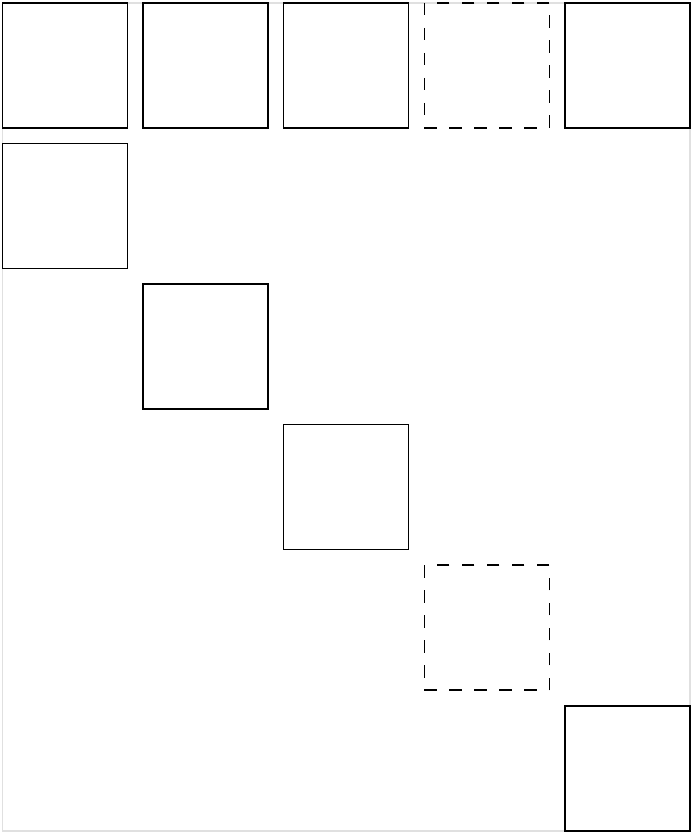}%
\begin{picture}(0,0)(0,0)
\put(-98,112){$v_1$}
\put(-77,112){$v_1$}
\put(-56.2,112){$v_1$}
\put(-37,111.2){$\cdots$}
\put(-14.3,112){$v_1$}
\put(-99,90.5){$w_1$}
\put(-78,70){$w_2$}
\put(-57.2,49.2){$w_3$}
\put(-37,26.1){$\ddots$}
\put(-15.3,8){$w_n$}
\end{picture}&
\includegraphics[height=4.35cm]{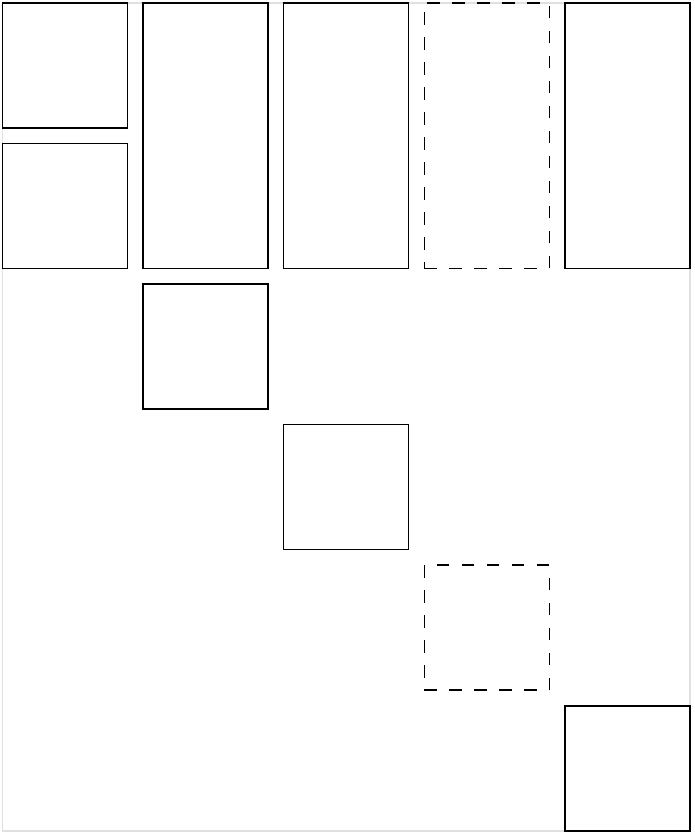}%
\begin{picture}(0,0)(0,0)
\put(-98,112){$v_1$}
\put(-77,101){$v_2$}
\put(-56.2,101){$v_2$}
\put(-37,100.2){$\cdots$}
\put(-14.3,101){$v_2$}
\put(-99,90.5){$w_1$}
\put(-78,70){$w_2$}
\put(-57.2,49.2){$w_3$}
\put(-37,26.1){$\ddots$}
\put(-15.3,8){$w_n$}
\end{picture}&
\includegraphics[height=4.35cm]{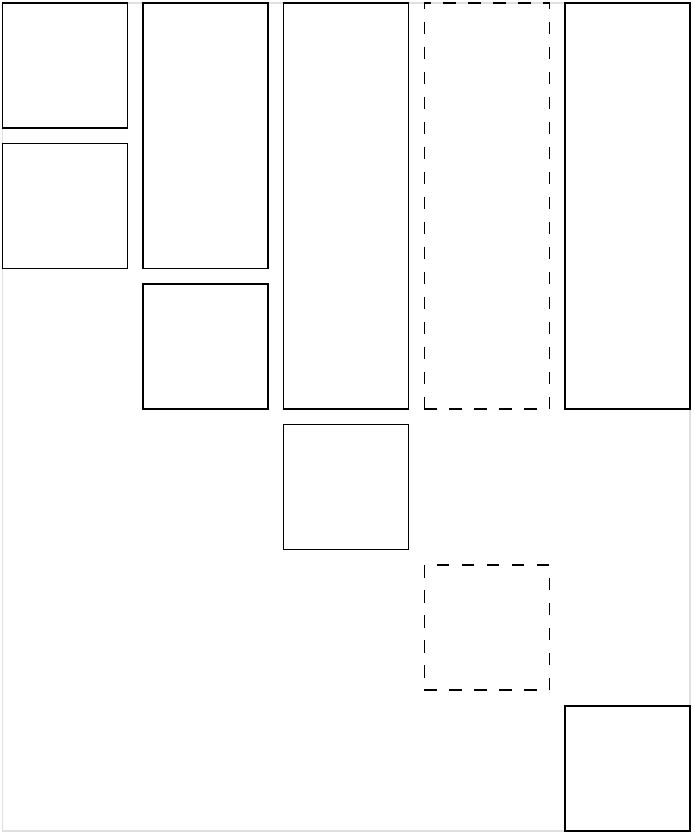}%
\begin{picture}(0,0)(0,0)
\put(-98,112){$v_1$}
\put(-77,101){$v_2$}
\put(-56.2,90.5){$v_3$}
\put(-37,89.7){$\cdots$}
\put(-14.3,90.5){$v_3$}
\put(-99,90.5){$w_1$}
\put(-78,70){$w_2$}
\put(-57.2,49.2){$w_3$}
\put(-37,26.1){$\ddots$}
\put(-15.3,8){$w_n$}
\end{picture}&
\includegraphics[height=4.35cm]{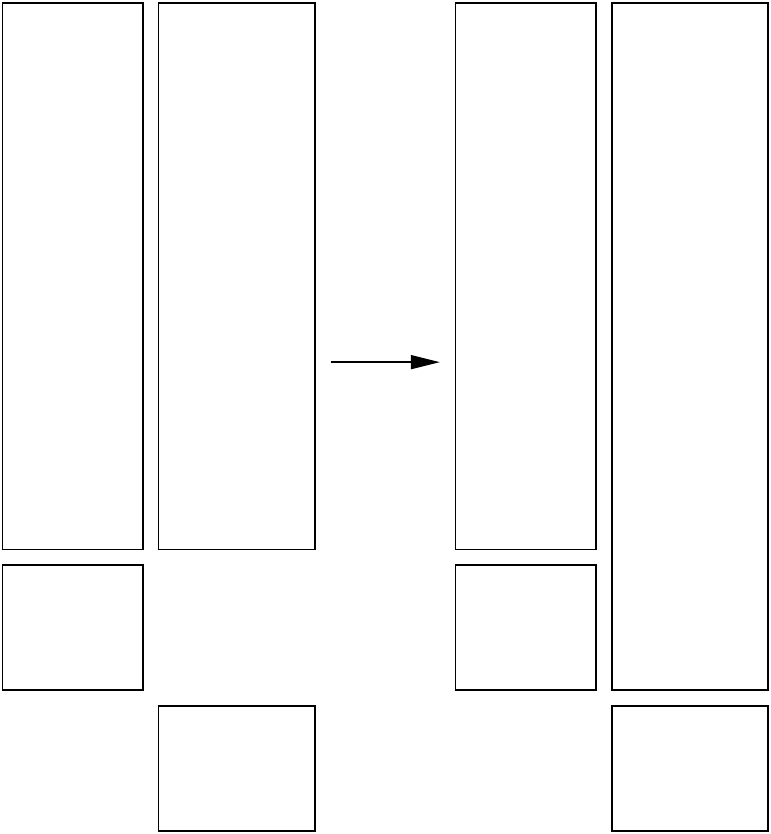}%
\begin{picture}(0,0)(0,0)
\put(-108.5,79){$v_k$}
\put(-84,79){$v_k$}
\put(-41,79){$v_k$}
\put(-22,70){$v_{k+1}$}
\put(-109.7,28.5){$w_k$}
\put(-90.5,8){$w_{k+1}$}
\put(-42.3,28.5){$w_k$}
\put(-23.2,8){$w_{k+1}$}
\end{picture}\\
(a) & (b) & (c) & (d)
\end{tabular}
\caption{Stages of the orthogonalization process.}\label{Fig:Sweep}
\end{figure}

For the weighted one-norm ball we consider a face given by
$\mathrm{conv}(w_0e_1,w_1e_2,\ldots,w_{n}e_{n})$, with nonzero weights
$w_0$ to $w_{n}$. (Throughout this paragraph it is more convenient to
work with a weight vector whose elements are the inverse of the
weights appearing in the weighted one norm; the actual vertices of the
weighted one norm ball are $\pm w_i^{-1}e_i$, not $\pm w_i e_i$.)  We
would again like to obtain an orthonormal basis corresponding to the
face. This can be done by applying QR factorization to the matrix of
differences between the vertices, as illustrated in
Figure~\ref{Fig:Sweep}(a) with $v_1 = -w_0$. The two operations in
this process are projecting out the contributions of all subsequent
columns and normalizing the columns to unit norm. We do not normalize
until the very end but do keep track of the squared two norm of the
completed columns. Given vectors $a$ and $b$ we obtain the component
in $b$ orthogonal to $a$ by evaluating $b - \frac{\langle
  a,b\rangle}{\langle a,a\rangle} a$. In the first step of the
factorization (we are interested only in $Q$) we orthogonalize with
respect to the first column $a$. The inner product of each column with
$a$ is identical and equal to $\alpha_1 = \langle v_1,v_1\rangle =
\norm{v_1}_2^2$. Using this we also compute the squared two norm of
the first column as $\gamma_1 = \alpha_1 + w_1^2$. After the sweep
with the first column we are left with the matrix shown in
Figure~\ref{Fig:Sweep}(b) where
\[
v_2 = \left[\begin{array}{r} v_1 - \frac{\alpha_1}{\gamma_1}v_1\\
-\frac{\alpha_1}{\gamma_1}w_1\end{array}\right]
= \left[\begin{array}{r}\frac{\gamma_1 - \alpha_1}{\gamma_1}v_1\\
    -\frac{\alpha_1}{\gamma_1}w_1\end{array}\right]
= \left[\begin{array}{r}\frac{w_1^2}{\gamma_1}v_1\\ -\frac{\alpha_1}{\gamma_1}w_1\end{array}\right].
\]


The next step is to sweep with the updated second column. For this we
compute the inner product with the remaining columns and itself,
yielding $\alpha_2 = \norm{v_2}_2^2$ and $\gamma_2 = \alpha_2 +
w_2^2$, respectively.  After this sweep we arrive at the matrix given
in Figure~\ref{Fig:Sweep}(c). Proceeding like this we successively
sweep with each of the column until we are done. Consider now the
sweep with some column $k$, illustrated in
Figure~\ref{Fig:Sweep}(d). Given $\alpha_k = \norm{v_k}_2^2$ and
$\gamma_k = \alpha_k + w_k^2$ we find
\[
v_{k+1} = \left[\begin{array}{r}v_k -
    \frac{\alpha_k}{\gamma_k}v_k\\-\frac{\alpha_k}{\gamma_k}w_k\end{array}\right]
=
\left[\begin{array}{r}\frac{w_k^2}{\gamma_k}v_k\\-\frac{\alpha_k}{\gamma_k}w_k\end{array}\right],
\]
from which we derive recurrence relations
\begin{equation}\label{Eq:RecurrenceAlpha}
\alpha_{k+1} = \left(\frac{w_k^2}{\gamma_k}\right)^2\norm{v_k}_2^2 +
\left(\frac{\alpha_k}{\gamma_k}\right)^2w_k^2 =
\frac{w_k^2w_k^2}{\gamma_k^2}\alpha_k +
\frac{\alpha_k^2}{\gamma_k^2}w_k^2 = \alpha_kw_k^2  \frac{ w_k^2 +
  \alpha_k}{\gamma_k^2} = \frac{\alpha_k w_k^2}{\gamma_k},
\end{equation}
and $\gamma_{k+1} = \alpha_{k+1} + w_{k+1}^2$. With $\alpha_1$ and
$\gamma_1$ as given above, this allows us to compute all $\alpha$ and
$\gamma$ values. Ultimately we are interested in the final orthonormal
$Q$ matrix. Defining scaling factors
\begin{equation}\label{Eq:Mu}
\mu_{i,j} = \prod_{k=i}^{j} \frac{w_k^2}{\gamma_k}\ \  \mbox{for}\ \ 
1 \leq i \leq j \leq n,
\end{equation}
 as well as factors $u_i = -\alpha_i / \gamma_i$ for $1 \leq i \leq n$
 and  $u_0 := -1$, it can be found based on the structure of the $v$
vectors that
\[
Q = \diag(w)\cdot\left[\begin{array}{cccccc}
u_0 & \mu_{1,1}u_0 & \mu_{1,2}u_0 & \mu_{1,3}u_0 & \cdots & \mu_{1,n-1}u_0\\
1 & u_1 & \mu_{2,2}u_1 & \mu_{2,3}u_1 &\cdots & \mu_{2,n-1}u_1\\
       & 1 & u_2 & \mu_{3,3}u_2 & \cdots & \mu_{3,n-1}u_2 \\
       &        & \ddots & \ddots& \ddots & \vdots \\
       &         &            & 1 & u_{n-2} & \mu_{n-1,n-1}u_{n-2}\\
       &         &            &             & 1 & u_{n-1}\\
   & & & & & 1
\end{array}\right]\cdot\diag(1/\sqrt{\gamma}).
\]
Multiplication with this matrix and its transpose may still seem
expensive but we now show how the structure enables $\mathcal{O}(n)$
algorithms for both operations. Multiplication with the diagonal
matrices is trivial so we focus only on multiplication with the
central part of the matrix. Looking at a small example we can
decompose this matrix as
\[
\left[\begin{array}{ccc}
u_0 & \mu_{1,1} u_0 & \mu_{1,2}u_0 \\
1    & u_1 & \mu_{2,2}u_1\\
& 1 & u_2\\
& & 1
\end{array}\right]
= \left[\begin{array}{c} 0 \\ I\end{array}\right] +
\diag(u)\left[\begin{array}{ccc}
1 & \mu_{1,1} & \mu_{1,2} \\
& 1 & \mu_{2,2} \\
&&1 \\
\\
\end{array}\right] = \left[\begin{array}{c} 0 \\ I\end{array}\right] +
\diag(u)\left[\begin{array}{c}M\\ 0\end{array}\right]
\]
The key part is multiplication with the last matrix $M$. To evaluate
$y = Mv$ we initialize $y_3 = v_3$ and then work upwards. Direct
evaluation gives $y_2 = v_2 + \mu_{2,2}v_3$, which can be rewritten as
$y_2 = v_2 + \mu_{2,2}y_2$. A pattern emerges when looking at the
computation of $y_1$:
\begin{eqnarray*}
y_1 &= & v_1 + \mu_{1,1} v_2 + \mu_{1,2}v_3 \\
& = & v_1 + \mu_{1,1}(v_2 + \mu_{2,2}v_3) \\
& = & v_1 + \mu_{1,1}y_2,
\end{eqnarray*}
where $\mu_{1,2} = \mu_{1,1}\mu_{2,2}$, or more generally $\mu_{i,k} =
\mu_{i,j}\mu_{j+1,k}$ for $i \leq j\leq k$, follows from the
definition of $\mu$ in \eqref{Eq:Mu}. Given $y_{n} = v_{n}$, we
therefore obtain the recurrence $y_k = v_k + \mu_{k,k}y_{k+1}$,
which allows us to evaluate $y = Mv$ in linear time. With $v$
appropriately redefined we now look at $y = M^Tv$:
\[
\left[\begin{array}{c}y_1\\y_2 \\ y_3\end{array}\right] =
\left[\begin{array}{ccc}
1 && \\ \mu_{1,1} & 1 & \\ \mu_{1,2} & \mu_{2,2} & 1 \end{array}\right]
\left[\begin{array}{c}v_1\\v_2 \\ v_3\end{array}\right].
\]
Starting with $y_1 = v_1$ we find $y_2 = \mu_{1,1} y_1 + v_2$ and $y_3
= \mu_{2,2}y_2 + v_3$, using $\mu_{1,2} = \mu_{1,1}\mu_{2,2}$. This
gives the recurrence $y_{k+1} = \mu_{k,k}y_k + v_{k+1}$. We summarize
the initialization and multiplication with $Q$ and $Q^T$ in
Algorithms~\ref{Alg:WeightedL1QInit}--\ref{Alg:WeightedL1QTv}. Note
that these algorithms use a different indexing scheme for a convenient
implementation. For practical implementations we can precompute and
store $1 / \sqrt{\gamma_k}$ instead of $\gamma_k$ and avoid storing
$\alpha$ since it is not used during the evaluation of matrix-vector
products. Alternatively, we can reduce the memory footprint at the
cost of increased computation by storing only $\alpha$ and
re-computing $\mu_k$, $u_k$, and $\gamma_k$ whenever they are needed.

\begin{figure*}[t!]
\centering
\begin{minipage}[t]{0.56\textwidth}
\vspace*{0pt}
\begin{algorithm}[H]
\dontprintsemicolon
\KwData{Weight vector $w = [w_1,\ldots,w_n]$}
\KwResult{Vectors $\alpha$, $\gamma$, $u$, and $\mu$}
Initialize $\alpha_1 = w_1^2$, $u_1 = -1$\;
\For{\rm{$k = 1$ to $n-1$}}{
$\gamma_{k} = \alpha_k + w_{k+1}^2$\;
$\mu_{k} = w_{k+1}^2 / \gamma_k$\;
$u_{k+1} = -\alpha_{k} / \gamma_{k}$\;
$\alpha_{k+1} = \alpha_k \mu_{k}$\;
} 
\caption{Initialization for multiplication with the orthogonal basis for
  a face of the weighted one-norm ball.}\label{Alg:WeightedL1QInit}
\end{algorithm}
\end{minipage}

\vspace*{10pt}
\begin{minipage}[t]{0.45\textwidth}
\vspace*{0pt}
\begin{algorithm}[H]
\dontprintsemicolon
\KwData{Vectors $\gamma$, $u$, $w$, $\mu$ and $v \in\mathbb{R}^{n-1}$}
\KwResult{Vector $y = Qv$}
Initialize $s = v_{n-1} / \sqrt{\gamma_{n-1}}$, $t=0$\;
$y_{n} = w_{n} s$\;
\For{\rm{$k = n\!-\!1$ down to $2$}}{
   $t \gets \mu_{k} t + s$\;
   $s \gets v_{k-1} / \sqrt{\gamma_{k-1}}$\;
   $y_{k} = w_{k} (u_k t + s)$\;
} 
$y_1 = w_1 u_1 (\mu_1 t + s)$\;
\caption{Multiplication with orthogonal basis $Q$ for a face of the
  weighted one-norm ball: $y=Qv$.}\label{Alg:WeightedL1Qv}
\end{algorithm}
\end{minipage}
\ \  
\begin{minipage}[t]{0.45\textwidth}
\vspace*{0pt}
\begin{algorithm}[H]
\dontprintsemicolon
\KwData{Vectors $\gamma$, $u$, $w$, $\mu$ and $v \in\mathbb{R}^{n}$}
\KwResult{Vector $y = Q^Tv$}
$t = w_1 v_1 u_1$\;
$s = w_2 v_2$\;
$y_1 = (t + s) / \sqrt{\gamma_1}$\;
\For{\rm{$k = 2$ to $n\!-\!1$}}{
   $t \gets \mu_{k-1}t + u_k s$\;
   $s \gets w_{k+1}v_{k+1}$\;
   $y_{k} = (t+s) / \sqrt{\gamma_{k}}$\;
} 
\caption{Multiplication with orthogonal basis $Q$ for
  a face of the weighted one-norm ball: $y=Q^Tv$.}\label{Alg:WeightedL1QTv}
\end{algorithm}
\end{minipage}
\end{figure*}

\subsubsection{Generation of the projection arc}\label{Sec:ProjectionArc}

\begin{figure}
\centering
\begin{tabular}{c}
\includegraphics[width=0.45\textwidth]{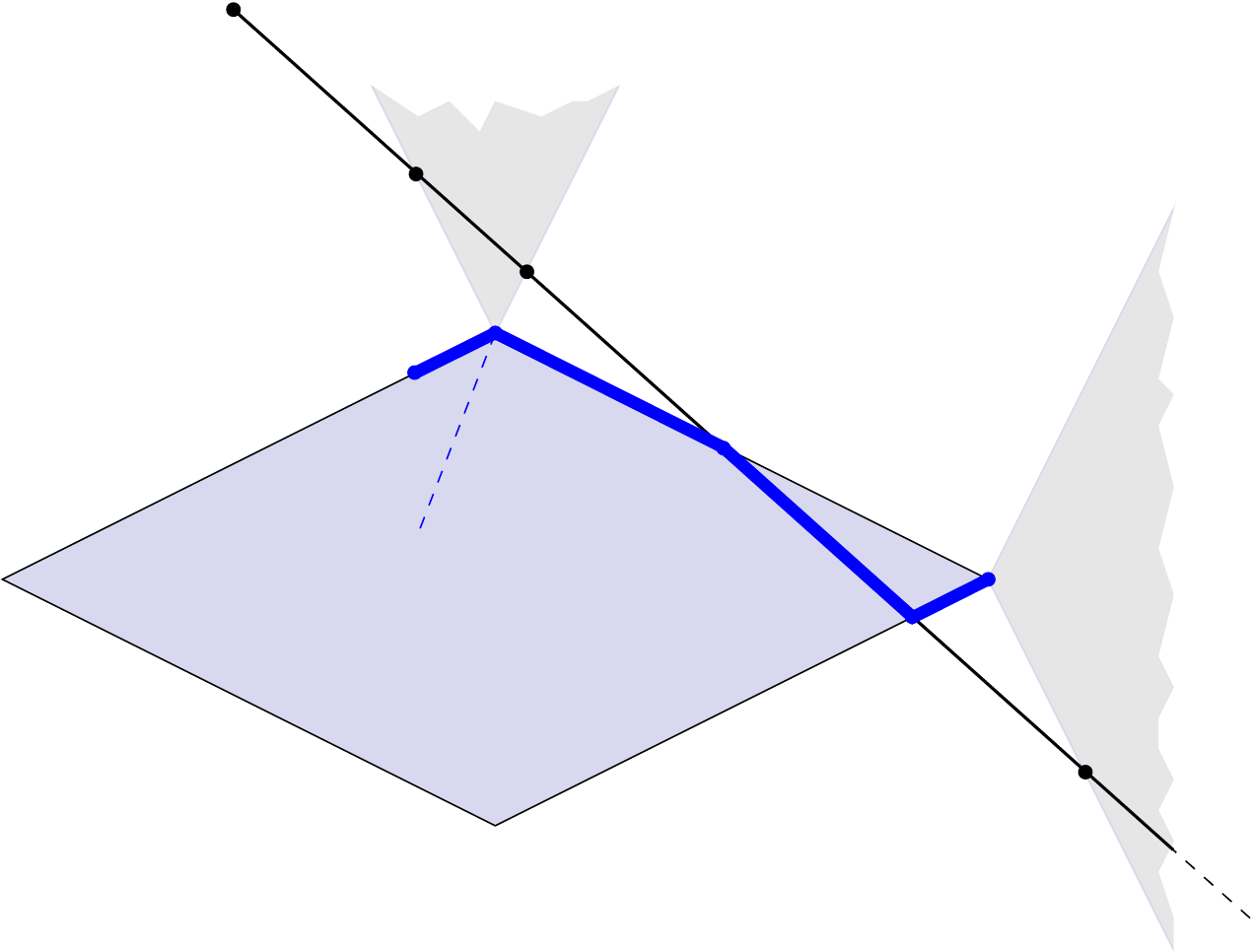}%
\begin{picture}(0,0)(0,0)
\put(-178,150){$s$}
\put(-136,130){$x_{(1)}$}
\put(-117,114){$x_{(2)}$}
\put(-85.5,85){$x_{(3)}$}
\put(-65.5,45.25){$x_{(4)}$}
\put(-45.5,25.25){$x_{(5)}$}
\put(-144,87){\color{blue}$p$}
\put(-156,67){\color{blue}$p_{(1)},p_{(2)}$}
\put(-104.5,79){\color{blue}$p_{(3)}$}
\put(-78.5,56.5){\color{blue}$p_{(4)}$}
\put(-38.5,60){\color{blue}$p_{(5)}$}
\end{picture}
\end{tabular}
\caption{Illustration of half line $x(\alpha) = s+\alpha d$ (black)
  with its projection arc (thick blue) on a weighted one-norm
  ball. The break points on the original line and their corresponding
  projection are indicated by $x_{(i)}$ and $p_{(i)}$,
  respectively. At the intersections with the domain we have $p_{(3)}
  = x_{(3)}$ and $p_{(4)} = x_{(4)}$. The light gray wedges indicate
  the normal cones at two of the
  vertices.}\label{Fig:WeightedL1ProjArc}
\end{figure}

In this section we consider the computation of the projection arc
$p(\alpha) = \mathcal{P}(x(\alpha))$ of half line $x(\alpha) =
s+\alpha d$ with $\alpha \geq 0$. We allow any starting point $s \in
\mathbb{R}^n$, even though $s \in \mathcal{C}_{w,1}$ always holds for
our application. From Section~\ref{Sec:PiecewiseProjectionArc} we know
that the projection arc is piecewise linear with discrete break points
at $\alpha = \alpha_i$. As illustrated in
Figure~\ref{Fig:WeightedL1ProjArc} there are three types of break
points: (1) the support reduces and we move to a lower dimensional
face ($\alpha_1$ and $\alpha_5$); (2) the support increases and we
move to a higher dimensional face ($\alpha_2$); and (3) we intersect
the polytope boundary ($\alpha_3$ and $\alpha_4$). The algorithm for
computing the projection arc thus proceeds as follows. Starting with
$\alpha_0 = 0$ we compute at each iteration $i \geq 0$ the minimal
$\alpha_{i+1} > \alpha_{i}$ for which one of the events
occurs. Once this is done we compute $x_{(i)}$ and, based on the type
of event, determine the corresponding $p_{(i)}$. The algorithm
completes when none of the three events happens for $\alpha$ exceeding
the current value. We now show how the next $\alpha$ for each event
type is computed.

\paragraph{Intersection with the boundary.} In order to determine
determine intersections with the boundary we need to keep track of
$\kappa(\alpha) = \norm{x(\alpha)}_{w,1}$. This function is linear
with break points occurring whenever one of the elements in $x$
crosses zero. Let $\mathcal{K} = \{j \mid s_j d_j < 0\}$ denote the
set of indices that will at some point cross zero, and define $r_j =
\abs{d_j}$ when $j\in\mathcal{K}$ and $r_j = -\abs{d_j}$ otherwise. We
maintain the directional derivative of $\kappa(\alpha)$ with respect
to increasing $\alpha$ as
\begin{equation}\label{Eq:ComputationSlope}
\rho = \sum_{j\not\in\mathcal{K}} w_j\abs{d_j} - \sum_{j \in\mathcal{K}}
w_j\abs{d_j} = \sum_{j}r_jw_j.
\end{equation}
Whenever some $x_j$ reaches zero we remove $j$ from the set
$\mathcal{K}$, set $r_j$ to $-r_j$, and update $\rho$ to the value
$\rho + 2w_j\abs{d_j}$. To deal with these updates we add a fourth type
of event corresponding to zero crossings. These happen at $\alpha$
values $-s_j/d_j$ for $j\in\mathcal{K}$, which that can be
pre-computed and sorted at the beginning of the algorithm. At the
beginning of iteration $i$ we are given the current slope $\rho$ and
have
\[
\kappa(\alpha_{i} + \delta) = \kappa(\alpha_{i}) + \delta \rho,
\]
for limited $\delta \geq 0$. Solving for the next boundary
intersection gives $\delta = + (\tau - \tau(\alpha_{i})) / \rho$ and
$\alpha = \alpha_i + \delta$. Whenever a zero crossing happens before
this value of $\alpha$ the algorithm will encounter that event first,
update $\rho$ accordingly, and recompute $\delta$ in the next
iteration. If $\alpha \leq \alpha_i$ we omit the boundary crossing
event from consideration for the current iteration. Regardless of the
event type we update $\kappa(\alpha_{i+1}) = \kappa(\alpha_i) +
(\alpha_{i+1}-\alpha_i)\rho$.  In the case of a boundary crossing
(there can at most be two such events) we set $\lambda = 0$ and
$\kappa = \tau$.

\paragraph{Changes in the support.}

Given $x = x_{(i)}$ with norm $\kappa = \norm{x}_{w,1}$ and
projection $p=\mathcal{P}(x)$ we want to find the smallest $\delta$
such that the support of $\mathcal{P}(x(\alpha_i+\delta))$ differs from
$\mathcal{I} := \{i \mid p_i \neq 0\}$. We assume without loss of
generality that either $\kappa > \tau$, or $\kappa = \tau$ with
corresponding rate of change $\rho > 0$, otherwise we either move
inside the polytope or along one of the faces, in which case the next
event is guaranteed to either be a zero crossing or a boundary
intersection. Many of the equations in this section, such as $\abs{x_i
  + \delta d_i} = \abs{x_i} + \delta r_i$, are valid only for
sufficiently small $\delta$ and break down whenever there is a sign
change or an intersection with the boundary of $\mathcal{C}_{w,1}$. We
nevertheless use all equations as if they hold for all $\delta$. The
rationale for this is that we are interested in the first event. If
the first event is a change in support then clearly there were no sign
changes or boundary intersections before that point, which means that
all equations used to compute $\delta$ were valid. Otherwise we do
have a sign change or boundary intersection, in which case we ignore
the incorrect $\delta$ value for the support change.

We now consider the point $x' = x(\alpha_i + \delta)$. For the support
of the corresponding projection to remain the same we must have that
the threshold parameter $\lambda'(\delta)$ computed based on the index
set $\mathcal{I}$ be consistent with $x'$. From the projection
operator it follows that
\begin{eqnarray*}
\tau & = & \sum_{i\in\mathcal{I}}w_i(\abs{x_i + \delta d_i} - \lambda'(\delta) w_i)
 =  \sum_{i\in\mathcal{I}}w_i(\abs{x_i} + \delta r_i - \lambda'(\delta) w_i) \\
& = & \kappa + \delta \sum_{i\in\mathcal{I}}w_i r_i - \lambda'(\delta) \sum_{i\in\mathcal{I}}w_i^2,
\end{eqnarray*}
from which we find
\begin{equation}\label{Eq:LambdaPrime}
\lambda'(\delta) = \frac{(\kappa -\tau) +
  \delta\sum_{i\in\mathcal{I}}w_ir_i}{\sum_{i\in\mathcal{I}}w_i^2}
= \lambda +
\delta\frac{\sum_{i\in\mathcal{I}}w_i r_i}{\sum_{i\in\mathcal{I}}w_i^2}
= \lambda + \delta\mu,
\end{equation}
where $\lambda$ is the threshold parameter for $x$ and $\mu$ is the
directional derivative with respect to $\delta$. The resulting
threshold parameter $\lambda'(\delta)$ is then consistent with $x'$ if
only if $\{i \mid \abs{x_i'} > w_i\lambda'(\delta)\} =\mathcal{I}$.

\paragraph{Additions to the support.}

\begin{figure*}[t]
\centering
\begin{minipage}[t]{0.75\textwidth}
\vspace*{0pt}
\begin{algorithm}[H]
\dontprintsemicolon
\KwData{Sets $\mathcal{I}$ and $\mathcal{J}$}
\KwResult{Set $\mathcal{I}'$}
Set $a = \sum_{i\in\mathcal{I}}w_ir_i$, $b =
\sum_{i\in\mathcal{I}}w_i^2$\;
Set $a_j = w_jr_j$, $b_j = w_jr_j$ for all $j\in\mathcal{J}$\;
Initialize $\mathcal{I}' = \mathcal{I}$\;
\While{\rm{($\mathcal{J} \neq \emptyset$)}}{
Find $j$ such that $a_j/p_j \geq a_k/p_k$ for all $k\in\mathcal{J}$\;
Set $a \gets a + a_j$, $b\gets b + b_j$\;
Set $\mathcal{J} = \{ j \in\mathcal{J} \setminus \{j\} \mid a_j/b_j > a/b\}$\;
} 
\caption{Selection of additions to the support given candidate set $\mathcal{J}$.}\label{Alg:SupportAddition}
\end{algorithm}
\end{minipage}
\end{figure*}

When considering additions to the support we assume, in addition to
zero crossings and boundary intersections, that there are no events
corresponding to variables leaving the support. A necessary condition
for a variable $j \not\in\mathcal{I}$ to enter the support is that
$\abs{x_j'} > w_j\lambda'(\delta)$, or equivalently $\abs{x_j} +
\delta r_j > w_j(\lambda + \delta \mu)$. From the definition of the
support we have $\abs{x_j} \leq w_j\lambda$, and for a variable to
enter it must therefore hold that $r_j> w_j\mu$, or equivalently
$r_j/w_j>\mu$. In this case, the value of $\delta$ at which variable
$j$ is about the enter is given by $\delta_j = (\abs{x_i} - \lambda
w_i) / (\mu w_i - r_i)\geq 0$, otherwise we set $\delta_j =
+\infty$. The smallest $\delta$ for which an addition to the support
is about to happen, if any, is then given by $\delta =
\min_{j\not\in\mathcal{I}}\delta_j$ with the set of variables staged
to enter given by $\mathcal{J} = \{j\not\in\mathcal{I} \mid \delta_i =
\delta\}$. Provided that no event occurs before this point, at least
one of the variables in $\mathcal{J}$ will enter. As it does, it
changes the rate of change in $\lambda$, which may mean that some of
the other variables in $\mathcal{J}$ never actually enter the
support. As such, care needs to be taken in determining which
variables enter and which do not; otherwise the same variable may
repeatedly enter and leave the support, causing the algorithm to cycle
forever. The following stage of the algorithm iteratively constructs
the desired new support set $\mathcal{I}'$, which is the largest
subset of $\mathcal{I} \cup \mathcal{J}$ such that $j \in\mathcal{J}$
is in $\mathcal{I}'$ if and only if $r_j > w_j\mu'$ with
\begin{equation}\label{Eq:MuPrime}
\mu' =
\frac{\sum_{i\in\mathcal{I}'}w_ir_i}{\sum_{j\in\mathcal{I}'}w_i^2}.
\end{equation}
The new set $\mathcal{I}'$ is formed iteratively starting from
$\mathcal{I}$ by successively adding more elements from $\mathcal{J}$
until it satisfies \eqref{Eq:MuPrime}.  Starting with $\mathcal{I}' = \mathcal{I}$ we define
\[
a = \sum_{i\in\mathcal{I}}w_ir_i,\ 
b = \sum_{i\in\mathcal{I}}w_i^2,\quad\mbox{and}\quad
a_j = w_jr_j,\ b_j = w_j^2\ \mbox{for $j\in\mathcal{J}$}.
\]
This allows us to rewrite $r_j > w_j\mu$ as $a_j/b_j > \mu =
a/b$. Given any $p_1 / q_1 \prec p_2/q_2$ with $q_1,q_2 > 0$ and
relational operator $\prec \in \{=,<,\leq\}$ it can be verified that
\begin{equation}\label{Eq:PQ}
\frac{p_1}{q_1} \prec \frac{p_1+p_2}{q_1+q_2} \prec \frac{p_2}{q_2}.
\end{equation}
Let $j \in\mathcal{J}$ be such that $a_j/b_j \geq a_k/b_k$ for all
$k\in\mathcal{J}$. We show that $j$ must be member of
$\mathcal{I}'$. Assume by contradiction that
$j\not\in\mathcal{I}'$. By repeated application of \eqref{Eq:PQ} we
find that
\[
\frac{a}{b}
< \frac{a + \sum_{i \in \mathcal{I}'\setminus\mathcal{I}}a_i}{b +
  \sum_{i\in\mathcal{I}'\setminus\mathcal{I}}b_i}
= \mu'
< \frac{\sum_{i \in \mathcal{I}'\setminus\mathcal{I}}a_i}{
  \sum_{i\in\mathcal{I}'\setminus\mathcal{I}}b_i}
\leq \frac{a_j}{b_j}.
\]
This shows that $r_j/w_j > \mu'$, which contradicts
$j\not\in\mathcal{I}'$. Given that $j$ must be part of $\mathcal{I}'$
we can add $j$ to $\mathcal{I}'$ and update $a$ and $b$ to $a+a_j$ and
$b+b_j$, respectively. After this we remove element $j$ from
$\mathcal{J}$ as well as all elements $k$ for which $a_k/b_k \leq
a/b$. This process is iterated until $\mathcal{J} = \emptyset$, at
which point we have the final $\mathcal{I}'$. This algorithm is
summarized in Algorithm~\ref{Alg:SupportAddition}.

\paragraph{Removal from the support.}

For the removal of items from the support, we start by finding the
smallest $\delta$ such that $\abs{x_i'} = w_i\lambda'(\delta)$ for
some $i\in\mathcal{I}$, if any. Next, we define the set $\mathcal{J}$
of entries staged to leave the support as the set of all all entries
$i\in\mathcal{I}$ for which the above equality holds for the given
$\delta$. Removal from the support causes the rate of change in
$\lambda$ to increase (see Section~\ref{AppSec:Projection} for more
details). This means that all entries staged in $\mathcal{J}$ leave
the support, and therefore that no further filtering is needed.

\begin{figure*}
\centering
\begin{minipage}[t]{0.75\textwidth}
\vspace*{0pt}
\begin{algorithm}[H]
\dontprintsemicolon
\KwData{Vectors $s$, $d$}
\KwResult{Break points of the projection arc}
Set $\mathcal{K} = \{j\mid x_jd_j < 0\}$ and compute $r$ according to \eqref{Eq:ComputationSlope}\;
Initialize $\alpha_0 = 0$, $i = 0$, $\kappa_0 = \norm{s}_{w,1}$\;
\While{\rm{true}}{
Determine the next event $\alpha_{i+1}$\;
Return if there is no next event\;
\hspace*{5pt}\;
Update $\kappa_{i+1} = \kappa_i + (\alpha_{i+1} - \alpha_i)r$\;
If zero crossing of index $j$: update $r\gets r + 2\abs{d_j}$\;
If boundary crossing: set $\kappa_{i+1} \gets \tau$, $\lambda_{i+1} = 0$.
} 
\caption{Outline of the general algorithm for computing all break points of the projection arc.}
\end{algorithm}
\end{minipage}
\end{figure*}

\subsubsection{Line search along the projection arc}

Given the set of break points $\alpha_k$ and the corresponding changes
to the support we can perform a line search along the piecewise
linear projection arc. When restricting all relevant vectors to the
support for a given segment, we can write $p(\alpha) = s + \alpha d -
\lambda(\alpha) v$, where $v$ is a vector of sign values for
soft-thresholding. For the objective value we need $Ap(\alpha)$ for
$\alpha_k \leq \alpha\leq \alpha_{k+1}$, which can then be written as
$Ax(\alpha) = As + \alpha Ad - \lambda(\alpha) Av$. Once we have the
three matrix-vector products with $A$, the objective function can
easily be converted into a quadratic function in $\alpha$ by
evaluating the appropriate inner products of these vectors. As the
support or signs change we need to update the products $As$, $Ad$, and
$Av$ by adding or subtracting multiples of the required columns of
$A$. Assuming that $A$ is explicitly available or that the columns of
$A$ can be extracted in $\mathcal{O}(m)$ time, each update takes
$\mathcal{O}(m)$ time.  The objective function over each segment is
quadratic and the minimum within the segment is therefore easily
determined. For the evaluation of the overall computational complexity
we need to know the maximum number of segments that a projection arc
can have, or likewise, the maximum number of faces that a line can
project onto. In Appendix~\ref{AppSec:Projection} we prove the
following result:
\begin{btheorem}\label{Thm:MaxProjectionFaces}
  The projection of the line $x + \alpha d$ onto a (weighted) one-norm
  ball in $\mathbb{R}^n$ is piecewise linear with at
  most $4n-1$ segments. For every $n \geq 1$ there exist parameters
  $x$, $d$, and weights $w$ for which this bound is achieved.
\end{btheorem}
Combined with this maximum number of possible segments possible, it
follows that the line search can be done in $\mathcal{O}(mn)$ time.
(In practice it may be necessary to recompute entirely the three
matrix-vector products at regular intervals to avoid numerical issues.)
For weighted one-norm balls the procedure requires slightly more
bookkeeping but is otherwise the same.

\subsubsection{Maximum step length along a face}

Given a feasible search direction $d$ it is useful to know the maximum
$\alpha$ for which $x+\alpha d \in \mathcal{C}_{w,1}$. When $x$ lies
in the interior of $\mathcal{C}_{w,1}$ or when
\eqref{Eq:NondecreasingNorm} is violated and $x+\alpha d$ moves into
the interior, we need to compute the first intersection with the
boundary. The procedure for doing this was described earlier in
Section~\ref{Sec:ProjectionArc}. When $x$ lies on a proper face of
$\mathcal{C}_{w,1}$ and $d$ moves along the face or onto a higher
dimensional face, the maximum step length is determined by the first
element to reach zero:
\[
\alpha_{\max} = \min_{i:x_id_i < 0} -x_i/d_i.
\]

\subsection{Stopping criteria}\label{Sec:StoppingCriteria}

We now look at stopping criteria for optimizing $f(x)$ as defined in
\eqref{Eq:ObjectiveF} over the weighted one-norm ball. A common
stopping criterion for problem of this type is to look at the relative
norm of the projected gradient:
\[
\rho(x) := \frac{\norm{\mathcal{P}_{\mathcal{C}}(x - \nabla f(x)) -
    x}_2}{\max\{1,\norm{\nabla f(x)}\}},
\]
which is zero if and only if $x$ is optimal. In addition to this we
can look at the relative duality gap, which we define as the
difference $\delta$ between $f(x)$ and any dual feasible objective,
divided by $\max\{1,f(x)\}$.

For the derivation of the dual problem we follow
\cite{BER2008Fb,BER2010Fa} and rewrite the original problem as:
\[
\minimize{x,r}\quad \half\norm{r}_2^2 + c^Tx +
\sfrac{\mu}{2}\norm{x}_2^2\quad\st\quad Ax+r-b = 0,\ \norm{x}_{w,1}\leq \tau.
\]
The dual of this problem is given by
\[
\maximize{y,\lambda}\quad \mathcal{L}(y,\nu)\quad\st\quad \lambda\geq 0,
\]
where the Lagrange dual function $\mathcal{L}$ is given by
\begin{eqnarray}
\mathcal{L}(y,\lambda)
& := & \inf_{x,r} \Big\{\half\norm{r}_2^2 + c^Tx +
\sfrac{\mu}{2}\norm{x}_2^2 - y^T(Ax+r-b) +
\lambda(\norm{x}_{w,1}-\tau) \Big\} \notag\\[2pt]
& = & y^Tb - \tau\lambda + \inf_{r}\Big\{ \half\norm{r}_2^2 - y^Tr\Big\} +
\inf_{x}\Big\{ (c - A^Ty)^Tx + \sfrac{\mu}{2}\norm{x}_2^2 +
  \lambda\norm{x}_{w,1}\Big\} \notag\\[2pt]
& = & y^Tb - \tau\lambda -\half\norm{y}_2^2+
\inf_{x}\left\{ (c - A^Ty)^Tx + \sfrac{\mu}{2}\norm{x}_2^2 +
  \lambda\norm{x}_{w,1}\right\}.\label{Eq:DualWeightedL1}
\end{eqnarray}
Here, the infimum over $r$ is solved by equating the gradient to zero,
giving $y=r$ and $y^Tr = \norm{y}_2^2$. For the infimum over $x$ we
consider two cases, based on the value of $\mu$.

\paragraph{Dual when $\mu = 0$.} With $c=0$ this is exactly to the
formulation considered in \cite{BER2010Fa}, and with minor changes it
can be shown that
\[
\inf_{x} \{ (c - A^Ty)^Tx + \lambda\norm{x}_{w,1} \} = \begin{cases} 0 &
  \norm{A^Ty - c}_{\frac{1}{w},\infty} \leq \lambda \\
-\infty & \mathrm{otherwise}.\end{cases}
\]
From this we then obtain the dual problem:
\begin{equation}\label{Eq:DualOneNorm1}
\maximize{y,\ \lambda\geq0}\quad y^Tb-\tau\lambda -\half\norm{y}_2^2\quad\st\quad
\norm{A^Ty - c}_{\frac{1}{w},\infty} \leq \lambda.
\end{equation}
As a dual-feasible point we can choose $y=r$. For any given $y$ it can be
verified that choosing $\lambda = \norm{A^Ty -
  c}_{\frac{1}{w},\infty}$ always gives the largest dual objective
value. Given $x$ and the corresponding residual $r = b - Ax$ we
therefore obtain the following duality gap:
\[
\delta = \norm{r}_2^2 + c^Tx - r^Tb + \tau\norm{A^Tr-c}_{\frac{1}{w},\infty}
\]

\paragraph{Dual when $\mu > 0$.} The simplest way of dealing with $\mu
> 0$ is to rewrite the problem as:
\[
\minimize{x}\quad \half\norm{\tilde{A} x - \tilde{b}}_2^2\quad\st\quad\norm{x}_{w,1}\leq\tau
\]
with $\tilde{A} = [A;\ \!\sqrt{\mu}I]$, and $\tilde{b} = [b;\
\!0]$. This reduces the problem to the form where $\mu=0$ and we can
therefore directly use dual formulation
\eqref{Eq:DualOneNorm1}. Choosing $y = \tilde{r}$, with $\tilde{r} =
[r;\ \!-\sqrt{\mu}x]$, and applying the same derivation as given above,
we obtain a dual objective value of
\begin{equation}\label{Eq:DualObjMu1}
r^Tb -\tau\lambda -\half\norm{r}_2^2
-\sfrac{\mu}{2}\norm{x}_2^2,\quad\mbox{with}\quad
\lambda=\norm{A^Tr-\mu x-c}_{\frac{1}{w},\infty},
\end{equation}
and a corresponding duality gap of
\begin{equation}\label{Eq:DualGapMu1}
\delta = \norm{r}_2^2 + c^Tx - r^Tb + \mu\norm{x}_2^2 + \tau\norm{A^Tr - \mu
  x - c}_{\frac{1}{w},\infty}.
\end{equation}

Another approach is to solve the original infimum over $x$ in
\eqref{Eq:DualWeightedL1} for the case where $\mu > 0$. For a fixed
$y$ and $\lambda$ we have
\begin{eqnarray}
m(y,\lambda) & := & \inf_{x}\left\{ (c - A^Ty)^Tx + \sfrac{\mu}{2}\norm{x}_2^2 +
  \lambda\norm{x}_{w,1}\right\} \notag\\
& = & 
\mu\inf_{x} \left\{ -\sfrac{1}{\mu}(A^Ty - c)^Tx + \half\norm{x}_2^2 +
  \sfrac{\lambda}{\mu}\norm{x}_{w,1}\right\} \label{Eq:Mylambda2}
\end{eqnarray}
When $\lambda = 0$ it is easily seen that $x^* =
\sfrac{1}{\mu}(A^Ty-c)$, thus giving $m(x) =
-\sfrac{1}{2\mu}\norm{A^Ty - c}_2^2$. For the more general case where $\lambda >
0$, we first reformulate \eqref{Eq:Mylambda2} as
\begin{equation}\label{Eq:Mylambda3}
m(y,\lambda) = \mu \inf_{x}\left\{ -v^Tx + \half{\norm{x}}_2^2 + h(x)\right\}.
\end{equation}
with $h(x) = \frac{\lambda}{\mu}\norm{x}_{w,1} =
\norm{x}_{\frac{\lambda w}{\mu},1}$ and $v =
\frac{1}{\mu}(A^Ty-c)$. For problems of the form \eqref{Eq:Mylambda3}
we have:

\begin{btheorem}\label{Thm:MinAsProxNorm}
Let $h(\cdot)$ be any norm then
\[
\inf_{x}\ -v^Tx + \half\norm{x}_2^2 + h(x) = -\half\norm{\prox_{h}(v)}_2^2
\]
\end{btheorem}
\begin{proof}
 Note that the objective is coercive and therefore attains the
  minimum. This allows us to rewrite and solve the objective as follows:
\[
u = \argmin_{x}\ \half\norm{x-v}_2^2 + h(x) = \prox_h(v).
\]
We then need to show that
\[
-v^Tu + \half\norm{u}_2^2 + h(u) = -\half\norm{u}_2^2
\]
From the Moreau decomposition we have $v =\prox_h(v) +
\prox_{h^*}(v)$, where $h^*$ is the conjugate of $h$. Using
$\prox_{h^*}(v) = v-u$ we have
\begin{eqnarray*}
-v^Tu + \half\norm{u}_2^2 + h(u) & = & -(u + (v-u))^Tu + \half u^T u
 + h(u)\\
& = & -\half\norm{u}_2^2 - (v-u)^Tu + h(u) \\
& = & -\half\norm{\prox_h(v)}_2^2 - \prox_{h^*}(v)^T\prox_{h}(v) +
h(\prox_h(v)) \\
& = & -\half\norm{\prox_h(v)}_2^2,
\end{eqnarray*}
where the last equality follows from Lemma~\ref{Lemma:NormOfProx}
given below.
\end{proof}

\begin{blemma}\label{Lemma:NormOfProx}
  Let $h(\cdot)$ be any norm with conjugate $h^*(\cdot)$,
  then
\[h(\prox_h(x)) = \prox_{h^*}(x)^Tprox_h(x).\]
\end{blemma}
\begin{proof}
Let $u = \prox_h(x)$, then it is well known that $\prox_{h^*}(x) =
x-u$ and $x -u \in \partial h(u)$. It thus remains to show that
$h(u) = (x-u)^Tu$.  The subgradient $\partial h(u)$ and norm $h$ are
defined in terms of the dual norm $h_*$ as
\[
\partial h(u) := \argmax_{w: h_*(w)\leq 1}\ w^Tu,\quad\mbox{and}\quad
h(u) = \max_{w: h_*(w)\leq 1}\ w^Tu,
\]
respectively, which means that $h(u) = w^Tu$ for any $w \in \partial h(u)$. Choosing
$w = x-u$ gives $h(u) = w^Tu = (x-u)^Tu$, as desired.
\end{proof}

Application of Theorem~\ref{Thm:MinAsProxNorm} to \eqref{Eq:Mylambda3}
with proximal operator (see also \eqref{Eq:ProxWeightedL1})
\[
\prox_h(v) = \mathrm{sign}(v)\big[\abs{v} - \sfrac{\lambda w}{\mu}\big]_+,
\]
we obtain
\[
m(y,\lambda) = -\frac{\mu}{2} \left\Vert \sign(\sfrac{1}{\mu}(A^Ty -
    c))\left[\big\vert{\sfrac{1}{\mu}(A^Ty - c)}\big\vert -
      \frac{\lambda w}{\mu}\right]_+\right\Vert_2^2
= -\frac{1}{2\mu}\left\Vert\left[\big\vert A^Ty - c\big\vert - \lambda w\right]_+\right\Vert_2^2.
\]
The same expression holds for $\lambda = 0$ and substitution into
\eqref{Eq:DualWeightedL1} therefore gives the following dual problem:
\begin{equation}\label{Eq:DualOneNorm2}
\maximize{y,\ \lambda \geq 0}\quad y^Tb  - \tau\lambda -\half\norm{y}_2^2 -\sfrac{1}{2\mu}\left\Vert\left[\big\vert A^Ty - c\big\vert - \lambda w\right]_+\right\Vert_2^2.
\end{equation}
Even when restricting $y$ to the current residual $r$ in the primal
formulation, we can show that the value of \eqref{Eq:DualOneNorm2} is
never smaller than that of \eqref{Eq:DualObjMu1} and, consequently,
that the duality gap never exceeds the value in \eqref{Eq:DualGapMu1}.
Choosing $\lambda = \norm{A^Ty + \mu x -c}_{\frac{1}{w},\infty}$,
means that for any index $i$ we have
\[
\lambda \geq \sfrac{1}{w_i}\left\vert [A^Tx - c]_i + \mu x_i\right\vert \geq
\sfrac{1}{w_i}\left\vert [A^Tx-c]_i\right\vert - \sfrac{\mu}{w_i}\abs{x_i}
\]
Multiplying either side by $w_i$ and rearranging gives
\[
\left\vert [A^Ty - c]_i\right\vert - \lambda w_i\leq \mu\abs{x_i}.
\]
Because the right-hand side is always nonnegative, this continues to
hold when applying the $[\cdot]_+$  operator on the left-hand side,
and as a result we have
\[
\sfrac{1}{2\mu}\left\Vert\left[\big\vert A^Ty -
    c\big\vert - \lambda w\right]_+\right\Vert_2^2
\leq \sfrac{\mu}{2}\norm{x}_2^2,
\]
from which the desired result immediately follows.

\paragraph{Finding a dual-feasible solution.} 

It follows from Slater's condition and strong duality that, at the
solution ($x$, $r$) for \eqref{Eq:DualOneNorm1}, we have $y = r$ and,
without loss of generality, $\lambda =
\norm{A^Tr-c}_{\frac{1}{w},\infty}$. When $r$ is not optimal, we
can still choose $y=r$ and obtain a dual-feasible solution. For
\eqref{Eq:DualOneNorm2} we can also take $y=r$, but finding $\lambda$
requires some more work. In general, given any $y$ we want to find a
$\lambda$ that maximizes the objective. Writing $z = \abs{A^Tr-c}$ and
ignoring constant terms, this is equivalent to solving
\[
\lambda^* := \argmin_{\lambda \geq 0}\quad \tau\lambda +
\sfrac{1}{2\mu}\norm{[z-\lambda w]_+}_2^2.
\]
With $\mathcal{I}(\lambda) := \{i \mid z_i \geq \lambda
w_i\}$ we can write the objective as
\[
f(\lambda) = \tau\lambda +
\sfrac{1}{2\mu}\sum_{i\in\mathcal{I}(\lambda)}(z_i - \lambda w_i)^2
\]
Discarding all zero terms with $z_i = 0$, this function is piecewise
quadratic with breakpoints at $\lambda_i = w_i / z_i$. We can write
the the sequence of breakpoints in non-decreasing order as
$\lambda_{[i]}$ for $i=0,\ldots,n$ , with $\lambda_{[0]} := 0$. The
gradient between successive breakpoints is linear and continuously
increases from $f'(0) = \tau - 1/\mu\sum_i w_iz_i$ to
$f'(\lambda_{[n]}) = \tau$. In order to find the optimal point
$\lambda^*$, we consider two cases. In the first case we have $f'(0)
\geq 0$, or equivalently $\tau \geq 1/\mu\sum_i w_iz_i$, which means
that the function is non-decreasing and we find $\lambda^* = 0$. In
the second case we need to find $\lambda$ for which the gradient
vanishes. This can be done by traversing the breakpoints until the
first breakpoint is found where the gradient is nonnegative. The
desired solution $\lambda^*$ is then found by linear interpolation
over the last segment. Including sorting this can be done in
$\mathcal{O}(n\log n)$ time. This problem is very similar to
projection onto the one-norm ball, and can also be evaluated in
expected $\mathcal{O}(n)$ time using an algorithm similar to that
proposed in \cite{DUC2008SSCa}.

\paragraph{Primal-dual pairs.}

Using the methods described above, we can compute an upper bound on
the duality gap given a feasible $x$ and the corresponding residual
$r$. We can do at least as good, and often better, by maintaining the
maximum dual objective found so far, and using this to determine the
relative duality gap. This way it is possible for the primal objective
for the current iterate $x$ to attain the desired optimality tolerance
while the corresponding dual estimate is far from optimal. Within the
root-finding framework this means that we cannot simply use the latest
residual $r$ to evaluate the gradient of the Pareto curve. Instead we
should maintain the value of $\lambda$ corresponding to the best dual
solution at any point, and use this as the gradient approximation. In
other words, we need to keep track of the best primal and dual
variables separately.


\section{Numerical Experiments}\label{Sec:NumericalExperiments}\label{Sec:Numerical}
In this section we evaluate the performance of the hybrid approach on
the Lasso problem \eqref{Eq:Lasso} both independently and within the
\spgl\ root-finding framework~\cite{BER2008Fb} described in the
introduction. The \spgl\ solver can be used both for stand-alone Lasso
problems, as well as for basis-pursuit denoise \eqref{Eq:BPDN}
problems. For the hybrid method we are mostly concerned with the
performance and the former and we therefore changed \spgl\ in two
stages. First we modified the stopping criteria used in the Lasso
mode, now declaring a solve successful only if the relative duality
falls below a certain tolerance level. We then added all modifications
needed for the implementation of the hybrid approach. To distinguish
between the different algorithms, we use the convention that \spgl\ is
used only to refer to the existing implementation provided by
\cite{BER2008Fb}. We refer to the version of \spgl\ with the more
stringent stopping criteria as the \spg{} method, which is then
extended with the techniques described in this paper to obtain the
hybrid method.

When used in the root-finding mode to solve \eqref{Eq:BPDN}, \spgl\
uses several different criteria to decide when to update $\tau$. Each
subproblems in \spgl\ is considered solved when the relative change in
objective is small, and at least one iteration was taken within the
current subproblem. The overall problem is declared solved when
$\norm{A^Ty}_{\infty}$, the relative difference between $\norm{r}_2$
and $\sigma$, or the relative duality gap is sufficiently small. For the
basis-pursuit denoise experiments based on the \spg{} and hybrid
algorithms, we use a separate implementation of the root-finding
framework in which each Lasso subproblem is fully solved before
updating $\tau$. The differences in stopping criteria, and especially
the lack of guarantees on the duality gap for the final subproblem in
\spgl, make it difficult to compare the performances directly. We
therefore focus predominantly on how the performance of the hybrid
method differs from the reference \spg{} method.

\subsection{Lasso on sparse problems}

\begin{table}[t]
\centering\small\setlength{\tabcolsep}{4.5pt}
\begin{tabular}{rp{5pt}>{\raggedleft\arraybackslash}p{21pt}>{\raggedleft\arraybackslash}p{21pt}p{0pt}>{\raggedleft\arraybackslash}p{21pt}>{\raggedleft\arraybackslash}p{21pt}p{0pt}>{\raggedleft\arraybackslash}p{21pt}>{\raggedleft\arraybackslash}p{21pt}p{5pt}>{\raggedleft\arraybackslash}p{16pt}>{\raggedleft\arraybackslash}p{16pt}p{0pt}>{\raggedleft\arraybackslash}p{16pt}>{\raggedleft\arraybackslash}p{16pt}p{0pt}>{\raggedleft\arraybackslash}p{16pt}>{\raggedleft\arraybackslash}p{16pt}}
&&\multicolumn{8}{c}{Runtime (sec.)}&&\multicolumn{8}{c}{Iterations}\\
\cline{1-1}\cline{3-10}\cline{12-19}
$k$&&\multicolumn{2}{c}{Backtrack}&&\multicolumn{2}{c}{Local}&&\multicolumn{2}{c}{Global}&&\multicolumn{2}{c}{Backtrack}&&\multicolumn{2}{c}{Local}&&\multicolumn{2}{c}{Global}\\
&&\multicolumn{2}{c}{}&&\multicolumn{2}{c}{trajectory}&&\multicolumn{2}{c}{trajectory}&&\multicolumn{2}{c}{}&&\multicolumn{2}{c}{trajectory}&&\multicolumn{2}{c}{trajectory}\\
\cline{1-1}\cline{3-4}\cline{6-7}\cline{9-10}\cline{12-13}\cline{15-16}\cline{18-19}
\\[-9pt]
 50&& 0.14& 0.17&& 0.46& 0.40&& 0.53& 0.46&&  25&  29&&  22&  24&&  22&  24\\
   && 0.13& 0.17&& 0.44& 0.39&& 0.52& 0.46&&  24&  28&&  20&  22&&  20&  23\\
   && 0.14& 0.18&& 0.45& 0.37&& 0.50& 0.46&&  25&  30&&  21&  22&&  21&  22\\[3pt]
100&& 0.22& 0.26&& 0.69& 0.58&& 0.93& 0.82&&  40&  44&&  35&  34&&  34&  34\\
   && 0.22& 0.25&& 0.68& 0.57&& 0.92& 0.82&&  40&  43&&  34&  33&&  33&  33\\
   && 0.23& 0.26&& 0.69& 0.57&& 0.96& 0.83&&  41&  45&&  34&  34&&  34&  34\\[3pt]
150&& 0.34& 0.38&& 0.97& 0.83&& 1.56& 1.39&&  59&  62&&  54&  54&&  53&  53\\
   && 0.33& 0.36&& 0.96& 0.80&& 1.50& 1.41&&  57&  59&&  51&  51&&  50&  51\\
   && 0.37& 0.40&& 0.94& 0.80&& 1.57& 1.39&&  62&  64&&  53&  53&&  50&  51\\[3pt]
200&& 0.55& 0.57&& 1.52& 1.32&& 2.81& 2.49&&  88&  88&&  93&  92&&  89&  88\\
   && 0.53& 0.55&& 1.50& 1.29&& 2.78& 2.55&&  86&  87&&  95&  93&&  88&  87\\
   && 0.59& 0.61&& 1.52& 1.31&& 2.88& 2.64&&  95&  94&& 100&  98&&  87&  88\\[3pt]
\multicolumn{19}{c}{({\bf{a}}) Tolerance $10^{-4}$}\\[10pt]
 50&& 0.16& 0.19&& 0.52& 0.42&& 0.59& 0.48&&  28&  32&&  25&  26&&  25&  26\\
   && 0.16& 0.19&& 0.51& 0.41&& 0.58& 0.48&&  28&  31&&  24&  25&&  24&  26\\
   && 0.24& 0.20&& 0.53& 0.39&& 0.57& 0.45&&  32&  33&&  25&  25&&  25&  25\\[3pt]
100&& 0.34& 0.28&& 0.75& 0.60&& 1.01& 0.84&&  48&  47&&  38&  38&&  37&  37\\
   && 0.25& 0.28&& 0.72& 0.59&& 0.97& 0.84&&  44&  47&&  36&  37&&  36&  37\\
   && 0.31& 0.29&& 0.76& 0.62&& 1.05& 0.88&&  48&  48&&  37&  38&&  37&  38\\[3pt]
150&& 0.50& 0.42&& 1.08& 0.89&& 1.69& 1.46&&  70&  67&&  59&  59&&  57&  58\\
   && 0.45& 0.39&& 1.00& 0.85&& 1.57& 1.38&&  67&  64&&  55&  56&&  54&  55\\
   && 0.60& 0.44&& 1.04& 0.82&& 1.61& 1.41&&  75&  69&&  58&  58&&  54&  56\\[3pt]
200&& 0.84& 0.60&& 1.70& 1.34&& 2.92& 2.52&& 107&  95&& 102&  98&&  97&  94\\
   && 0.80& 0.61&& 1.65& 1.33&& 2.96& 2.58&& 104&  94&& 102& 100&&  96&  94\\
   && 0.88& 0.65&& 1.70& 1.37&& 3.06& 2.66&& 112& 101&& 107& 105&&  95&  95\\[3pt]
\multicolumn{19}{c}{({\bf{b}}) Tolerance $10^{-6}$}\end{tabular}
\caption{Comparison of different line-search methods. From top to bottom for each sparsity
  level are the results for random Gaussian, uniform, and discrete
  sign support. The left and right columns correspond to the \spg{}
  and the hybrid methods, respectively. Runtimes and iterations are averaged over twenty
  problem instances.}\label{Table:Trajectory}
\end{table}

In the first set of experiments we compare the performance of
backtracking line search and line search along the entire projection
trajectory for both the \spg{} and hybrid method. For the trajectory
line search we return either with the first local minimum or with the
global minimum. For the test problems we follow a conventional
compressed-sensing scenario where $A$ is a random $1024\times 2048$
matrix with i.i.d.~normal entries with columns normalized to unit
norm. We set $b = Ax_0$, for $k$-sparse vectors $x_0$ with random
support and entries generated i.i.d.~from (1) the normal distribution;
(2) the uniform distribution over $[-1,1]$; and (3) the discrete set
$\{-1,+1\}$ with equal probability. We set $\tau = 0.995\cdot
\norm{x_0}_1$ and terminate the algorithm whenever the relative
duality gap, computed as $(f(x) - f_{\mathrm{dual}}) / \max\{f(x),
10^{-3}\}$, falls below a given tolerance level.

Table~\ref{Table:Trajectory} shows the average runtime and number of
iterations for twenty random instances of each test problem for
different sparsity levels $k$, and each of the three distributions
used for the on-support values in $x_0$. Comparing across the line
search methods we see that the two trajectory line search methods
require fewer iterations than backtracking to converge. The
backtracking line search, on the other hand, has a lower per-iteration
cost and overall outperforms the trajectory line search uniformly in
terms of runtime.  Looking at the difference between the \spg{}
and hybrid methods, we see that the number of iterations required by
the hybrid method is larger than that of the \spg{} method for the
lower optimality tolerance and smaller values of $k$. Combined with
the backtracking line search this means that the runtime of the hybrid
method is slightly larger compared to the \spg{} method. For the
trajectory line search we see, perhaps somewhat surprisingly, that the
runtime of the hybrid method is uniformly lower than the \spg{} 
method, despite the larger number of iterations. The reason for this
is that the line search for each quasi-Newton iteration taken by the
hybrid method is much faster than the trajectory line search, thereby
reducing the overall runtime.  Comparing between the two optimality
tolerance levels we note that the hybrid method does well for the
lower tolerance level of $10^{-6}$. The \spg{} method is not only
slower for these problems, but also suffers from a problem where the
line search fails to find a feasible step length before the desired
optimality tolerance is reached, thereby terminating the optimization
prematurely. This problem did not occur for a tolerance level of
$10^{-4}$, but for $10^{-6}$ this happened on one of the twenty
problems for $k=50$ and gradually went up to five out of twenty for
$k=200$. No such line-search errors occurred in the hybrid method.

\begin{table}
\centering\small\setlength{\tabcolsep}{5.5pt}
\begin{tabular}{rcrrrrrcrrrrrcr}
$k$ && \multicolumn{3}{c}{\ \ Runtime \spg{} (s)}& rel.gap&\ding{51}&&\multicolumn{3}{c}{\ \ Runtime hybrid (s)}&rel.gap&\ding{51}&&\multicolumn{1}{c}{(\%)}\\
\cline{1-1}\cline{3-7}\cline{9-13}\cline{15-15}
 50&&  0.17&  0.16&  0.16& 3.8e-7&100&&  0.20&  0.19&  0.19& 4.0e-7&100&& -17\\
 75&&  0.24&  0.21&  0.21& 5.8e-7& 99&&  0.24&  0.23&  0.23& 4.7e-7&100&& -7\\
100&&  0.28&  0.27&  0.31& 5.8e-7& 95&&  0.28&  0.27&  0.28& 5.3e-7&100&& 3\\
125&&  0.39&  0.36&  0.37& 6.5e-7& 91&&  0.34&  0.32&  0.31& 6.2e-7&100&& 14\\
150&&  0.44&  0.42&  0.45& 7.3e-7& 92&&  0.39&  0.39&  0.37& 6.7e-7&100&& 13\\
175&&  0.64&  0.53&  0.55& 7.7e-7& 81&&  0.46&  0.44&  0.45& 6.9e-7&100&& 21\\
200&&  0.74&  0.66&  0.73& 8.2e-7& 75&&  0.58&  0.54&  0.53& 7.1e-7&100&& 23\\
225&&  0.90&  0.88&  0.86& 8.5e-7& 71&&  0.71&  0.68&  0.65& 7.5e-7&100&& 23\\
250&&  1.30&  1.10&  1.17& 9.7e-7& 55&&  0.93&  0.86&  0.84& 7.6e-7&100&& 26\\
275&&  1.65&  1.54&  1.42& 1.0e-6& 50&&  1.30&  1.14&  1.07& 8.3e-7&100&& 24\\
300&&  2.39&  2.01&  1.85& 2.1e-6& 35&&  1.86&  1.53&  1.44& 8.3e-7&100&& 23\\
325&&  3.60&  2.76&  2.66& 4.7e-6& 26&&  3.08&  2.04&  2.00& 8.8e-7&100&& 22\\
350&&  6.61&  4.31&  3.91& 7.6e-6& 15&&  5.68&  3.30&  3.08& 9.0e-7&100&& 19\\
375&& 31.82&  8.27&  7.15& 1.5e-5&  7&& 20.09&  6.19&  5.68& 9.1e-7&100&& 28\\
400&& 218.46& 23.16& 16.16& 3.0e-5&  9&& 129.42& 14.65& 10.93& 9.5e-7& 89&& 37\\
\end{tabular}
\caption{Comparison between the \spg{} and the hybrid method on
  exact sparse problems with $k$ non-zero entries. The first three
  columns for either method give the average runtime over 50
  instances when the non-zero entries are sampled i.i.d.~from
  respectively the discrete $\{-1,1\}$, uniform(-1,1), and
  the normal distribution.  The fourth column gives the median relative duality
  gap at the final iteration taken over all 150 problem
  instances and should be compared with the optimality tolerance, which was set to
  $10^{-6}$. The fifth column for each of the two blocks, indicated by
  the check mark, gives the percentage of runs that completed
  successfully, that is, completed without a line-search error. The right-most
  column gives the average of the speed up values for each of the three
  distributions.
}\label{Table:SparseX0}
\end{table}

From the results in Table~\ref{Table:Trajectory}, along with various
other experiments not shown here, it was found that backtracking line
search outperforms the trajectory line search. As a result, we only
consider the former throughout the remainder of this section. We now
take a closer look at the occurrence of line-search errors and the
speed up obtained using the hybrid method. For this, we modify the
earlier setup by increasing the range of sparsity levels $k$ and
choosing $\tau = 0.99\cdot\norm{x}_1$. We run 50 instances for each of
the three distributions used above and report in
Table~\ref{Table:SparseX0} the results obtained with an optimality
tolerance of $10^{-6}$. In general we see that the runtime goes up
considerable as we keep increasing $k$. Moreover, the results show
clear differences in the runtime for the three distributions with a
much higher runtime for problems based on sparse vectors with $\pm 1$
entries. For sparsity levels up to around one hundred the number of
iterations in the \spg{} method is relatively small (between 25 and
50). For these problems the hybrid method may complete before or soon
after the first quasi-Newton step is taken. The slight overhead of the
method and occasionally a small number of additional iterations make
the hybrid method somewhat slower on average for these problems than
the \spg{} method.  For larger values of $k$, the number of
iterations goes up, and the effect of the quasi-Newton steps in the
hybrid method becomes apparent with average speed up values between 20
and 30\%. Aside from reduced runtime we see from
Table~\ref{Table:SparseX0} that the hybrid method also manages to
solve problems to the desired accuracy level much better than the
\spg{} method. The number of solved problems steadily falls to
around 9\% with increasing $k$ for the \spg{} method, but remains at
100\% for all but the largest $k$ for the hybrid method. The median
relative duality gap provides further information about the level of
accuracy reached before the algorithm completes or terminates with a
line-search error. For the largest values of $k$, the \spg{} method
fails to complete with a relative duality gap of even $10^{-5}$ for at
least half of the problems.

\subsection{Root finding}

Given that most of the runtimes that appear in
Tables~\ref{Table:Trajectory} and \ref{Table:SparseX0} are of the
order of seconds, it is valid to question whether these problems are
too idealized and well behaved to give a good idea about the practical
performance of the algorithms. In this section we therefore look at
two different types of problems. First we introduce a class of random
problems that better reflect conditions found in practical
problems. Second we evaluate the performance on the
Sparco~\cite{BER2009FHHa} collection of test problems for sparse
reconstruction.

\subsubsection{Coherent test problem generation}\label{Sec:CoherentA}

\begin{figure}[t]
\centering
\begin{tabular}{ccc}
\multicolumn{3}{c}{\includegraphics[width=0.9\textwidth]{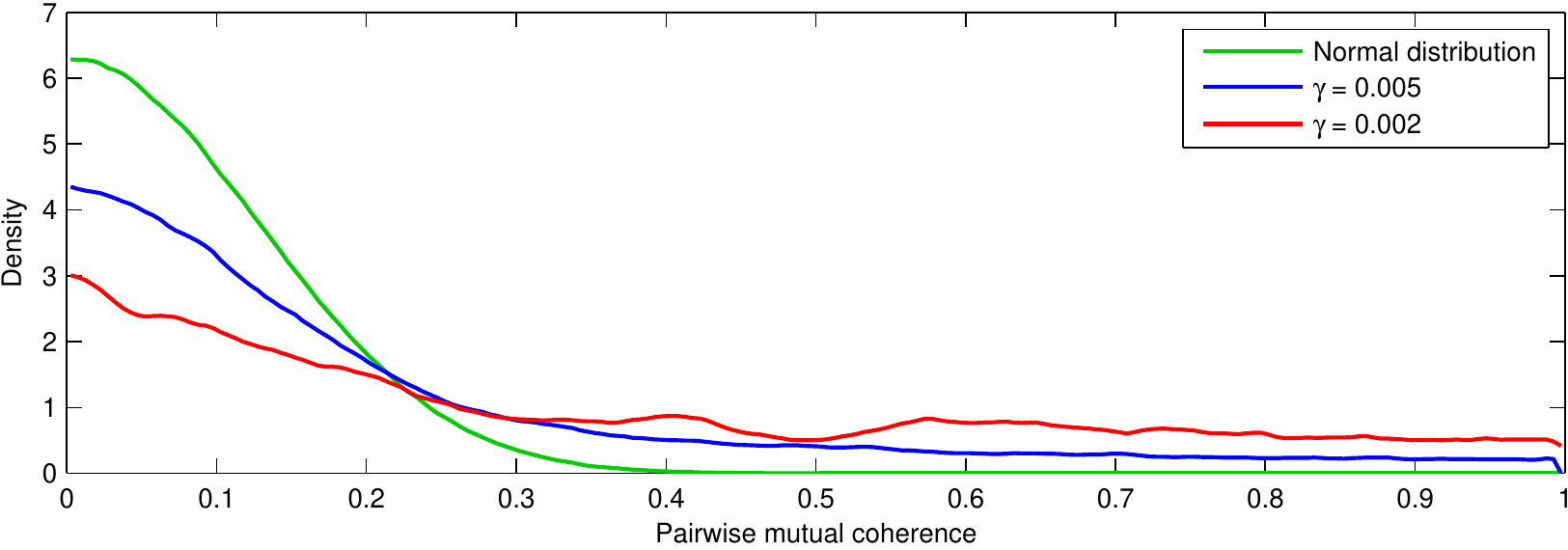}}\\
\multicolumn{3}{c}{({\bf{a}})}\\[3pt]
\includegraphics[height=0.23\textwidth]{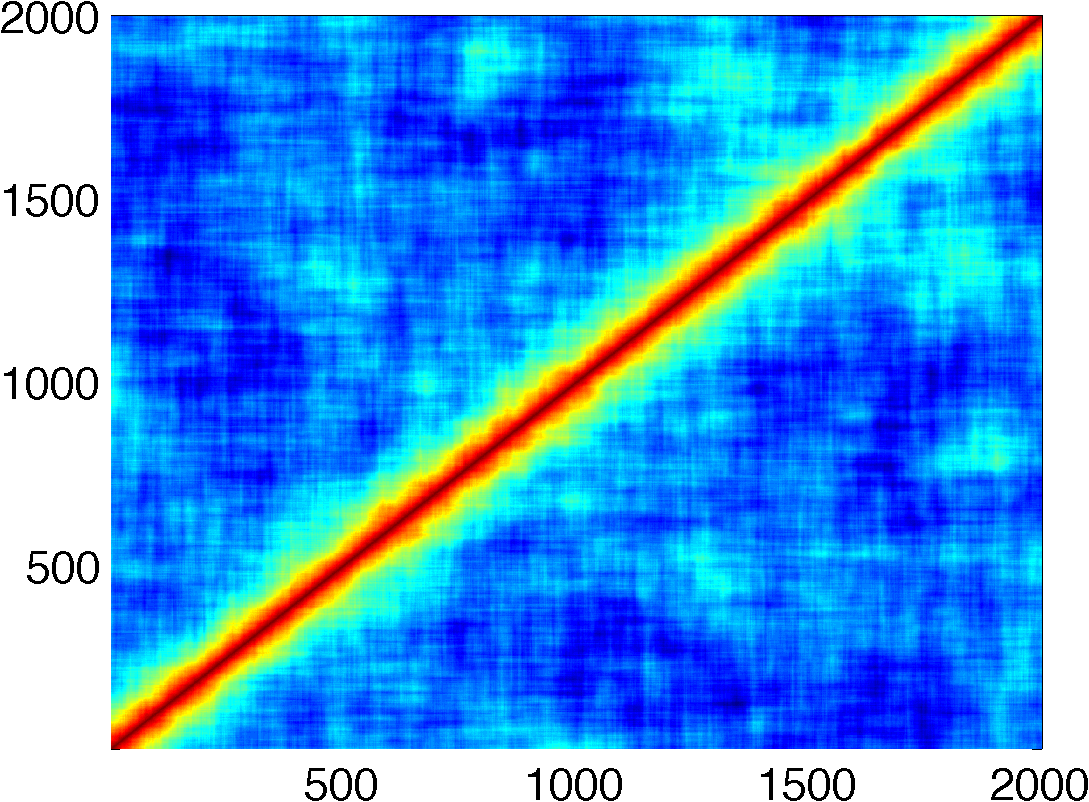}&
\includegraphics[height=0.23\textwidth]{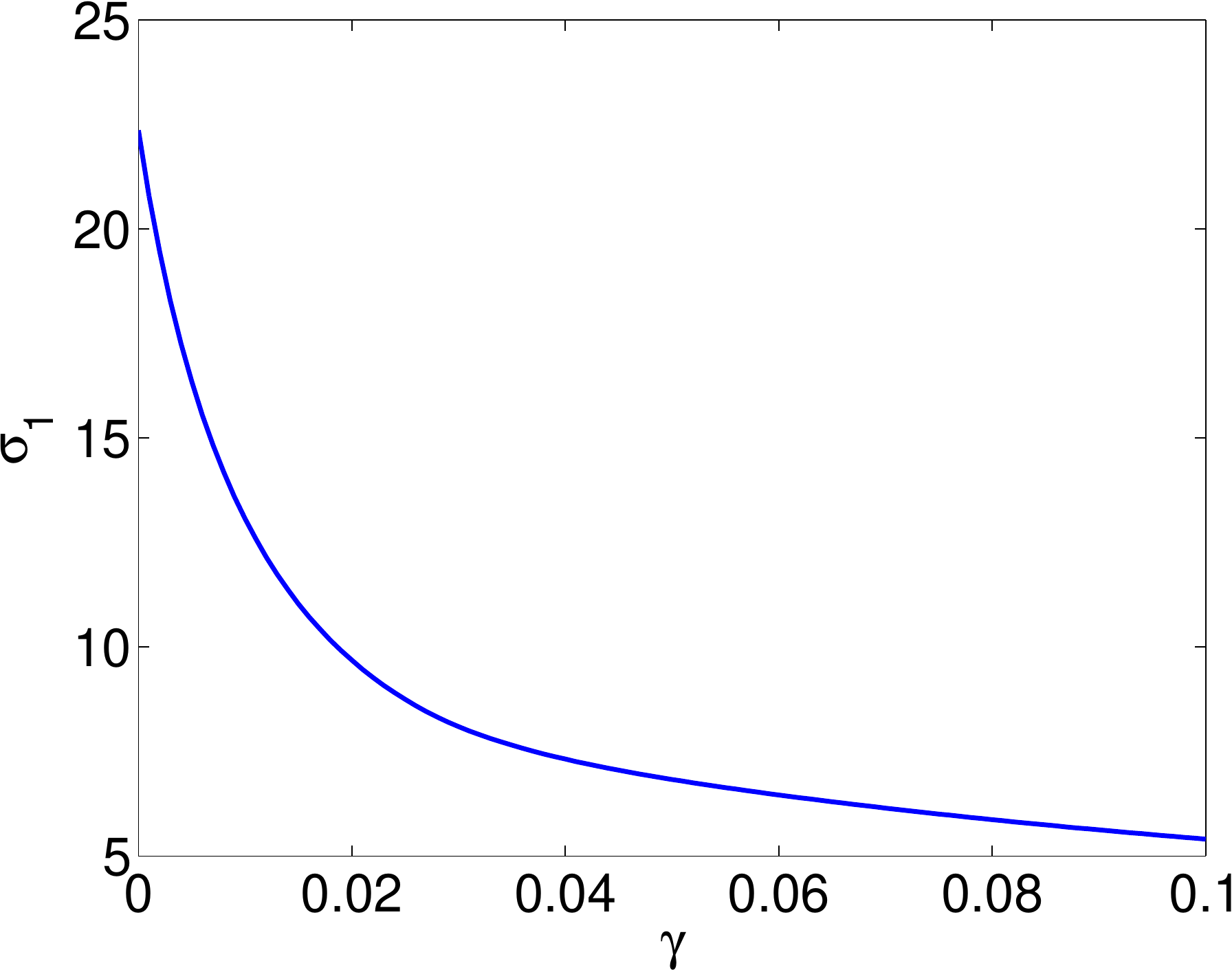}&
\includegraphics[height=0.23\textwidth]{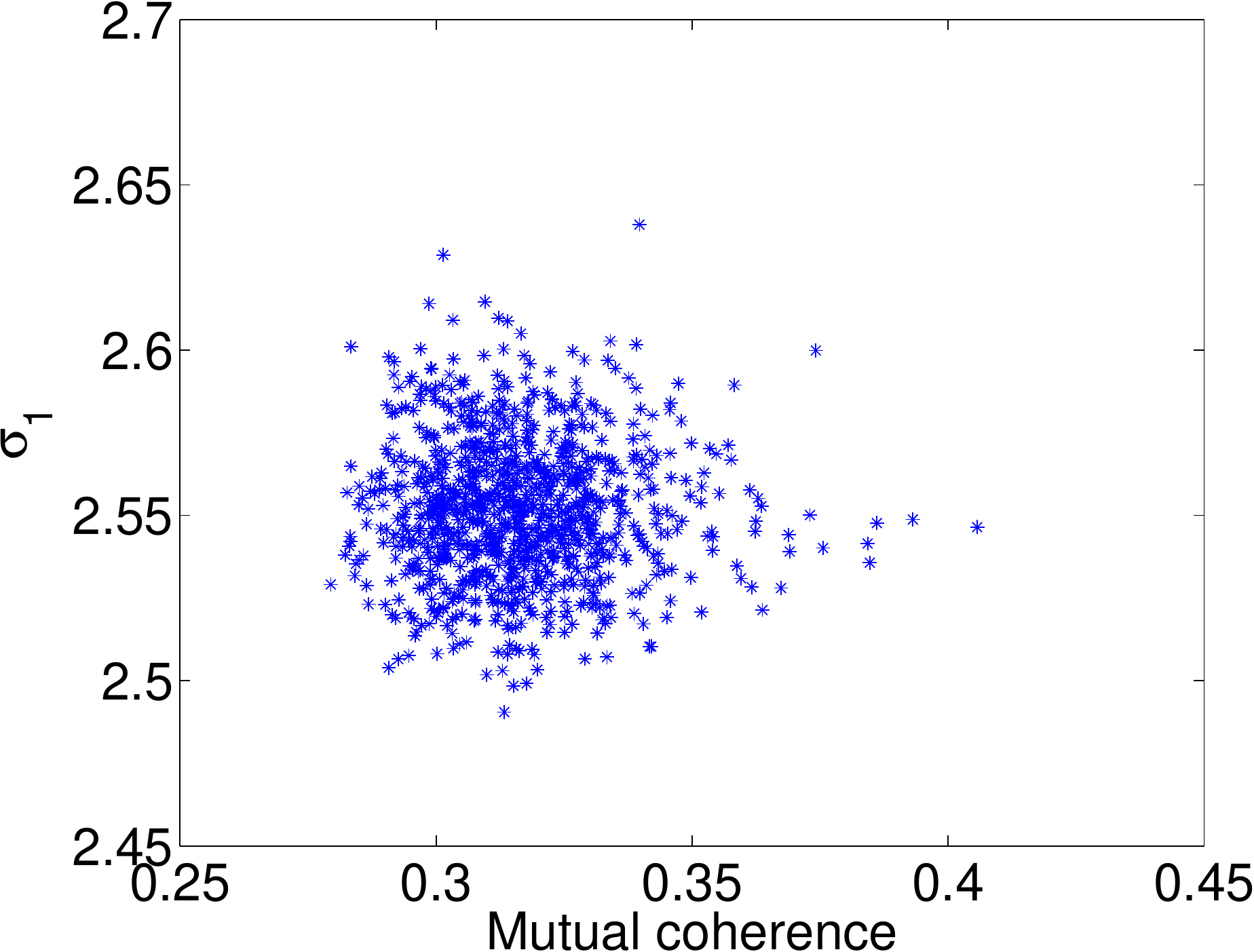}\\
({\bf{b}}) & ({\bf{c}}) & ({\bf{d}})
\end{tabular}
\caption{(a) Distribution of pairwise mutual coherence between vectors
  of two types of $64\times 2048$ matrices with unit-normalized
  columns generated as: (i) random vectors with i.i.d.~normal entries,
  and (ii) a random walk over the sphere with the mutual coherence
  between successive columns equal to $1-\gamma$; (b) Gram matrix of a
  $200\times 2000$ matrix generated with $\gamma = 0.01$; (c) the top
  singular value of a $200\times 500$ matrix generated with the same
  $a_1$ and $v_k$ for different values of $\lambda$; and (d) the
  mutual coherence and top singular value for 1000 random Gaussian
  $200\times 500$ matrices with columns scaled to unit
  norm.}\label{Fig:RandomWalk}
\end{figure}

In the compressed-sensing literature it is well known that a random
Gaussian matrix satisfies with high probability that all sufficiently
small subsets of columns form a near-orthogonal basis for the subspace
spanned by these columns---a property known as the restricted
isometry~\cite{CAN2005Tb}. Another quantity used to characterize
matrices is the mutual coherence, defined as the maximum absolute
pairwise cosine distance between the columns. In practical
applications matrix $A$ is often more coherent
\cite{CAN2011ENRa}. Although there are no theoretical results on how
this affects the complexity of one-norm minimization, it has been
observed empirically that more coherent problems are harder to
solve. The construction we propose for generating such problems is by
means of a random walk on the $(m-1)$-sphere with a step size
parameterized by $\gamma$. Starting with a unit norm column $a_1$ we
construct successive columns by sampling a vector $v_k$ with
i.i.d.~Gaussian entries and setting $a_{k+1} = \alpha_1 a_k + \alpha_2
v_k$, where $\alpha_1$ and $\alpha_2$ are chosen such that
$\norm{a_{k+1}}_2 = 1$ and $\langle a_k, a_{k+1}\rangle =
1-\gamma$. In other words, $a_{k+1}$ lies on the boundary of a
spherical cap with center center $a_k$ and angle $\theta$ such that
$\mathrm{cos}(\theta) = 1-\gamma$. The mutual coherence of the
resulting matrix is lower bounded by $1-\gamma$, and an example of the
distribution of the pairwise cosine distance between the columns is
given in Figure~\ref{Fig:RandomWalk}(a). An example Gram matrix,
plotted in Figure~\ref{Fig:RandomWalk}(b), shows that aside from the
banded structure, there are regions of increased coherence whenever
the random walk approaches earlier locations. From
Figure~\ref{Fig:RandomWalk}(c) we see that lowering $\gamma$ while
keeping $a_1$ and $v_k$ fixed leads to an increase of the top singular
value $\sigma_1$ as the columns become more and more
similar. Figure~\ref{Fig:RandomWalk}(d) illustrates that the maximum
pairwise coherence $\mu$ does not necessarily have a relationship with
the top singular value.

\subsubsection{Highly coherent measurement matrices}

 We apply the \spg{} and hybrid method to solve \eqref{Eq:BPDN}
 using the root-finding framework explained in
 Section~\ref{Sec:Introduction}. Each Lasso subproblem \eqref{Eq:Lasso}
 is optimized to a certain optimality tolerance, and the overall
 problem is considered solved whenever the relative misfit $\abs{\sigma
   - \norm{r}_2} / \max(\sigma, 10^{-3})$ falls below $10^{-5}$. For
 completeness we also compare the performance with the
\spgl\ algorithm as provided by~\cite{BER2008Fb}.

\begin{table}
\centering\footnotesize
\setlength{\tabcolsep}{4.5pt}
\begin{tabular}{rcrcrrrcrrrcrrrcrrr}
$\gamma$ && SPGL1 && \multicolumn{3}{c}{Root finding}&& \multicolumn{3}{c}{Root finding}&& \multicolumn{3}{c}{Lasso}&& \multicolumn{3}{c}{Lasso}\\
&& runtime && \multicolumn{3}{c}{tolerance $10^{-4}$}&& \multicolumn{3}{c}{tolerance $10^{-6}$}&& \multicolumn{3}{c}{tolerance $10^{-4}$}&& \multicolumn{3}{c}{tolerance $10^{-6}$}\\
\cline{1-1}\cline{3-3}\cline{5-7}\cline{9-11}\cline{13-15}\cline{17-19}
 .100&&   0.7&&   0.7 &   0.8 & {\color{blue}{  -4}}&&   0.7 &   0.8 & {\color{blue}{  -4}}&&   0.7 &   0.7 & {\color{blue}{  -4}}&&   0.7 &   0.7 & {\color{blue}{  -3}}\\
 .050&&   1.6&&   1.7 &   1.7 & {\color{blue}{   2}}&&   1.9 &   1.7 & {\color{blue}{  12}}&&   1.8 &   1.7 & {\color{blue}{   1}}&&   1.9 &   1.8 & {\color{blue}{   7}}\\
 .020&&   6.8&&   6.6 &   5.8 & {\color{blue}{  12}}&&   7.6 &   5.9 & {\color{blue}{  22}}&&   6.5 &   6.0 & {\color{blue}{   7}}&&   7.1 &   6.1 & {\color{blue}{  14}}\\
 .010&&  17.8&&  17.7 &  14.7 & {\color{blue}{  17}}&&  21.4 &  15.0 & {\color{blue}{  30}}&&  16.7 &  15.5 & {\color{blue}{   7}}&&  18.3 &  15.6 & {\color{blue}{  15}}\\
 .005&&  51.9&&  43.6 &  36.8 & {\color{blue}{  16}}&&  52.6 &  38.2 & {\color{blue}{  27}}&&  40.6 &  37.9 & {\color{blue}{   7}}&&  42.7 &  38.2 & {\color{blue}{  11}}\\[4pt]
\multicolumn{19}{c}{({\bf{a}}) Sparsity $k=10$}\\[8pt]
 .100&&   1.8&&   3.1 &   2.7 & {\color{blue}{  11}}&&   5.0 &   3.2 & {\color{blue}{  37}}&&   1.8 &   1.7 & {\color{blue}{   7}}&&   2.5 &   1.8 & {\color{blue}{  29}}\\
 .050&&   5.1&&   8.2 &   7.1 & {\color{blue}{  13}}&&  11.3 &   7.9 & {\color{blue}{  30}}&&   4.7 &   4.3 & {\color{blue}{   8}}&&   5.8 &   4.6 & {\color{blue}{  21}}\\
 .020&&  19.8&&  29.4 &  23.2 & {\color{blue}{  21}}&&  36.5 &  25.8 & {\color{blue}{  29}}&&  17.5 &  15.5 & {\color{blue}{  11}}&&  18.9 &  16.1 & {\color{blue}{  15}}\\
 .010&&  48.5&&  77.7 &  58.9 & {\color{blue}{  24}}&&  92.5 &  63.0 & {\color{blue}{  32}}&&  45.9 &  40.0 & {\color{blue}{  13}}&&  47.9 &  41.2 & {\color{blue}{  14}}\\
 .005&& 145.7&& 221.9 & 144.9 & {\color{blue}{  35}}&& 263.4 & 153.4 & {\color{blue}{  42}}&& 127.7 & 101.9 & {\color{blue}{  20}}&& 134.9 & 104.7 & {\color{blue}{  22}}\\[4pt]
\multicolumn{19}{c}{({\bf{b}}) Sparsity $k=50$}\\[8pt]
 .100&&   3.2&&   5.2 &   4.5 & {\color{blue}{  13}}&&   7.8 &   5.3 & {\color{blue}{  32}}&&   2.8 &   2.6 & {\color{blue}{   8}}&&   3.6 &   2.8 & {\color{blue}{  23}}\\
 .050&&   8.6&&  14.9 &  12.4 & {\color{blue}{  17}}&&  19.4 &  13.8 & {\color{blue}{  29}}&&   7.6 &   6.9 & {\color{blue}{  10}}&&   8.8 &   7.3 & {\color{blue}{  17}}\\
 .020&&  29.9&&  55.1 &  41.7 & {\color{blue}{  24}}&&  65.8 &  45.2 & {\color{blue}{  31}}&&  29.3 &  25.3 & {\color{blue}{  14}}&&  31.0 &  26.1 & {\color{blue}{  16}}\\
 .010&&  72.4&& 158.9 & 110.5 & {\color{blue}{  30}}&& 188.0 & 115.3 & {\color{blue}{  39}}&&  84.4 &  68.7 & {\color{blue}{  19}}&&  88.3 &  70.5 & {\color{blue}{  20}}\\
 .005&& 224.2&& 502.7 & 302.4 & {\color{blue}{  40}}&& 596.0 & 312.6 & {\color{blue}{  48}}&& 251.3 & 187.1 & {\color{blue}{  26}}&& 261.5 & 189.8 & {\color{blue}{  27}}\\[4pt]
\multicolumn{19}{c}{({\bf{c}}) Sparsity $k=100$}\\[8pt]
 .100&&   2.2&&   4.7 &   4.2 & {\color{blue}{  12}}&&   8.3 &   5.3 & {\color{blue}{  37}}&&   2.3 &   2.1 & {\color{blue}{   9}}&&   2.9 &   2.3 & {\color{blue}{  22}}\\
 .050&&   7.3&&  12.5 &  10.0 & {\color{blue}{  20}}&&  18.1 &  12.0 & {\color{blue}{  34}}&&   5.8 &   5.4 & {\color{blue}{   8}}&&   6.7 &   5.7 & {\color{blue}{  15}}\\
 .020&&  28.7&&  43.6 &  33.1 & {\color{blue}{  24}}&&  54.9 &  37.7 & {\color{blue}{  31}}&&  23.0 &  19.9 & {\color{blue}{  14}}&&  24.4 &  20.6 & {\color{blue}{  15}}\\
 .010&&  68.8&& 118.0 &  84.1 & {\color{blue}{  29}}&& 140.5 &  93.4 & {\color{blue}{  34}}&&  59.4 &  51.5 & {\color{blue}{  13}}&&  62.3 &  52.9 & {\color{blue}{  15}}\\
 .005&& 233.1&& 367.5 & 221.5 & {\color{blue}{  40}}&& 447.5 & 234.4 & {\color{blue}{  48}}&& 180.1 & 140.6 & {\color{blue}{  22}}&& 189.9 & 143.0 & {\color{blue}{  25}}\\[4pt]
\multicolumn{19}{c}{({\bf{d}}) Sparsity $k=50$, noise level 1\%}\\[8pt]
 .100&&   1.8&&   3.9 &   3.3 & {\color{blue}{  13}}&&   5.6 &   4.2 & {\color{blue}{  26}}&&   1.6 &   1.5 & {\color{blue}{   4}}&&   2.0 &   1.6 & {\color{blue}{  16}}\\
 .050&&   4.6&&   9.3 &   7.7 & {\color{blue}{  17}}&&  13.1 &   9.3 & {\color{blue}{  29}}&&   3.8 &   3.6 & {\color{blue}{   4}}&&   4.5 &   3.9 & {\color{blue}{  13}}\\
 .020&&  15.2&&  30.8 &  23.7 & {\color{blue}{  23}}&&  42.3 &  27.1 & {\color{blue}{  36}}&&  13.0 &  11.4 & {\color{blue}{  12}}&&  15.1 &  12.1 & {\color{blue}{  20}}\\
 .010&&  32.4&&  78.5 &  55.2 & {\color{blue}{  30}}&& 109.6 &  62.3 & {\color{blue}{  43}}&&  33.2 &  27.9 & {\color{blue}{  16}}&&  39.2 &  28.7 & {\color{blue}{  27}}\\
 .005&&  66.5&& 219.2 & 131.4 & {\color{blue}{  40}}&& 334.7 & 139.7 & {\color{blue}{  58}}&&  90.2 &  68.6 & {\color{blue}{  24}}&& 115.8 &  70.1 & {\color{blue}{  39}}\\[4pt]
\multicolumn{19}{c}{({\bf{e}}) Sparsity $k=50$, noise level 5\%}\\[8pt]
 .100&&   1.5&&   3.4 &   3.0 & {\color{blue}{  10}}&&   4.5 &   3.7 & {\color{blue}{  18}}&&   1.4 &   1.3 & {\color{blue}{   8}}&&   1.7 &   1.4 & {\color{blue}{  13}}\\
 .050&&   3.3&&   7.9 &   6.7 & {\color{blue}{  16}}&&  10.9 &   8.2 & {\color{blue}{  25}}&&   3.2 &   2.8 & {\color{blue}{  13}}&&   3.8 &   3.1 & {\color{blue}{  18}}\\
 .020&&   9.5&&  24.7 &  19.1 & {\color{blue}{  23}}&&  34.3 &  22.2 & {\color{blue}{  35}}&&  10.0 &   8.3 & {\color{blue}{  17}}&&  12.0 &   8.8 & {\color{blue}{  27}}\\
 .010&&  19.5&&  61.7 &  44.6 & {\color{blue}{  28}}&&  85.9 &  50.9 & {\color{blue}{  41}}&&  24.5 &  19.2 & {\color{blue}{  22}}&&  29.1 &  20.2 & {\color{blue}{  31}}\\
 .005&&  33.2&& 154.5 &  98.1 & {\color{blue}{  37}}&& 263.1 & 107.7 & {\color{blue}{  59}}&&  60.5 &  43.1 & {\color{blue}{  29}}&&  84.7 &  44.8 & {\color{blue}{  47}}\\[4pt]
\multicolumn{19}{c}{({\bf{f}}) Sparsity $k=50$, noise level 10\%}\\
\end{tabular}
\caption{Comparison between \spgl\ and root finding with strict
  tolerance levels using the \spg{} and hybrid method. The columns
  within the root finding and Lasso blocks are respectively the
  runtime in seconds of the \spg{} and hybrid method, and (in blue) the
  reduction in runtime in percent
  of the hybrid method compared to the \spg{} method.}\label{Table:Correlated}
\end{table}

For the first set of experiments we use the highly coherent matrices
described in Section~\ref{Sec:CoherentA}. As before we create a
$k$-sparse vector $x_0$ with non-zero entries sampled from different
distributions, and set $b = Ax_0 + v$, where the entries in $v$ are
zero in the noiseless case, and sampled i.i.d.~from the normal
distribution and scaled to the desired noise level otherwise. For the
noiseless results in Tables~\ref{Table:Correlated}(a)--(c) we run ten
instances for each of the three distributions and report the average
run time over all thirty runs. The percentage time reduction is
computed based on the total runtime and matches the percentage
obtained for each of the three signal classes independently. For the
root-finding columns we solve \eqref{Eq:BPDN} with $\sigma =
0.01\norm{b}_2$ and optimality tolerance levels of $10^{-4}$ and
$10^{-6}$. For the Lasso columns we solve \eqref{Eq:Lasso} on
equivalent problems with $\tau$ set to the value obtained using the
root-finding procedure. The results in
Tables~\ref{Table:Correlated}(d)--(f) apply to noisy problems where
$\norm{v}_2$ is scaled to the given percentage of $\norm{Ax_0}_2$, and
$\sigma$ is set accordingly. For these experiments we only consider
sparse $x_0$ with random $\pm 1$
entries. Table~\ref{Table:CoherentAccuracy} summarizes the total
runtime for the different solvers along with percentage of solutions
that have a relative duality gap within the given ranges.

\begin{table}
\small
\begin{tabular}{lllrrrrrrr}
\hline
&&&&\multicolumn{6}{c}{Relative duality gap}\\
\cline{5-10}
\\[-10pt]
Type & Method & Tol. & \multicolumn{1}{l}{Time} & $\leq 10^{-6}$& $10^{(-6,-5]}$& $10^{(-5,-4]}$& $10^{(-4,-3]}$& $10^{(-3,-2]}$& $> 10^{-2}$\\
\hline
BP$_{\sigma}$ & \spgl{}   & $10^{-6}$&  7h25&   0.2&   1.5&   4.6&    15&    49&    29\\
              & \spg{}    & $10^{-6}$& 17h29&    39&    33&    28&   0.2& --& --\\
              & hybrid    & $10^{-6}$&  9h53&   100& --& --& --& --& --\\
              & \spg{}    & $10^{-4}$& 14h02&   0.2&   0.6&    99&   0.5& --& --\\
              & hybrid    & $10^{-4}$&  9h18&   0.2&   0.8&    99& --& --& --\\
LS$_{\tau}$   & \spg{}    & $10^{-6}$&  7h55&    36&    33&    30&   0.5& --& --\\
              & hybrid    & $10^{-6}$&  6h02&   100& --& --& --& --& --\\
              & \spg{}    & $10^{-4}$&  7h20& --&   1.2&    98&   0.5& --& --\\
              & hybrid    & $10^{-4}$&  5h54& --&   1.2&    99& --& --& --\\
\hline
\end{tabular}
\caption{Total runtime for the coherent problems with different
  methods and optimality tolerances, along with the percentage of
  instances that attain a relative duality gap  in the given intervals at the
  final iteration. The reduction in runtime for the successive
  \spg{}-hybrid pairs are 43, 34, 24, and 20\%,
  respectively.}\label{Table:CoherentAccuracy}
\end{table}

The first thing to note from the results in
Table~\ref{Table:Correlated} is that the problems generated with lower
$\gamma$ values are indeed more difficult to solve for both the
\spg{} and the hybrid method. Compared to the \spg{} method, the
hybrid method reduces the average runtime for nearly all problems, and
does so by a percentage that increases as the problems get
harder. From Table~\ref{Table:CoherentAccuracy} we see that the hybrid
method with optimality tolerances of $10^{-4}$ and $10^{-6}$ reduces
the total runtime respectively by 34\% and 43\% for the basis-pursuit
problems, and 20\% and 24\% for the Lasso problems. The larger
relative reduction in runtime for basis pursuit is due to the use of
warm starting in the root-finding procedure, which removes a
substantial number of iterations that would otherwise be identical for
the hybrid and \spg{} methods. Despite the improvements, the hybrid
method still has a larger runtime than \spgl\ on most
problems. However, from Table~\ref{Table:CoherentAccuracy} we see that
\spgl\ does not even reach a relative duality gap of $10^{-3}$ for
nearly 80\% of the problems, as a result of the relaxed stopping
criteria.  Tightening these criteria, as done in what we label the
\spg{} method, increases the number of solutions that attain the
desired optimality tolerance. Nevertheless, the \spg{} method still
fails to reach an optimality of $10^{-6}$ for some 60\% of the
problems. Finally, we see that the hybrid method not only improves the
runtime of the \spg{} method, but also manages to reach the
requested optimality on all problems from
Table~\ref{Table:Correlated}.

\subsubsection{Sparco test problems}

Sparco \cite{BER2009FHHa} provides a standard collection of test
problems for compressed sensing and sparse recovery. The problems in
Sparco are of the form $b = Ax+v$, where $A$ is represented as a
linear operator rather than an explicit matrix. After excluding
problems that are too easy to solve or require access to third-party
software, we obtain the problem selection listed in
Table~\ref{Table:SparcoProblems}. For some problems we scale the
original $b$ to avoid a very small objective value at the solution,
which causes the duality gap relative to $\max(f(x),1)$ to be satisfied
more easily. The table also lists the one-norm of the solutions found
when solving with  $\sigma = 0.01\norm{b}_2$ and $\sigma =
0.001\norm{b}_2$, respectively, for the scaled $b$.

We run the \spg{} and hybrid methods with optimality tolerances
ranging from $10^{-1}$ down to $10^{-4}$. Beyond that, some of the
problems simply took too long to finish. For \spgl\ we use optimality
tolerance values set to $10^{-6}$ and $10^{-9}$. By comparison these
may seem excessively small, and we certainly do not expect the
relative duality gap to reach these levels. Instead, we choose the
small values to help control the other stopping criteria, such as the
relative change in the objective value, which are parameterized using
the same tolerance parameter. The results of the experiments with the
two choices of $\sigma$, are summarized in Tables~\ref{Table:Sparco1}
and~\ref{Table:Sparco2}. The hybrid method reduces the runtime of the
\spg{} method in 42 out of the 56 settings, often considerably
so. For a tolerance level of $10^{-4}$ the hybrid method consistently
outperforms the \spg{} method with an average time reduction of
38\%. 
The required optimality level is reached on all problems except for
problem 903 with the smaller $\sigma$ and optimality tolerance
$10^{-4}$. For this problem the \spg{} method stops with a relative
duality gap of $2\cdot 10^{-4}$ following a line-search error.  The
runtime for \spgl\ with optimality tolerance $10^{-6}$ is very low
overall, but comes at the cost of a rather large relative duality gap
at the solution. Lowering the tolerance to $10^{-9}$ reduces the gap,
but also leads to a considerable increase in runtime. In either case
the number of root-finding iterations can be very large, especially if
the target value of $\tau$ is exceeded and gradual reduction
follows. The lowest relative duality gap reached by \spgl\ over all
problems in Tables~\ref{Table:Sparco1} and~\ref{Table:Sparco2} is
$4\cdot 10^{-3}$.
The varying optimality levels make it difficult to compare results, so
of special interest are problem instances where \spgl\ simultaneously
has a lower runtime and relative duality gap with either the
\spg{} or hybrid method, or vice versa. From the tables we see that \spgl\
outperforms the \spg{} method on both instances of problem
702. For problem 401 in Table~\ref{Table:Sparco1}, \spgl\ with an
optimality tolerance of $10^{-6}$ is better, but aside from this
problem, \spgl{} consistently has the lowest runtime, but also the largest
duality gap. The \spg{} method with more stringent root-finding
iterations dominates \spgl\ with a tolerance level of $10^{-9}$ on all
remaining problems aside from the instance of problem 701 in
Table~\ref{Table:Sparco1}.  As we saw earlier, the hybrid method
performs especially well when the desired relative duality gap is
small. Nevertheless, even for the large duality gaps in question it
still dominates \spgl\ on nine out of the fourteen problem instances and is
dominated on only one.

\begin{table}[t]
\centering
\begin{tabular}{lrrrrrrr}
\hline
Problem & ID & $m$ & $n$ & $\norm{b}_2$ & scale & $\norm{x^*_{1e-2}}_1$  & $\norm{x^*_{1e-3}}_1$\\
\hline
blurrycam & 701 & 65536 & 65536 & 1.3e+2 & 100 & 2.8e+5 & 9.1e+5\\
blurspike & 702 & 16384 & 16384 & 2.2e+0 & 100 & 2.5e+4 & 3.4e+4\\
soccer1 & 601 & 3200 & 4096 & 5.5e+4 & 1 & 7.9e+1 & 3.1e+2\\
spiketrn & 903 & 1024 & 1024 & 5.7e+1 & 100 & 1.3e+3 & 1.3e+3\\
yinyang & 603 & 1024 & 4096 & 2.5e+1 & 100 & 2.5e+4 & 2.6e+4\\
srcsep1 & 401 & 29166 & 57344 & 2.2e+1 & 1 & 1.0e+3 & 1.0e+3\\
srcsep2 & 402 & 29166 & 86016 & 2.3e+1 & 1 & 1.1e+3 & 1.1e+3\\
\hline
\end{tabular}
\caption{Selected sparco problems.}\label{Table:SparcoProblems}
\end{table}

\begin{table}
\centering\small\setlength{\tabcolsep}{4.5pt}
\begin{tabular}{lrrcrrcrrcrrcrr}
&\multicolumn{2}{c}{\spgl{}}&&\spg{}&hybr.&&\spg{}&hybr.&&\spg{}&hybr.&&\spg{}&hybr.\\
&$10^{-6}$ & $10^{-9}$ &&\multicolumn{2}{c}{Tol = $10^{-1}$}&&\multicolumn{2}{c}{Tol = $10^{-2}$}&&\multicolumn{2}{c}{Tol = $10^{-3}$}&&\multicolumn{2}{c}{Tol = $10^{-4}$}\\
\cline{2-3}\cline{5-6}\cline{8-9}\cline{11-12}\cline{14-15}
Runtime (s) &       18.0 &       32.2 &  &       24.3 &       24.9 &  &       44.6 &       54.0 &  &       84.1 &       75.9 &  &       91.1 &       66.7\\
Rel.gap & 3e-1 & 2e-2 &  & 9e-2 & 9e-2 &  & 9e-3 & 9e-3 &  & 9e-4 & 1e-3 &  & 1e-4 & 8e-5\\
Outer iterations &        144 &        246 &  &         12 &         13 &  &         12 &         12 &  &         12 &         10 &  &         10 &         10\\
\multicolumn{3}{l}{\color{blue}{\raisebox{1pt}{\footnotesize$\blacktriangleright$} Problem 701 -- blurrycam}}&&&{\color{blue}{-2.6}}&&&{\color{blue}{-21.2}}&&&{\color{blue}{9.7}}&&&{\color{blue}{26.8}}\\
\\[-7pt]
Runtime (s) &        9.7 &       16.7 &  &       17.5 &       18.2 &  &       27.8 &       28.0 &  &       38.9 &       38.2 &  &       38.2 &       35.4\\
Rel.gap & 6e-1 & 3e-2 &  & 9e-2 & 1e-1 &  & 9e-3 & 1e-2 &  & 1e-3 & 1e-3 &  & 1e-4 & 1e-4\\
Outer iterations &        191 &        296 &  &         11 &         11 &  &         10 &         10 &  &         10 &         10 &  &          8 &          8\\
\multicolumn{3}{l}{\color{blue}{\raisebox{1pt}{\footnotesize$\blacktriangleright$} Problem 702 -- blurspike}}&&&{\color{blue}{-3.9}}&&&{\color{blue}{-0.5}}&&&{\color{blue}{1.8}}&&&{\color{blue}{7.2}}\\
\\[-7pt]
Runtime (s) &        133 &        206 &  &       24.6 &       11.5 &  &       28.0 &       20.1 &  &       34.0 &       17.7 &  &       38.0 &       17.2\\
Rel.gap & 1.5 & 6e-2 &  & 5e-2 & 5e-2 &  & 7e-3 & 6e-3 &  & 1e-3 & 4e-4 &  & 1e-5 & 2e-5\\
Outer iterations &       1274 &       1945 &  &         13 &         11 &  &          9 &          9 &  &         12 &          9 &  &          8 &          9\\
\multicolumn{3}{l}{\color{blue}{\raisebox{1pt}{\footnotesize$\blacktriangleright$} Problem 601 -- soccer1}}&&&{\color{blue}{53.2}}&&&{\color{blue}{28.2}}&&&{\color{blue}{47.8}}&&&{\color{blue}{54.7}}\\
\\[-7pt]
Runtime (s) &        2.8 &        6.6 &  &        3.2 &        2.9 &  &        5.6 &        3.8 &  &        7.6 &        5.6 &  &       11.2 &        5.3\\
Rel.gap & 3.6 & 1e-1 &  & 9e-2 & 1e-1 &  & 8e-3 & 8e-3 &  & 1e-3 & 8e-4 &  & 1e-4 & 1e-4\\
Outer iterations &        448 &       1120 &  &         17 &          8 &  &          6 &          6 &  &          6 &         10 &  &          5 &          5\\
\multicolumn{3}{l}{\color{blue}{\raisebox{1pt}{\footnotesize$\blacktriangleright$} Problem 903 -- spiketrn}}&&&{\color{blue}{9.0}}&&&{\color{blue}{31.3}}&&&{\color{blue}{27.2}}&&&{\color{blue}{53.0}}\\
\\[-7pt]
Runtime (s) &        1.8 &        2.7 &  &        2.1 &        2.1 &  &        2.5 &        2.9 &  &        2.9 &        3.1 &  &        4.2 &        4.0\\
Rel.gap & 3e-1 & 2e-2 &  & 8e-2 & 9e-2 &  & 7e-3 & 6e-3 &  & 8e-4 & 8e-4 &  & 7e-5 & 8e-5\\
Outer iterations &         38 &         60 &  &         10 &         11 &  &          7 &          8 &  &          7 &          7 &  &          7 &          7\\
\multicolumn{3}{l}{\color{blue}{\raisebox{1pt}{\footnotesize$\blacktriangleright$} Problem 603 -- yinyang}}&&&{\color{blue}{-2.1}}&&&{\color{blue}{-13.5}}&&&{\color{blue}{-4.5}}&&&{\color{blue}{3.7}}\\
\\[-7pt]
Runtime (s) &       79.0 &        418 &  &       88.1 &       75.6 &  &        299 &        331 &  &       1404 &       1047 &  &       6351 &       2447\\
Rel.gap & 5e-2 & 1e-1 &  & 7e-2 & 8e-2 &  & 9e-3 & 9e-3 &  & 1e-3 & 1e-3 &  & 1e-4 & 1e-4\\
Outer iterations &         37 &        213 &  &         10 &         16 &  &          9 &         20 &  &          8 &          9 &  &          9 &          8\\
\multicolumn{3}{l}{\color{blue}{\raisebox{1pt}{\footnotesize$\blacktriangleright$} Problem 401 -- srcsep1}}&&&{\color{blue}{14.2}}&&&{\color{blue}{-10.5}}&&&{\color{blue}{25.5}}&&&{\color{blue}{61.5}}\\
\\[-7pt]
Runtime (s) &        114 &        622 &  &        151 &        146 &  &        441 &        409 &  &       1114 &       1023 &  &       2447 &       1565\\
Rel.gap & 5e-1 & 1e-1 &  & 8e-2 & 7e-2 &  & 1e-2 & 9e-3 &  & 9e-4 & 1e-3 &  & 1e-4 & 6e-5\\
Outer iterations &         30 &        198 &  &         10 &         11 &  &          9 &          9 &  &          8 &          7 &  &          7 &          9\\
\multicolumn{3}{l}{\color{blue}{\raisebox{1pt}{\footnotesize$\blacktriangleright$} Problem 402 -- srcsep2}}&&&{\color{blue}{3.3}}&&&{\color{blue}{7.1}}&&&{\color{blue}{8.2}}&&&{\color{blue}{36.1}}\\
\end{tabular}
\caption{Sparco problems with $\sigma=0.01\norm{b}_2$, (a) runtime in seconds, (b)
  relative duality gap, (c) number of root-finding
  iterations. The percentage reduction in time of the hybrid method
  over \spg{} is given in blue next to the problem
  index.}\label{Table:Sparco1} 
\end{table}

\begin{table}
\centering\small\setlength{\tabcolsep}{4.5pt}
\begin{tabular}{lrrcrrcrrcrrcrr}
&\multicolumn{2}{c}{\spgl{}}&&\spg{}&hybr.&&\spg{}&hybr.&&\spg{}&hybr.&&\spg{}&hybr.\\
&$10^{-6}$ & $10^{-9}$ &&\multicolumn{2}{c}{Tol = $10^{-1}$}&&\multicolumn{2}{c}{Tol = $10^{-2}$}&&\multicolumn{2}{c}{Tol = $10^{-3}$}&&\multicolumn{2}{c}{Tol = $10^{-4}$}\\
\cline{2-3}\cline{5-6}\cline{8-9}\cline{11-12}\cline{14-15}
Runtime (s) &        121 &        200 &  &        184 &        195 &  &        272 &        300 &  &        316 &        342 &  &        430 &        333\\
Rel.gap & 2.1 & 1e-1 &  & 7e-2 & 7e-2 &  & 9e-3 & 6e-3 &  & 1e-3 & 8e-4 &  & 9e-5 & 9e-5\\
Outer iterations &        485 &        640 &  &         13 &         13 &  &         12 &         12 &  &         11 &         12 &  &         11 &         11\\
\multicolumn{3}{l}{\color{blue}{\raisebox{1pt}{\footnotesize$\blacktriangleright$} Problem 701 -- blurrycam}}&&&{\color{blue}{-6.4}}&&&{\color{blue}{-10.4}}&&&{\color{blue}{-8.2}}&&&{\color{blue}{22.5}}\\
\\[-7pt]
Runtime (s) &       28.8 &       34.9 &  &       35.4 &       32.7 &  &       45.8 &       49.3 &  &       56.1 &       59.8 &  &       74.4 &       68.1\\
Rel.gap & 4e-1 & 7e-2 &  & 9e-2 & 9e-2 &  & 9e-3 & 1e-2 &  & 9e-4 & 9e-4 &  & 9e-5 & 9e-5\\
Outer iterations &        127 &        151 &  &         13 &         13 &  &         11 &         11 &  &         10 &         10 &  &         10 &         10\\
\multicolumn{3}{l}{\color{blue}{\raisebox{1pt}{\footnotesize$\blacktriangleright$} Problem 702 -- blurspike}}&&&{\color{blue}{7.5}}&&&{\color{blue}{-7.7}}&&&{\color{blue}{-6.5}}&&&{\color{blue}{8.5}}\\
\\[-7pt]
Runtime (s) &        278 &        368 &  &       84.3 &       47.6 &  &        155 &       56.4 &  &        203 &       53.2 &  &        257 &       47.8\\
Rel.gap & 1.1 & 7e-2 &  & 8e-2 & 7e-2 &  & 9e-3 & 8e-3 &  & 1e-3 & 7e-4 &  & 1e-4 & 7e-5\\
Outer iterations &       2059 &       2707 &  &         19 &         10 &  &         11 &         10 &  &          9 &          9 &  &          9 &         10\\
\multicolumn{3}{l}{\color{blue}{\raisebox{1pt}{\footnotesize$\blacktriangleright$} Problem 601 -- soccer1}}&&&{\color{blue}{43.5}}&&&{\color{blue}{63.6}}&&&{\color{blue}{73.8}}&&&{\color{blue}{81.4}}\\
\\[-7pt]
Runtime (s) &        4.8 &        9.4 &  &        6.4 &        5.0 &  &        9.2 &        6.2 &  &       11.7 &        6.6 &  &       13.5 &        6.6\\
Rel.gap & 22.5 & 1.3 &  & 5e-2 & 9e-2 &  & 1e-2 & 1e-2 &  & 1e-3 & 9e-4 &  & 2e-4 & 1e-4\\
Outer iterations &        700 &       1253 &  &          6 &          6 &  &          8 &          8 &  &          6 &          8 &  &          5 &          5\\
\multicolumn{3}{l}{\color{blue}{\raisebox{1pt}{\footnotesize$\blacktriangleright$} Problem 903 -- spiketrn}}&&&{\color{blue}{21.9}}&&&{\color{blue}{33.1}}&&&{\color{blue}{43.7}}&&&{\color{blue}{51.0}}\\
\\[-7pt]
Runtime (s) &        9.0 &       19.0 &  &       17.8 &       21.7 &  &       37.4 &       26.8 &  &       55.1 &       29.4 &  &       71.8 &       35.8\\
Rel.gap & 1.2 & 9e-1 &  & 1e-1 & 8e-2 &  & 1e-2 & 7e-3 &  & 9e-4 & 1e-3 &  & 1e-4 & 1e-4\\
Outer iterations &         90 &        190 &  &         16 &         15 &  &          8 &         10 &  &          8 &          8 &  &          8 &          8\\
\multicolumn{3}{l}{\color{blue}{\raisebox{1pt}{\footnotesize$\blacktriangleright$} Problem 603 -- yinyang}}&&&{\color{blue}{-21.4}}&&&{\color{blue}{28.4}}&&&{\color{blue}{46.7}}&&&{\color{blue}{50.1}}\\
\\[-7pt]
Runtime (s) & \color{gray}{       241} &       3309 &  &        108 &       85.0 &  &        394 &        360 &  &       2803 &       2253 &  &      21829 &      11193\\
Rel.gap & \color{gray}{3e-3} & 4e-3 &  & 8e-2 & 6e-2 &  & 1e-2 & 9e-3 &  & 1e-3 & 1e-3 &  & 1e-4 & 1e-4\\
Outer iterations & \color{gray}{        32} &        264 &  &         46 &         21 &  &         13 &         12 &  &         20 &         11 &  &         10 &         10\\
\multicolumn{3}{l}{\color{blue}{\raisebox{1pt}{\footnotesize$\blacktriangleright$} Problem 401 -- srcsep1}}&&&{\color{blue}{21.6}}&&&{\color{blue}{8.5}}&&&{\color{blue}{19.6}}&&&{\color{blue}{48.7}}\\
\\[-7pt]
Runtime (s) & \color{gray}{       272} &       5999 &  &        191 &        178 &  &        518 &        477 &  &       2041 &       1787 &  &       8631 &       6486\\
Rel.gap & \color{gray}{6e-3} & 2e-2 &  & 6e-2 & 8e-2 &  & 8e-3 & 8e-3 &  & 1e-3 & 9e-4 &  & 1e-4 & 9e-5\\
Outer iterations & \color{gray}{        18} &        324 &  &        240 &        192 &  &         11 &         14 &  &          9 &         10 &  &         10 &         10\\
\multicolumn{3}{l}{\color{blue}{\raisebox{1pt}{\footnotesize$\blacktriangleright$} Problem 402 -- srcsep2}}&&&{\color{blue}{6.5}}&&&{\color{blue}{8.0}}&&&{\color{blue}{12.5}}&&&{\color{blue}{24.9}}\\
\end{tabular}
\caption{Sparco problems with $\sigma = 0.001\norm{b}_2$, (a) runtime in seconds, (b)
  relative duality gap, (c) number of root-finding
  iterations. The percentage reduction in time of the hybrid method
  over \spg{} is given in blue next to the problem index. Entries marked
  in gray indicate a solution that is not a root.}\label{Table:Sparco2}
\end{table}

\subsubsection{Primal-dual gap}

We now consider the formulation
\begin{equation}\label{Eq:Abmu}
\minimize{x}\quad \half\norm{Ax-b}_2^2 + \sfrac{\mu}{2}\norm{x}_2^2
\quad\st\quad x\in\mathcal{C},
\end{equation}
for $\mu > 0$. In Section~\ref{Sec:StoppingCriteria} we described two
different ways of deriving a dual formulation.  In the first approach
we augment $A$ and $b$ to account for the $\sfrac{\mu}{2}\norm{x}_2^2$
term and reduce the problem to the standard Lasso formulation. The
derivation of the dual for this formulation in
\cite{BER2008Fb,BER2010Fa} provides a way of generating a
dual-feasible point $(\bar{y},\bar{\lambda})$ from a primal-feasible
$x$ by choosing $\bar{y} = \bar{A}x - \bar{b}$ and solving a trivial
optimization problem for $\bar{\lambda}$. In the second approach we
deal with formulation \eqref{Eq:Abmu} directly and obtain a dual
problem parameterized in $(y,\lambda)$. As before we can choose $y$ to
be equal to the residual, now in terms of the original $A$ and $b$,
and remain with a non-trivial optimization problem for $\lambda$ that
is nevertheless easily solved using the algorithm described in
Section~\ref{Sec:StoppingCriteria}. We refer to the two derivations as
the augmented derivation and the optimized derivation. The term
`optimized' refers to the need to solve for $\lambda$, but more
importantly, to the fact that the dual objective generated from any $x$
using the optimized derivation is never smaller than that using the
augmented derivation, as shown in Section~\ref{Sec:StoppingCriteria}.

\begin{figure}[t!]
\centering
\begin{tabular}{cc}
\includegraphics[width=0.475\textwidth]{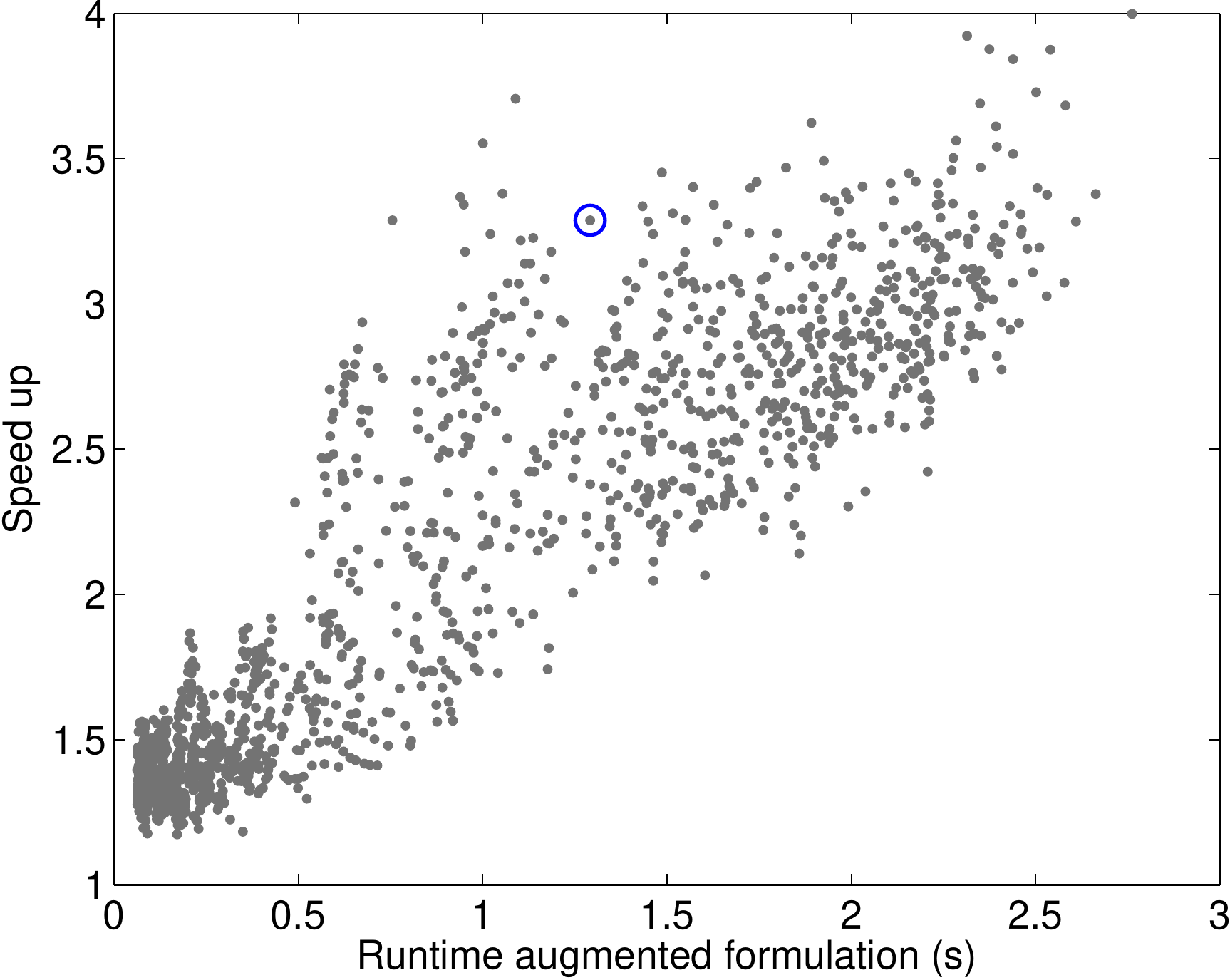}&
\includegraphics[width=0.475\textwidth]{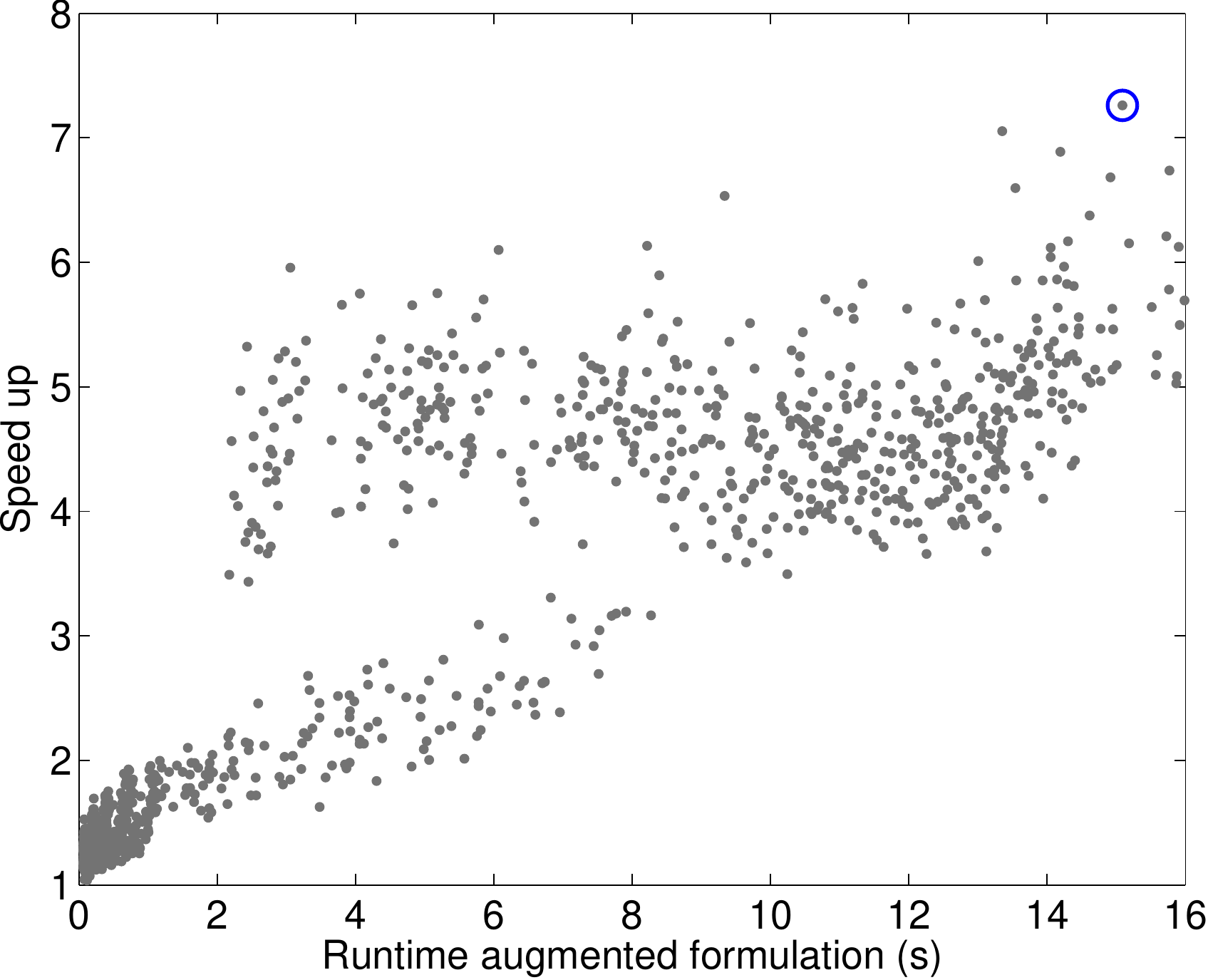}\\
({\bf{a}}) $\mu = 0.01$ & ({\bf{b}}) $\mu = 0.001$\\[3pt]
\includegraphics[width=0.475\textwidth]{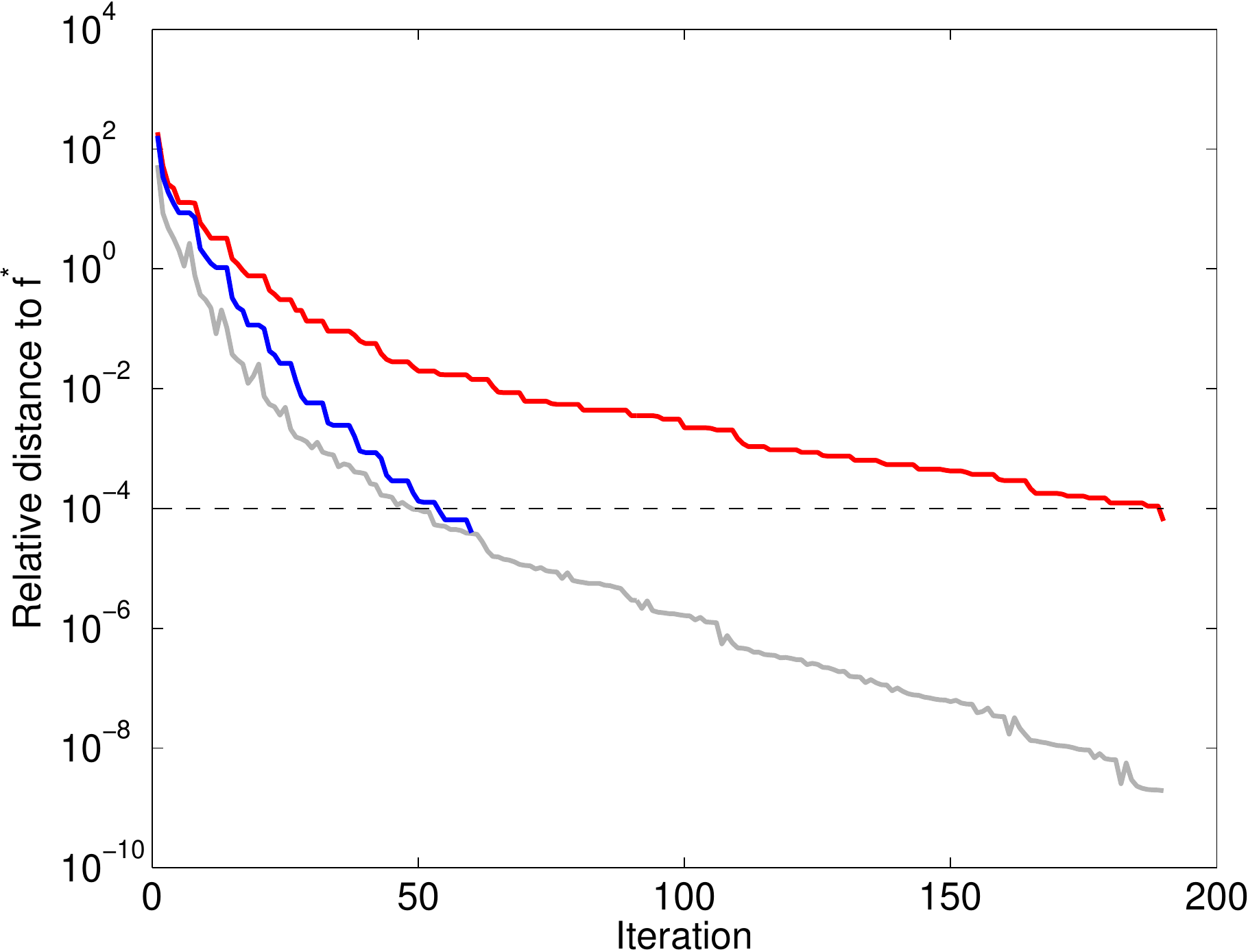}&
\includegraphics[width=0.475\textwidth]{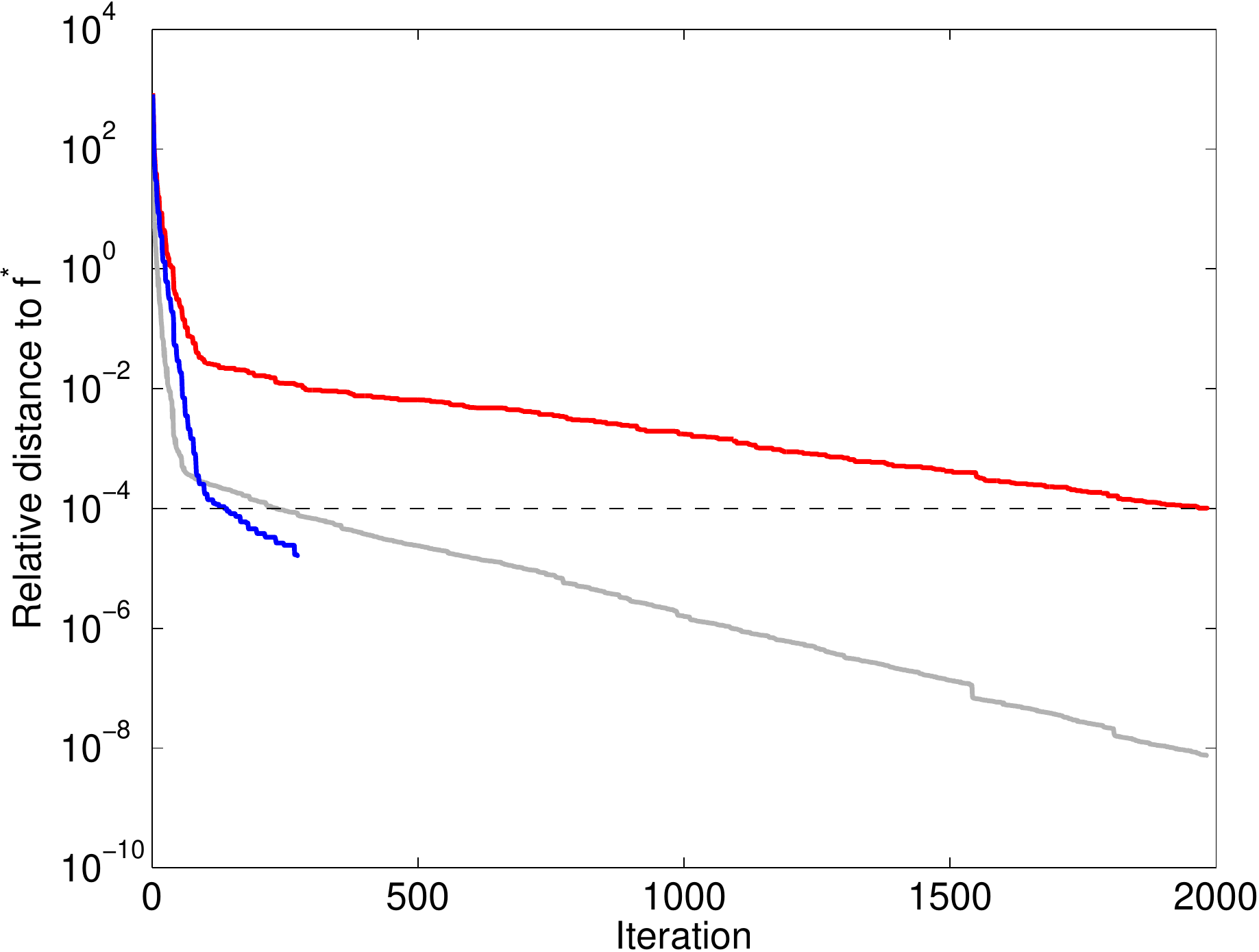}
\end{tabular}
\caption{Plots of (a,b) the time required to reach the desired
  relative optimality gap of $10^{-4}$ using the augmented
  formulation, and the corresponding speed up obtained using the
  optimized formulation. Each point indicates a single problem
  instance; and (c,d) the relative distance of the primal (gray) and
  dual objective (red for the augmented formulation and blue for the
  optimized formulation) to the optimal objective as a function of
  iterations for the two circled problem
  instances.}\label{Figure:Duality}
\end{figure}

To evaluate the practical difference between the two approaches we
generate a large number of randomized test problems of the form $b =
Ax_0 + v$, where $x_0$ are random vectors with sparsity levels ranging
from 50 to 350 in steps of 50 and on-support entries draw i.i.d.~from the
normal distribution. The $m\times n$ measurement matrices $A$ are
drawn i.i.d.~from $\mathcal{N}(0,1/m)$, with $m=1000$ and $n=2000$.
Finally, the additive noise vectors $v$ are drawn from the normal
distribution and scaled to have Euclidean norm equal to $0$, $0.01$,
$0.1$, and $1$ percent of that of the clean observation $Ax_0$. For
the problem formulation we set $\tau$ to $\norm{x_0}_1$ scaled by 0.7
through 1.2 in 0.1 increments, and range $\mu$ log-linearly from
$10^{-1}$ to $10^{-4}$ in four steps. As a result of the additive
$\sfrac{\mu}{2}\norm{x}_2$ term in the objective, the solutions are no
longer sparse. As a result the hybrid method tends to coincide with
the \spg{} method, and we therefore only consider the latter for
these experiments.

\begin{figure}[t!]
\centering
\begin{tabular}{cc}
\includegraphics[height=0.325\textwidth]{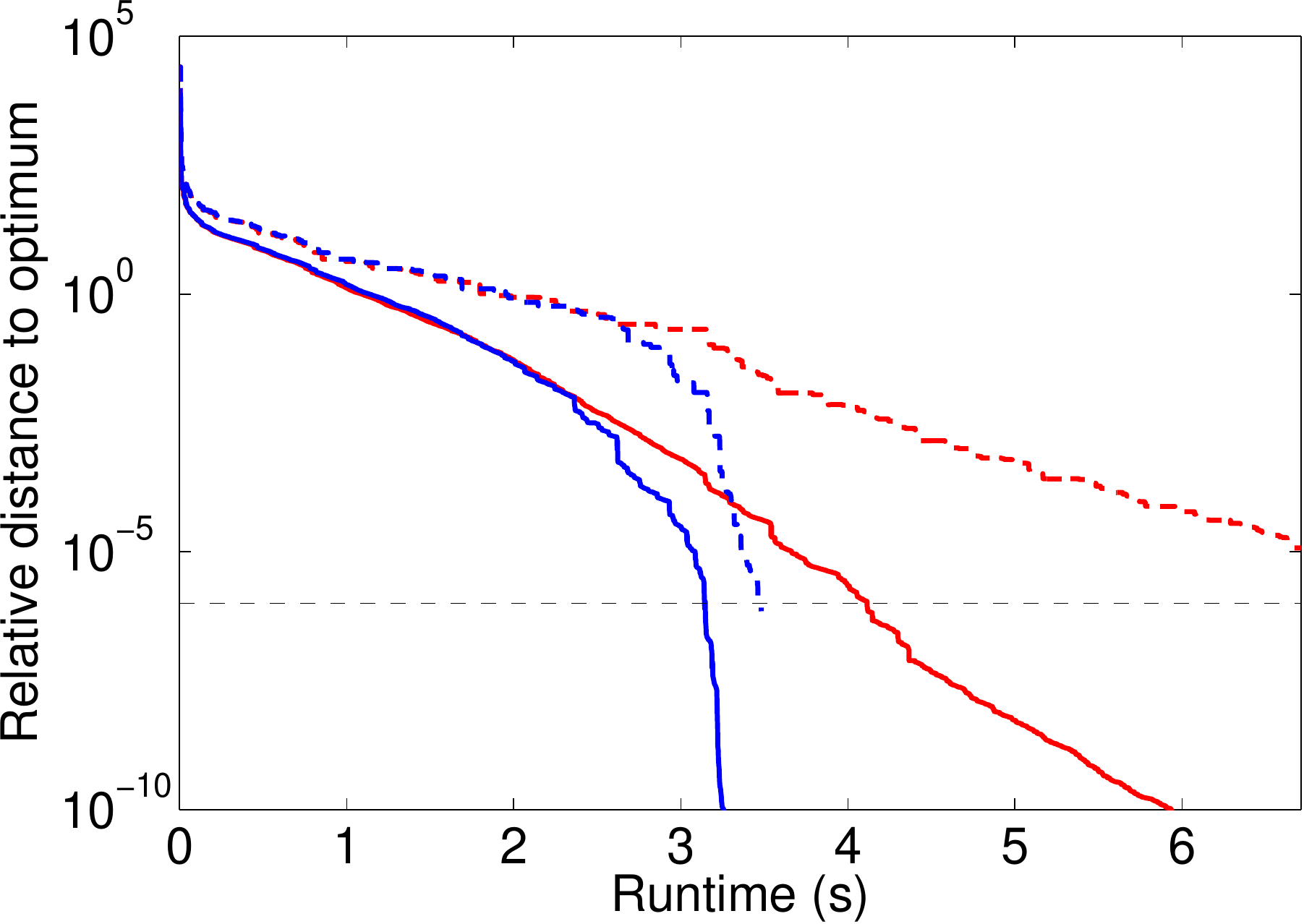}&
\includegraphics[height=0.325\textwidth]{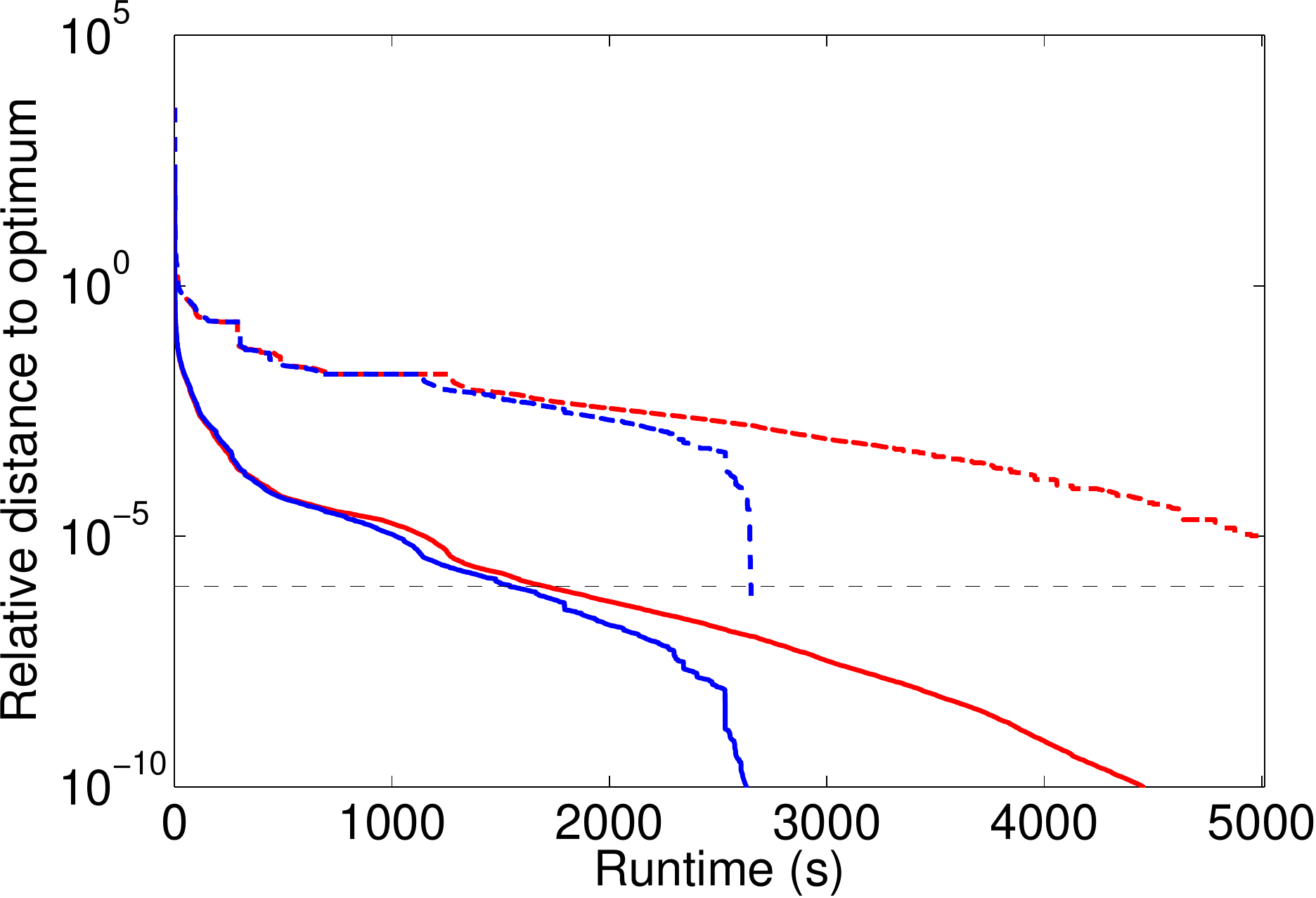}\\
({\bf{a}}) spiketrn  ($10^{-2}$) & ({\bf{b}}) srcsep2 ($10^{-2}$) \\[2pt]
\includegraphics[height=0.325\textwidth]{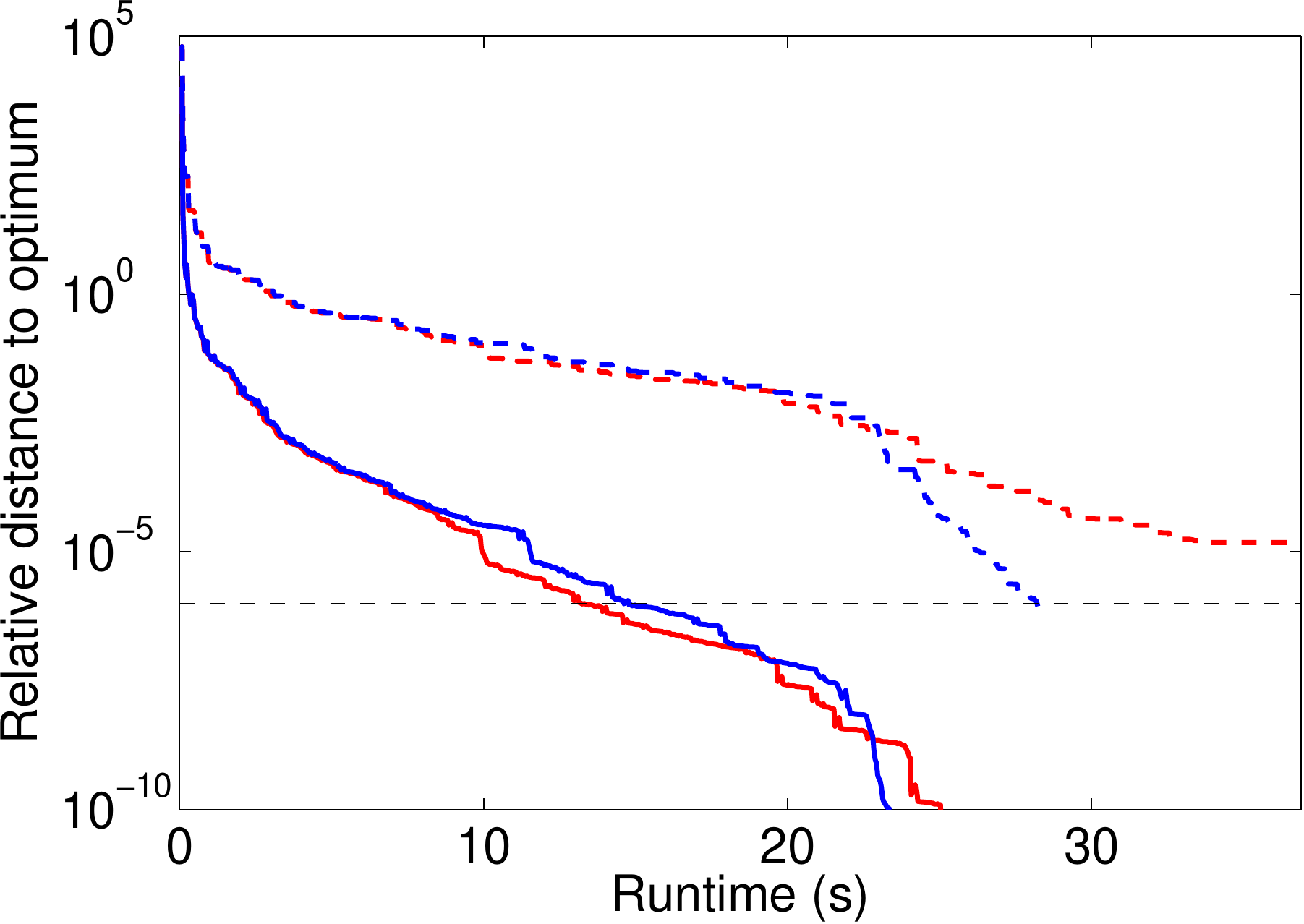}&
\includegraphics[height=0.325\textwidth]{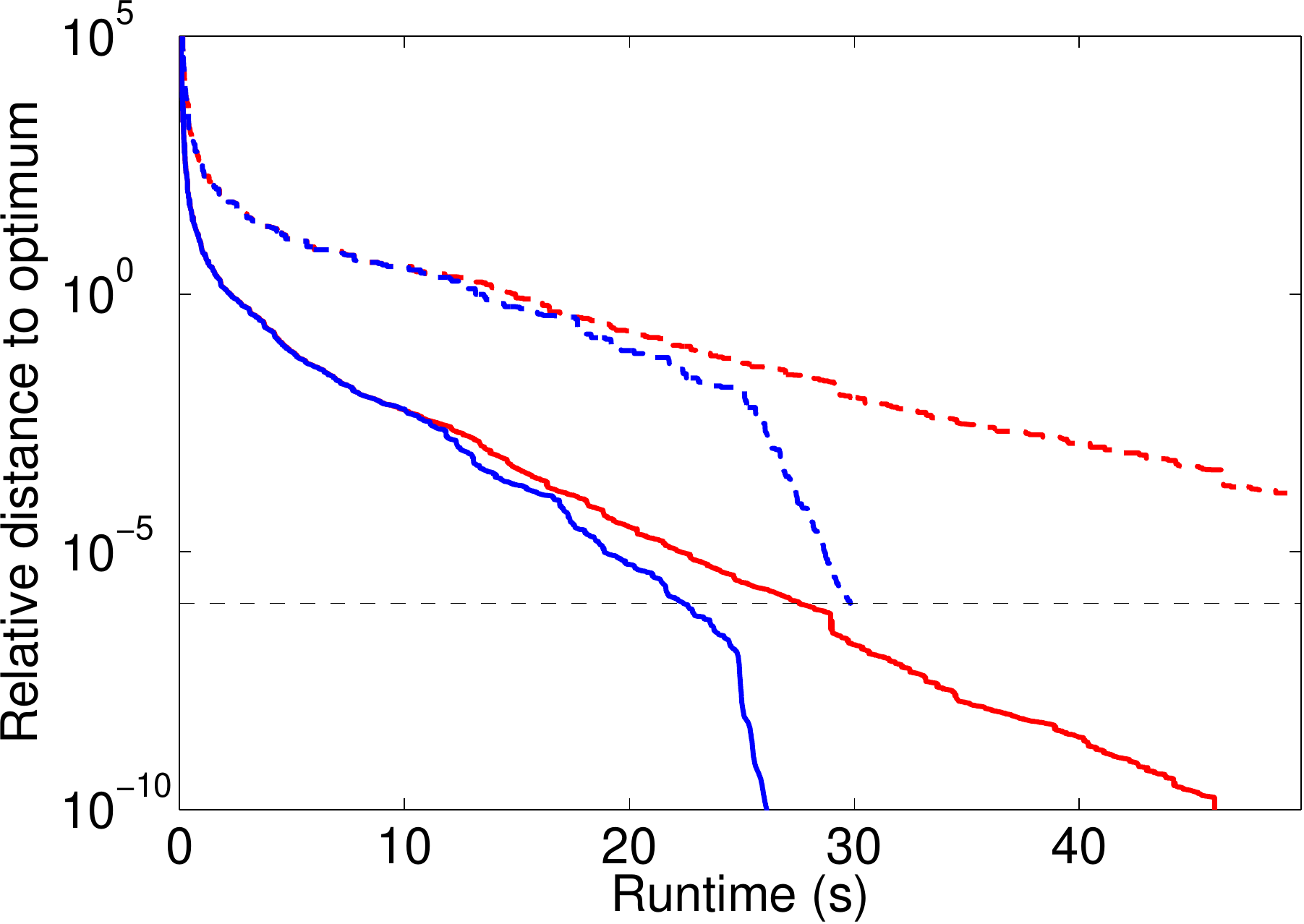}\\
({\bf{c}}) blurrycam ($10^{-2}$) & ({\bf{d}}) yinyang ($10^{-3}$)
\end{tabular}
\caption{Relative distance of primal (solid) and dual iterates
  (dashed) to the optimum
  as a function of time using the \spg{} (red) and hybrid method
  (blue) on four Sparco problems for fixed $\tau$ corresponding to the
  root for  $\sigma$ equal to the given multiple of $\norm{b}_2$.}\label{Fig:DualityGap2}
\end{figure}

For each of the settings we evaluate the time required by the
augmented and optimized formulations to reach a relative duality gap
of $10^{-4}$. Figures~\ref{Figure:Duality}(a,b) plots the speed up
obtained using the optimized formulation against the runtime of the
augmented formulation for two levels of $\mu$. Despite the slightly
more expensive evaluation of the dual, we see that the optimized
formulation is around 1.5 to 4 times faster for $\mu=0.01$, and up to
7 times faster for $\mu=0.001$. For $\mu=0.1$ (not shown in the plot)
the speedup ranges from 1.2 to 3, and for $\mu=0.001$ the speed up
exceeds 10 on many problem instances and reaches a maximum of around
30. 

We now take a closer look at the relative distance of the primal and
dual objective to the optimum for the two circled problems in
Figures~\ref{Figure:Duality}(c,d). The progress of the primal
objective over the iterations, indicated by the gray line, is the same
for both formulations. For the dual objective there is a marked
difference between the two. Notably, the augmented formulation
converges much slower than the optimized formulation, thereby
preventing the stopping criterion from being satisfied for many more
iterations.

We now take another look at the Sparco problems from
Tables~\ref{Table:Sparco1} and \ref{Table:Sparco2}. For each of the
settings we record the optimal $\tau$ and then run the hybrid solver
with a target optimality tolerance of $10^{-8}$ to obtain a
best-effort optimum (for some problems the line search failed before
reaching the desired tolerance). We then run the \spg{} and hybrid
solvers with a target accuracy of $10^{-5}$ and record the relative
distance of the primal and dual objective to the optimum at every
iteration. The results for four representative problems are plotted in
Figure~\ref{Fig:DualityGap2}. From the plots we see that the iterates
of the hybrid method initially coincide or otherwise closely follow
those of the \spg{} method. Once the hybrid method starts using
quasi-Newton iterates increasingly often we see a sharp decrease in
the relative distance to the optimum of the primal and dual
iterates. The iterates of the \spg{} method, by contrast, continue
to decrease very slowly. Indeed, of the fourteen problem settings, the
\spg{} method managed to solve only two to the desired level of
accuracy. Of the remaining problems, two reach the default iteration
limit of ten times the number of rows in $A$, while all other problems
fail with a line-search error. The hybrid method manages to solve all
problems except for problem 401 with multiplier $10^{-3}$. This
problem reached the iteration limit, but could otherwise be solved
successfully to a tolerance level of even $10^{-8}$.

As before, we see that the dual objective converges to the optimum
much slower than the primal, and unfortunately, there is no clear way
to extend the optimized dual formulation from
Section~\ref{Sec:StoppingCriteria} to the standard Lasso formulation
where $\mu=0$. Given that the satisfaction of the optimality condition
is controlled almost entirely by the dual objective value, it makes
sense to look at the potential speed up if the optimal objective value
was known and optimality was instead driven by the primal
objective. In Table~\ref{Table:DualityGap2} we provide this speed up
for the different Sparco problems with varying optimality tolerance
levels. Clearly, both the \spg{} and hybrid methods would benefit
from an improved dual, although the effect is less for the hybrid
method, due to the already fast convergence of the dual objective in
the final iterations.

\begin{table}
\centering\small
\setlength{\tabcolsep}{5.75pt}
\begin{tabular}{lrrrrrcrrrrr}
 & \multicolumn{5}{c}{($\sigma = 0.01$)} &&\multicolumn{5}{c}{($\sigma = 0.001$)}\\
 & $10^{-2}$ & $10^{-3}$ & $10^{-4}$ & $10^{-5}$ & $10^{-6}$ &&$10^{-2}$ & $10^{-3}$ & $10^{-4}$ & $10^{-5}$ & $10^{-6}$\\
\cline{2-6}\cline{8-12}
701 - blurrycam  &  9.15  &  5.93  &  3.79  &  \raisebox{2pt}{\tiny$>$}3.66  &  \raisebox{2pt}{\tiny$>$}2.79  &  &  3.55  &  2.73  &  \raisebox{2pt}{\tiny$>$}2.21  &  \raisebox{2pt}{\tiny$>$}1.94  &  \raisebox{2pt}{\tiny$>$}1.61 \\[1pt]
702 - blurspike  &  2.83  &  2.39  &  2.16  &  \raisebox{2pt}{\tiny$>$}1.93  &  \raisebox{2pt}{\tiny$>$}1.66  &  &  3.22  &  2.48  &  2.20  &  2.07  &  \raisebox{2pt}{\tiny$>$}1.81 \\[1pt]
601 - soccer1  &  1.80  &  1.25  &  1.07  &  1.02  &  1.02  &  &  2.59  &  2.84  &  2.88  &  \raisebox{2pt}{\tiny$>$}2.57  &  \raisebox{2pt}{\tiny$>$}2.09 \\[1pt]
903 - spiketrn  &  1.62  &  1.64  &  1.73  &  \raisebox{2pt}{\tiny$>$}1.83  &  \raisebox{2pt}{\tiny$>$}1.63  &  &  1.45  &  1.53  &  \raisebox{2pt}{\tiny$>$}1.64  &  \raisebox{2pt}{\tiny$>$}1.56  &  \raisebox{2pt}{\tiny$>$}1.44 \\[1pt]
603 - yinyang  &  2.00  &  1.70  &  1.62  &  1.64  &  1.52  &  &  3.54  &  3.02  &  \raisebox{2pt}{\tiny$>$}2.76  &  \raisebox{2pt}{\tiny$>$}2.24  &  \raisebox{2pt}{\tiny$>$}1.81 \\[1pt]
401 - srcsep1  &  68.59  &  40.72  &  26.71  &  \raisebox{2pt}{\tiny$>$}9.49  &  \raisebox{2pt}{\tiny$>$}4.15  &  &  118.54  &  75.78  &  \raisebox{2pt}{\tiny$>$}49.75  &  \raisebox{2pt}{\tiny$>$}12.30  &  \raisebox{2pt}{\tiny$>$}2.38 \\[1pt]
402 - srcsep2  &  20.44  &  15.96  &  9.81  &  \raisebox{2pt}{\tiny$>$}4.35  &  \raisebox{2pt}{\tiny$>$}2.94  &  &  54.46  &  56.53  &  22.74  &  9.88  &  \raisebox{2pt}{\tiny$>$}2.78 \\[3pt]
\multicolumn{12}{c}{({\bf{a}}) \spg{}}\\[10pt]
701 - blurrycam  &  9.45  &  5.51  &  3.20  &  2.27  &  1.93  &  &  3.96  &  3.62  &  2.73  &  1.90  &  1.48 \\[1pt]
702 - blurspike  &  2.74  &  2.20  &  1.81  &  1.56  &  1.45  &  &  2.72  &  2.30  &  1.82  &  1.61  &  1.42 \\[1pt]
601 - soccer1  &  1.24  &  1.19  &  1.08  &  1.01  &  1.02  &  &  1.35  &  1.11  &  1.03  &  1.02  &  1.02 \\[1pt]
903 - spiketrn  &  1.34  &  1.23  &  1.14  &  1.09  &  1.11  &  &  1.19  &  1.13  &  1.13  &  1.11  &  1.09 \\[1pt]
603 - yinyang  &  2.23  &  1.80  &  1.64  &  1.55  &  1.45  &  &  2.99  &  2.12  &  1.63  &  1.52  &  1.33 \\[1pt]
401 - srcsep1  &  38.12  &  24.83  &  11.83  &  3.72  &  1.81  &  &  119.99  &  65.28  &  46.18  &  \raisebox{2pt}{\tiny$>$}13.44  &  \raisebox{2pt}{\tiny$>$}2.61 \\[1pt]
402 - srcsep2  &  17.46  &  11.62  &  6.45  &  2.60  &  1.72  &  &  50.98  &  44.70  &  16.92  &  6.07  &  1.99 \\[3pt]
\multicolumn{12}{c}{({\bf{b}}) Hybrid}\\
\end{tabular}
\caption{Projected speed up when the optimal objective value is known
  and satisfaction of the optimality condition depends only on the
  primal objective value. We give a lower bound (indicated by the
  `\,{\tiny{$^{_>}$}}\!' sign) when the dual objective failed to reach the given
  optimality level, either because the maximum number of iterations
  was reached, or because a line-search error occurred.}\label{Table:DualityGap2}
\end{table}


\section{Conclusions}\label{Sec:Conclusions}
In this paper we have presented a hybrid algorithm for minimization of
quadratic functions over weighted one-norm balls. The method extends
the spectral projected gradient method with \hbox{\abbrv{L-BFGS}}
iterations applied to reparameterizations of the objective function
over active faces of the one-norm ball. For the decision of the
iteration type we introduce the self-projection cone of a face and
provide a complete characterization of this cone for weighted one-norm
balls. The reparameterization uses an implicit orthonormal basis for
the current face, and we provide an efficient algorithm for
matrix-vector multiplication with this basis and its transpose. Our
regular first-order iterations use backtracking line search of
projected gradient steps. In addition to this we investigate the use
of a trajectory line search over the entire projection curve of $x +
\alpha d$ with $\alpha \geq 0$. We show that this curve is piecewise
linear with at most $4n-1$ segments, and that a local or global
minimum of the objective along this curve can be determined
efficiently. Despite this, the computational cost was still found to
be high relative to projected backtracking line search, which showed
overall better performance.

As part of the numerical experiments we propose a challenging class of
test problems in which the columns of the $m\times n$ measurement
matrix $A$ are generated based on a random walk over the
$(m\!-\!1)$-sphere.  Based on extensive numerical experiments on these
and other test problems we showed that the hybrid method outperforms
the original spectral projected gradient methods on a large fraction
of the problems. Especially for medium to high accuracy solves and
more challenging problems the \spg{} method was found to either take
much more time to reach the desired level of accuracy, or fail
prematurely due to line-search problems. The current stopping
criterion of both methods relies on the generation of a dual feasible
point from the primal iterate to determine the relative optimality of
the iterate. From the experiments we found that the the primal
objective converges to the optimum much faster than the dual
objective, and that satisfaction of the stopping criterion therefore
depends entirely on the dual objective reaching the critical
threshold. The performance of both methods could therefore be improved
substantially given a better dual estimate.

In this paper we have studied the application of the hybrid method to
the Lasso problem. Other important problems that may benefit from the
approach but not discussed in this paper include box-constrained
optimization and minimization of quadratic functions over the simplex.


\appendix

\section{Proof of Theorem~\ref{Thm:MaxProjectionFaces}}\label{AppSec:Projection}

\begin{figure}[t!]
\centering
\begin{tabular}{cc}
\includegraphics[width=0.455\textwidth]{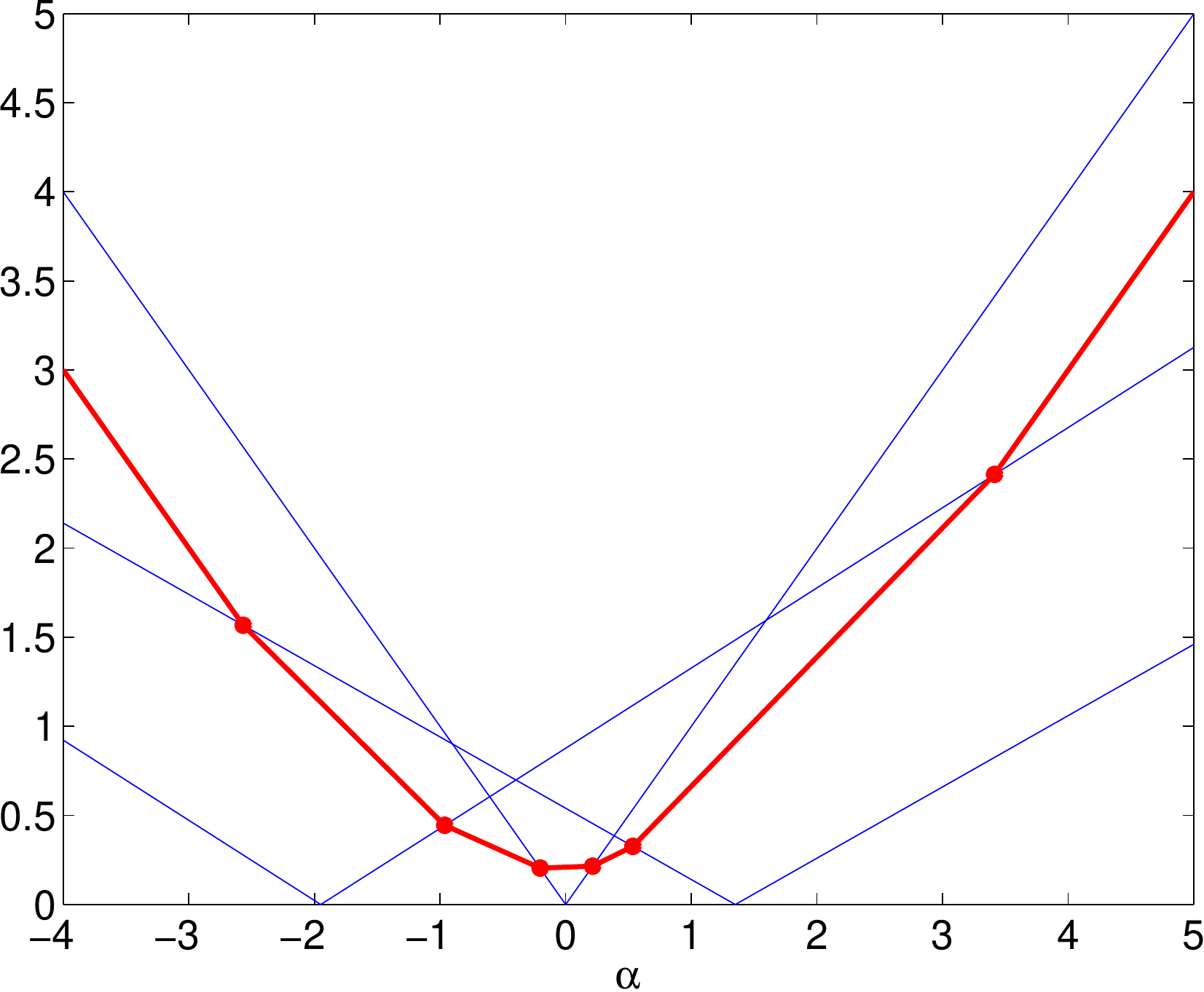}%
\begin{picture}(0,0)(0,0)
\put(-62,54){$\lambda(\alpha)$}
\end{picture}&
\includegraphics[width=0.455\textwidth]{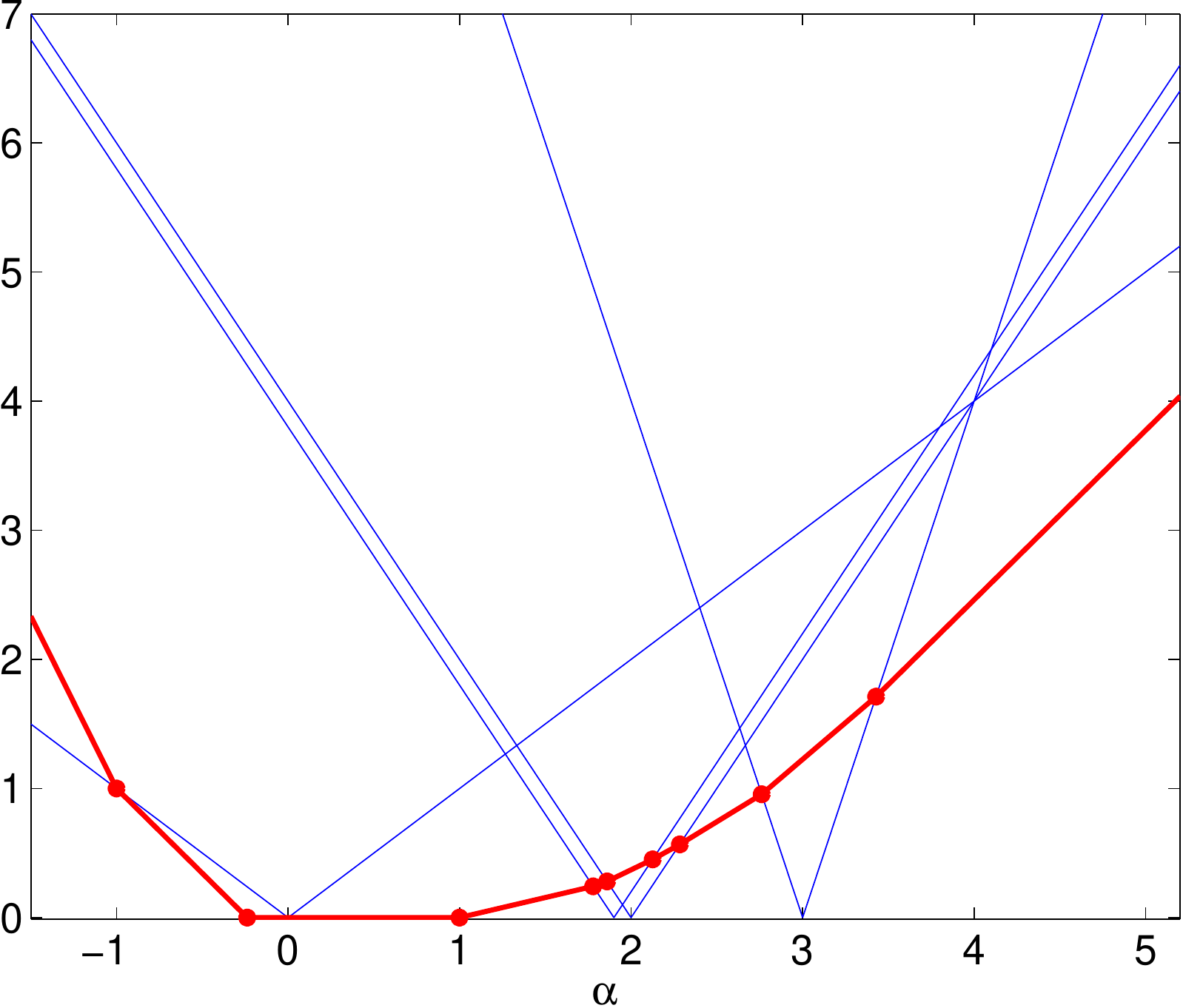}\\
({\bf{a}}) & ({\bf{b}})
\end{tabular}
\caption{Plot of the absolute value of $x_i + \alpha d_i$ (blue
  curves) and soft-thresholding parameter $\lambda(\alpha)$ (red) as a
  function of $\alpha$, with (a) a generic situation, and (b) a
  configuration for $n=4$ that attains the maximum of $4n-2 = 14$
  intersections (five of these appear outside the plotted
  region).}\label{Fig:ProjectionPath}
\end{figure}

In this section we study the combinatorial properties of the
projection of the line $x(\alpha) = o - \alpha d$, onto an
$n$-dimensional one-norm ball $\mathcal{C}_{w,1}$ of radius
$\tau$. Without loss of generality we can assume that the slopes $d_i$
satisfy $d_i \geq 0$, by changing the signs of $o_i$ and $d_i$ if
needed (this consistently reverses the signs of individual
coordinates, but does not otherwise affect the projection). When
plotting $v_i(\alpha) = \abs{o_{i} - \alpha d_{i}} / w_i$ as a
function of $\alpha$ we get horizontal lines when $d_i = 0$, and
v-shaped curves when $d_i > 0$, as illustrated in
Figure~\ref{Fig:ProjectionPath}(a). Also plotted in this figure is
$\lambda(\alpha)$: the soft-thresholding parameter $\lambda$ used to
project $x(\alpha)$ onto $\mathcal{C}_{w,1}$. For any $\alpha$ we can
see that the support $\mathcal{I}$ in the projection consists of all
indices $i$ for which $v_i(\alpha) > \lambda(\alpha)$. Over intervals of
$\alpha$ where the support remains fixed, we have
\[
\lambda'(\alpha) = \frac{\sum_{i\in\mathcal{I}}
  w_i r_i(\alpha)}{\sum_{j\in\mathcal{I}} w_i^2},\qquad
r_i(\alpha) = \begin{cases} -d_i &\ o_i - \alpha d_i > 0\\ \phantom{-}d_i
  & \ \mathrm{otherwise.}\end{cases}
\]
Changes in the support occur whenever $\lambda(\alpha)$ intersects
some $v_i(\alpha)$. Letting $\alpha'$ be a critical value where such
an intersection occurs we can consider the linear segments of
$\lambda(\alpha)$ that end, respectively start at this value of
$\alpha$. Within each of these segments we can choose arbitrary points
$\alpha^{ -} < \alpha$ and $\alpha^+ > \alpha$. For a single addition of entry $i$ to
the support we must have $r_i(\alpha^-)/w_i > \lambda'(\alpha^-)$, and
it therefore follows from \eqref{Eq:PQ} that
\[
\lambda'(\alpha^-) = \frac{\sum_{j\in\mathcal{I}}
  r_j(\alpha)}{\sum_{j\in\mathcal{I}} w_j^2}
<
\frac{w_ir_i + \sum_{j\in\mathcal{I}}
  w_jr_j(\alpha)}{w_i^2\ + \sum_{j\in\mathcal{J}} w_j^2\ \ \ \ \ } = \lambda'(\alpha^+).
\]
For the removal of a single entry $i$ from the support we must have
$r_i(\alpha^-)/w_i < \lambda'(\alpha^-)$ it likewise follows that
\[
\lambda'(\alpha^-) =
\frac{\sum_{j\in\mathcal{I}} w_jr_j(\alpha)}{\sum_{j \in \mathcal{I}}}
=
 \frac{w_ir_i(\alpha) +
  \sum_{j\in\mathcal{I}\setminus i}
  w_jr_j(\alpha)}{w_i^2 \ \ \ \ \!+ \sum_{j \in {\mathcal{I}\setminus i}}
  w_j^2\ \ \ \ }
<
\frac{\sum_{j\in\mathcal{I}\setminus i }
  w_jr_j(\alpha)}{\sum_{j \in \mathcal{I} \setminus i}w_j^2} = \lambda'(\alpha^+).
\]
Multiple simultaneous changes to the support can be dealt with one at
a time in a similar manner, resulting in $\lambda'(\alpha^-) <
\lambda'(\alpha^+)$. When $x(\alpha) \in\mathcal{C}_{w,1}$ we have
$\lambda(\alpha) = 0$. At the entry point we must have
$\lambda'(\alpha^-) < 0 = \lambda'(\alpha^+)$, and for the exit point
we have $\lambda(\alpha^-) = 0 < \lambda'(\alpha^+)$. Let slope $s_i =
d_i / w_i$, $s_{\max} = \max_j s_j$, and $k$ be any index such that
$s_k = s_{\max}$, then for sufficiently negative $\alpha$ we have
$v_k(\alpha) - v_i(\alpha) > \tau$ for all $i$ such that $s_i <
s_{\max}$. This implies that $[x_i(\alpha) - w_i\lambda(\alpha)]_+=0$,
and therefore that only those entries with the maximum slope can be in
the support, thus giving $\lambda'(\alpha) = -s_{\max}$. Similarly, we
have $\lambda'(\alpha) = s_{\max}$ for sufficiently large $\alpha$.
Summarizing we have that $\lambda(\alpha)$ is piecewise linear with
slopes strictly increasing from $-s_{\max}$ to $s_{\max}$, and we
therefore conclude that $\lambda(\alpha)$ is convex.

\subsection{Upper bound}

From the convexity of $\lambda(\alpha)$ it immediately follows that
the maximum number of intersections of $\lambda(\alpha)$ with any
$v_i(\alpha)$ is four whenever $s_i < s_{\max}$ and two whenever $s_i
= s_{\max}$. If $\lambda(\alpha)$ reaches zero there are two more
break points, but between these points there must be at least one
$v_i$ that reaches zero, thereby removing two possible intersections
with that curve. Since there is at least one index for which $s_i =
s_{\max}$, the maximum possible number of break points is therefore
$4(n-1) + 2 = 4n-2$. As each break point corresponds to a transition
from one face to the next, it follows that the maximum number of faces
of $\mathcal{C}_{w,1}$ that a line can project onto is $4n-1$.

\subsection{Constructions for the weighted one-norm ball}

For $n=1$ the maximum number of three faces is reached whenever $d_1
\neq 0$. For any $n \geq 2$ we can use the construction illustrated in
Figure~\ref{Fig:ProjectionPath}(b), consisting of two individual
curves and a bundle of $n-2$ curves in between. Working with zero
crossings $z_i$ instead of origin values $o_i$ we define the first
curve by $z_1 = 0$, $s_1 = 1$, and weight $w_1 = \omega$ to be
specified later. The second curve has $z_2 = 3$, $s_2 = s_{\max} = 4$,
and $w_2 = 1$. For each of the remaining $n-2$ curves we sample the
zero crossing $z_k$ i.i.d.~from $\mathcal{U}(1.9,2)$, and choose $s_k
= 2$ and $w_k = 1$. These values are chosen such that $v_i(4) \geq 4$.
The only two parameters that remain to be chosen are $\omega$ and
$\tau$, and the approach is then as follows. By choosing $\omega$
sufficiently large, the minimum of $\norm{x(\alpha)}_{w,1}$ occurs at
$\alpha = 0$.  We can then choose $\tau$ such that $\lambda(1) = 0$
forms a break point.  Along with a second zero crossing for some
$\alpha < 0$ and additional intersections for sufficiently small and
large values of $\alpha$, this gives a total of four
break points from the first curve. It then remains to ensure that
$\lambda(4) < 4$, in which case we have two intersections of
$\lambda(\alpha)$ with each of the remaining curves on the interval
$\alpha \in [1,4]$. The random sampling of $z_k$ for $k \geq 2$
ensures with probability one that no three curves cross at the same
point. For $k\geq 2$ we have $s_k < s_{\max}$, and each of these
curves will therefore have an additional two intersections for
sufficiently large positive and negative values of $\alpha$. This
gives the desired $4n-2$ break points and $4n-1$ faces.

The weighted one-norm $\norm{x(\alpha)}_{w,1}$ can be verified to be
equal to $\sum_i s_i\abs{\alpha - z_i} w_i^2$. The directional
derivative at $\alpha = 0$ is equal to $\omega^2 - 2n$, and for this
to be positive we need to choose $\omega > \sqrt{2n}$. It now remains
to choose an $\omega$ such that $\lambda(4) < 4$. A sufficient
condition for this is that the directional derivative
$\lambda'(\alpha) < 4/3$ for all $\alpha \in [1,4]$. Initially we have
$\lambda'(1) = \frac{\omega^2 - 2n}{\omega^2 + (n-1)} < 1$, and we
must therefore first intersect a curve with index $i \geq 2$. At the
first break point we have $\lambda' = \frac{c_1 + \omega^2}{c_2 + \omega^2}$
for some $c_1$ and $c_2$, and $\lambda'$ can therefore be kept smaller
than $4/3$ by choosing a sufficiently large $\omega$. This process can
be repeated until we arrive at $\alpha = 4$ by taking the largest
necessary $\omega$ over all steps.  Using \eqref{Eq:PQ} we find that
the largest combination of including or excluding each of the curves
$i\geq 2$ is attained by including all indices, giving $\lambda' =
 \frac{2n + \omega^2}{(n-1) + \omega^2}$. Choosing $\omega > \sqrt{2n +
  4}$ gives $\lambda' = \lambda'(4) < 4/3$, as
needed. Figure~\ref{Fig:ProjectionPath}(b) plots the result for
$\omega = \sqrt{2n+5}$.

\subsection{Constructions for the canonical one-norm ball}

For $n=1$ the maximum number of three faces is reached whenever $d_1
\neq 0$.  For $n=2$ the construction in
Figure~\ref{Fig:WeightedL1ProjArc} for a weighted one-norm ball
attains the upper bound of seven. For the canonical one-norm ball, it
is easily seen that the number of faces that $\lambda(\alpha)$ can
projected onto is at most five. An example construction that attains
the maximum for $n=3$ is plotted in
Figure~\ref{Fig:ProjectionPath3D4D}(a). The construction uses slopes
$d = [1.0,\ 0.5,\ 1.0\!-10^{-3}]$, and offsets $s$ such that the zero
crossings are at $z = [-4,0,4]$. The slopes of the right- and
left-most curves are chosen nearly identical to keep the sum of these
components nearly constant for $\alpha \in [-4,4]$. The middle curve
has a smaller slope and forces the one-norm to grow from the middle
outwards. Choosing $\tau = \norm{x-\alpha d}_1$ at $\alpha = 3$ causes
$\lambda$ to be zero between $\alpha \approx -3$ and $\alpha =3$. This
is close enough to the point where the outer curves reach zero to
ensure two intersections of $\lambda(\alpha)$ with each of these two
curves before it intersects the middle curve again (intersections
occur outside of the plotted range). The slope of the right-most curve
is chosen slightly less than the maximum to ensure that it will
eventually be intersected again by $\lambda(\alpha)$ at large positive
and negative values of $\alpha$.

\subsubsection{Four dimensions} 
\begin{figure}[t!]
\centering
\begin{tabular}{cc}
\includegraphics[width=0.455\textwidth]{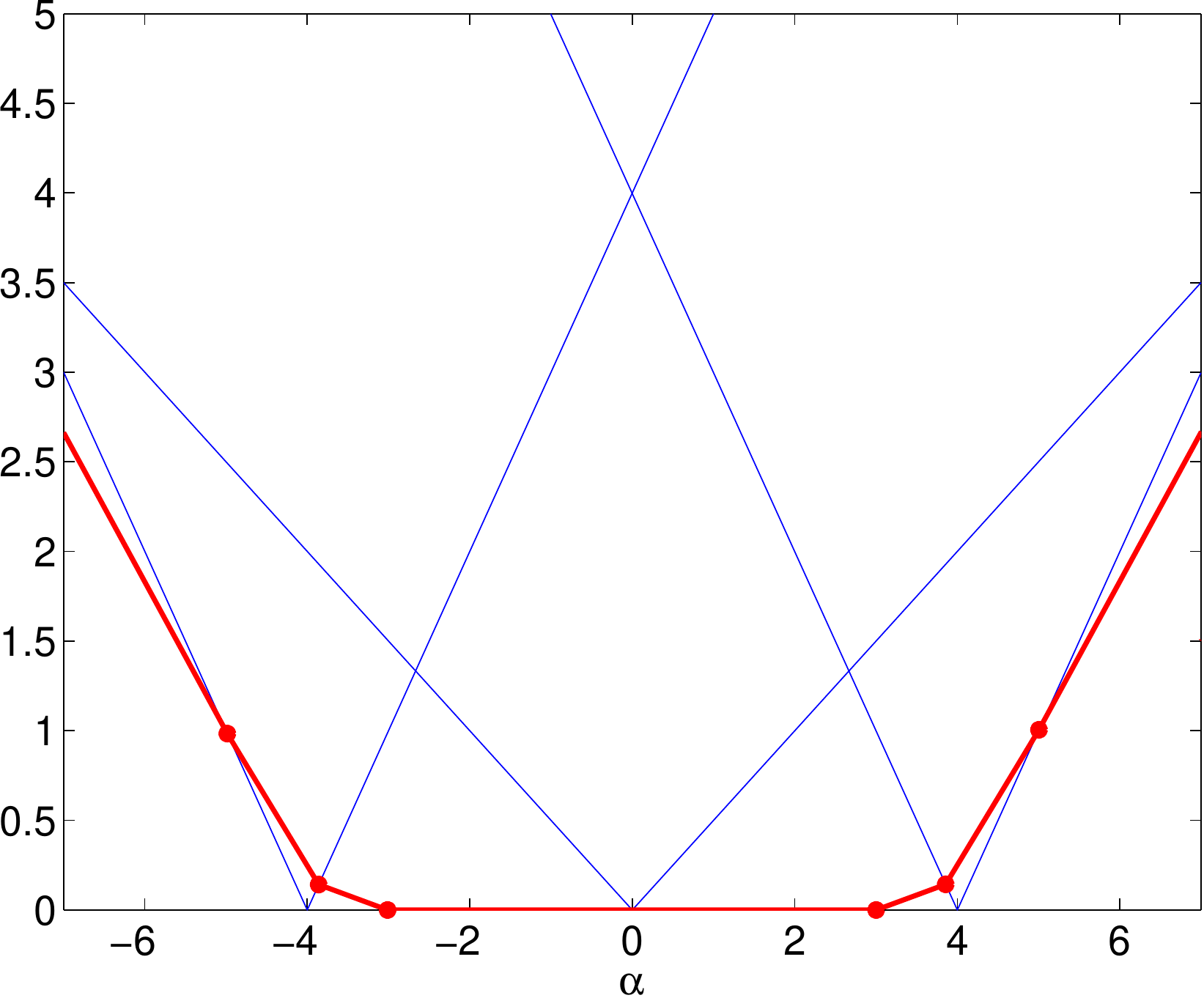}&
\includegraphics[width=0.455\textwidth]{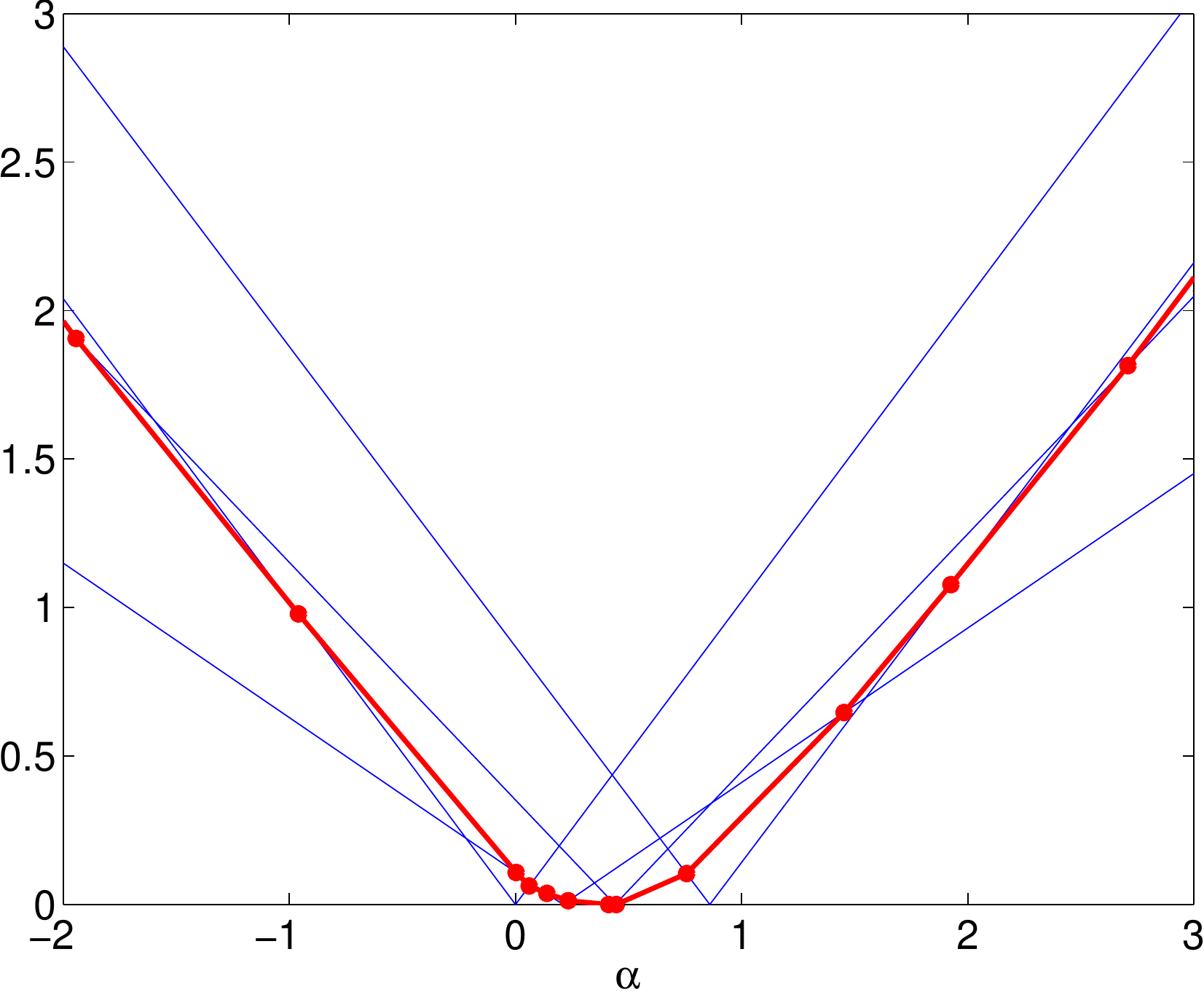}\\
({\bf{a}}) & ({\bf{b}})
\end{tabular}
\caption{ Examples that attains the upper bound of $4n-2$
  intersections for (a) three dimensions and (b) four dimensions (some
  intersections occur outside the plotted
  region).}\label{Fig:ProjectionPath3D4D}
\end{figure}

The construction of $n=4$ uses
parameters $d=[1.02, 0.52, 0.80, 1.01]$ with zero crossings at $z=
[0.00, 0.21, 0.44, 0.86]$ and is plotted in
Figure~\ref{Fig:ProjectionPath3D4D}(b). Some of the intersections occur
outside of the plotted range, but note that for curves with a slope
that is smaller than the maximum, it suffices to have
$\lambda(\alpha)$ below the curve on either side of two intersections;
the slope of $\lambda(\alpha)$ eventually has to match the maximum of
$d$ and must therefore cross. The right-most curve has a slope that is
just below the maximum slope and has its first two intersections at
$\alpha$ near $0.75$ and $2$.  The two intersections with left-most
and steepest curve occur at $\alpha$ close to $0.05$ and $-1$.

\subsubsection{Higher dimensional construction} 
\begin{figure}[t!]
\centering
\begin{tabular}{cc}
\multicolumn{2}{c}{\includegraphics[width=0.95\textwidth]{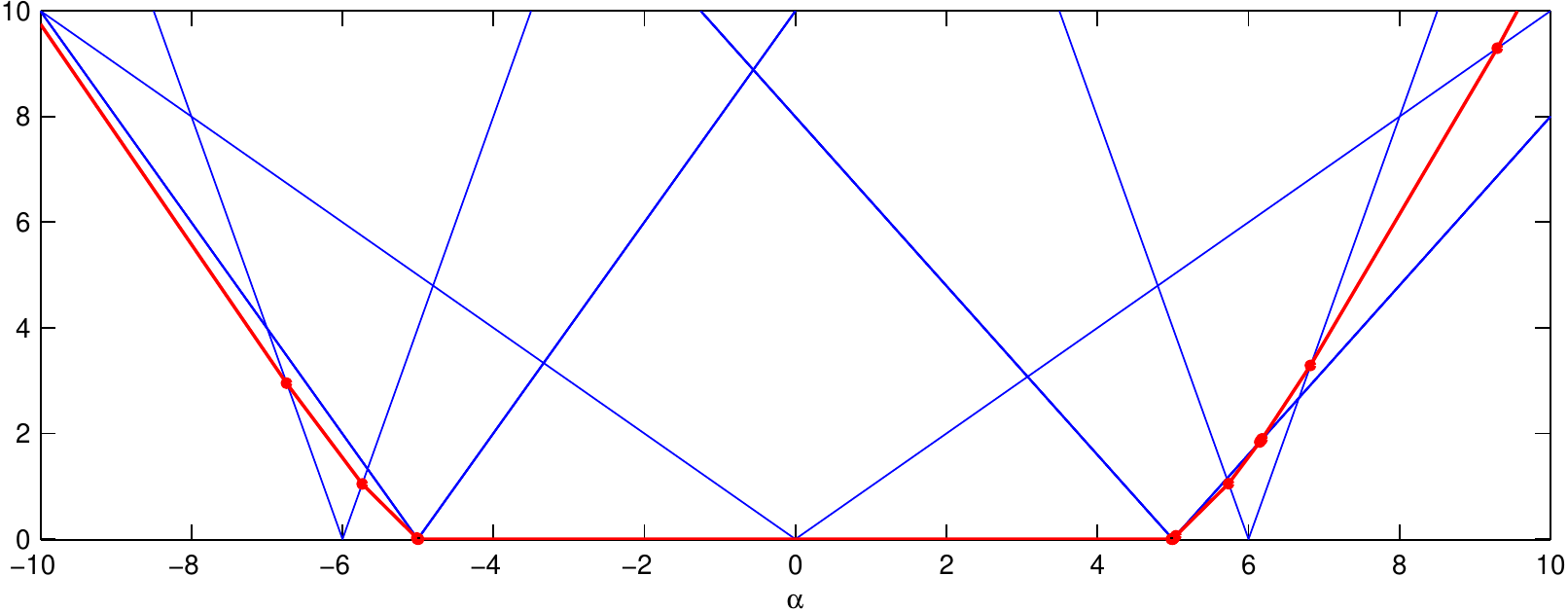}}\\
\multicolumn{2}{c}{({\bf{a}})}\\
\includegraphics[width=0.455\textwidth]{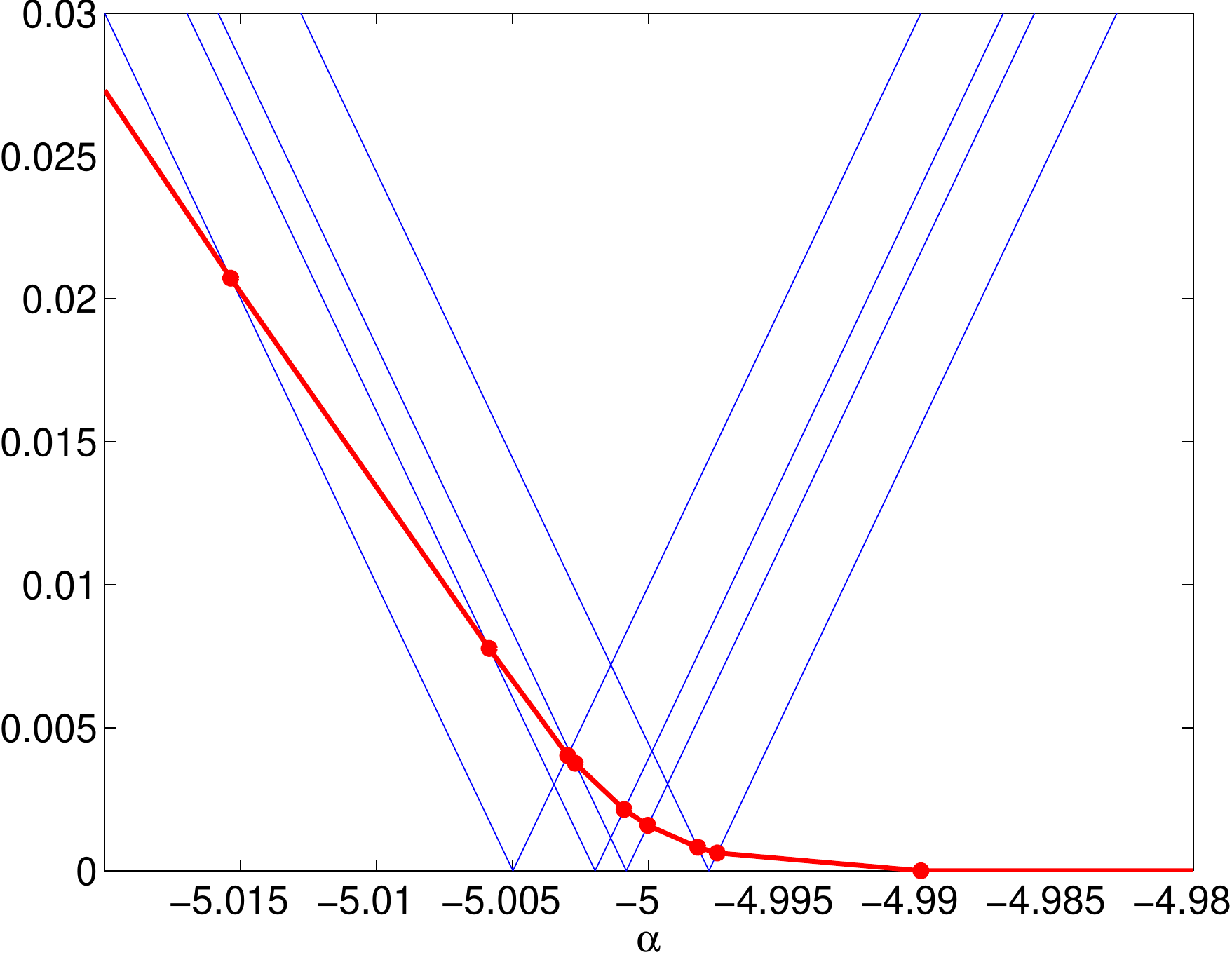}&
\includegraphics[width=0.455\textwidth]{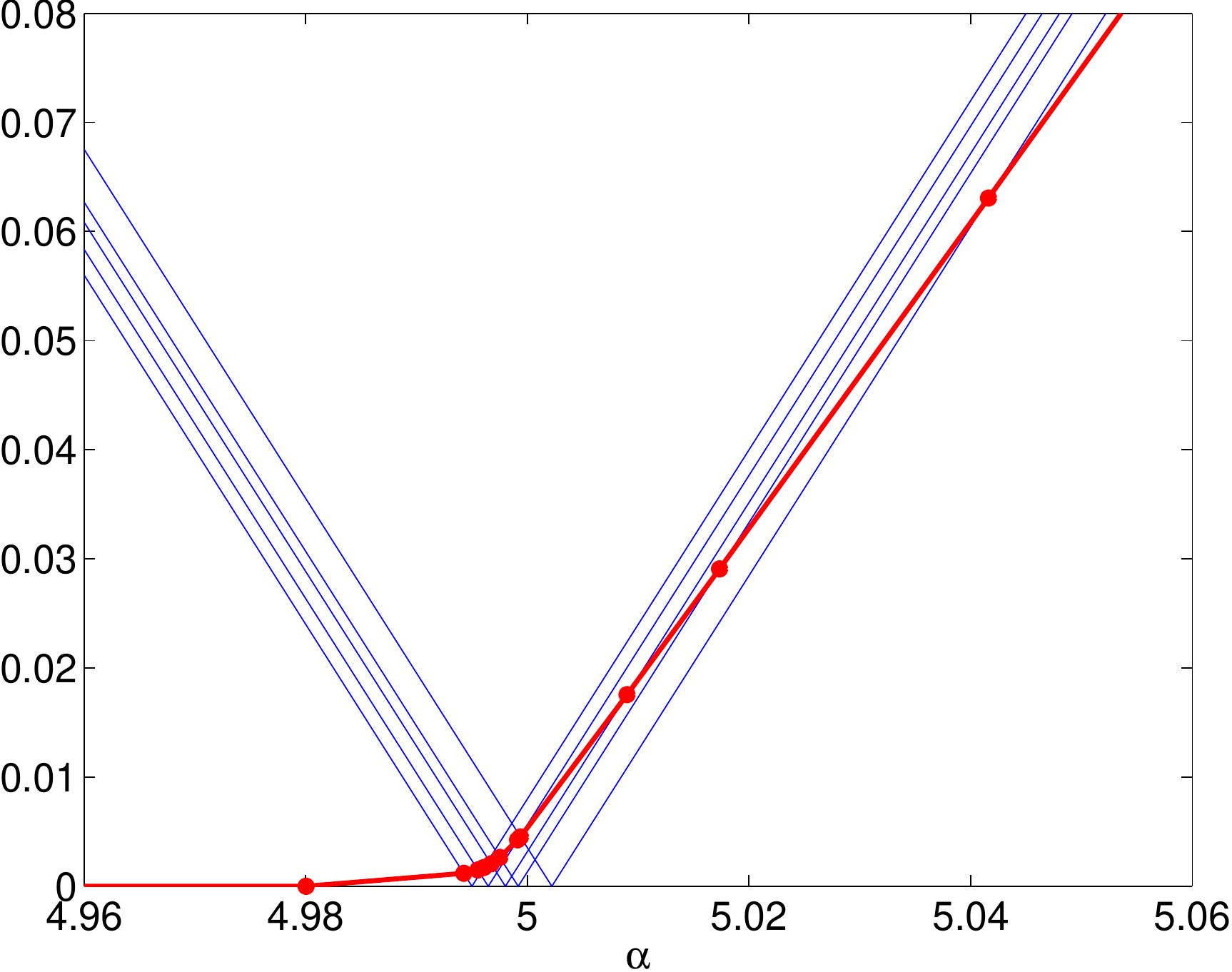}\\
({\bf{b}}) & ({\bf{c}})
\end{tabular}
\caption{Plots showing (a) the configuration for $n=12$, with $\beta =
  \sfrac{1}{2}$ and $\delta = 0.01$; and (b,c) the intersections of
  $\lambda(\alpha)$ with the center bundles around
  $\pm\mu_1$.}\label{Fig:ProjectionPathHighDim}
\end{figure}

For $n=5$ we can use the
general setup illustrated in
Figure~\ref{Fig:ProjectionPathHighDim}. It consists of two outer
curves crossing at $-\mu_2$ and $\mu_2$ with slopes $4$ and
$4-\epsilon$, respectively, for some small $\epsilon > 0$. Next there
is a bundle of $k_1 = \lfloor (n-3)/2\rfloor$ uniformly sampled i.i.d.~in
the interval $[-\sigma,\sigma]$ around $-\mu_1$ with slopes
$2$. Another bundle of $k_2 = \lceil (n-3)/2\rceil$ is uniformly
sampled i.i.d.~in the interval $[-\sigma,\sigma]$ around $\mu_1$, with
slopes $2k_1/k_2$. Finally there is a central curve with zero crossing
at $0$ and slope $1$. Aside from the $-\epsilon$ term (which is there
to ensure that only one curve attains the maximum slope) and excluding
the central curve, the one-norm remains constant between
$-\mu_1+\sigma$ and $\mu_1-\sigma$. This enables us to force $\lambda$
to be zero between $-\mu_1+2\delta$ and $\mu_1-2\delta$ and ensure a
total of four break points for the central curve. By choosing $\delta$
and $\sigma$ sufficiently close we then force two crossings of
$\lambda(\alpha)$ with the curves in each of the two bundles (and an
additional two crossings with each curve for sufficiently large
positive and negative values of $\alpha$). We then make sure to place
$\mu_2$ close to $\mu_1$ such that the two outer curves too are
intersected twice before crossing the central line again (this can be
done by choosing $\mu_1$ large enough and give the central curve
enough space to grow sufficiently large). The proof is done in a
number of steps:
\begin{description}
\item[Step 1.] Derive conditions such that $\lambda(\alpha)$ changes
  to or from $0$ at $-\mu_1 + \delta$ and in the interval $[\mu_1 -
  2\delta, \mu_1 - \delta]$;
\item[Step 2.] Determine conditions on $\sigma$ and $\delta$ under
  which $\lambda(\alpha)$ crosses all curves in each of the bundles
  and remains below them at a distance $\beta$ from $\pm\mu_1$;
\item[Step 3.] Determine $\mu_2$ relative to $\mu_1$ such that
  $\lambda(\alpha)$ crosses the outer curves;
\item[Step 4.] The entire construction allows us to change $\mu_1$
  (and $\mu_2$, $\epsilon$ accordingly) without changing the
  intersections of the bundles and outer curves (aside from minor
  effects due to $\epsilon$). In the last step we therefore choose
  $\mu_1$ to ensure that the central curve remains above
  $\lambda(\alpha)$ until after the outer two curves have been
  intersected.
\end{description}
We now consider each of the different steps.

\paragraph{Step 1 -- zero crossing}

We choose $\tau$ such that $\lambda(\alpha)$ reaches $0$ at $\alpha =
-\mu_1 + \delta$. This is done simply by equation $\tau$ to the sum of
the values of the curves at this point.  Because $\epsilon > 0$ the
sum of the curves at $\mu_1 - \delta$ will exceed $\tau$ and the zero
crossing of $\lambda(\alpha)$ must therefore occur before this
point. We now choose $\epsilon$ such that the one norm at
$\mu_1-2\delta$ is no greater than $\tau$.  The contribution of the
curves within each of the bundles affect the value of $\tau$, but
their contribution to the sum remains constant over the entire
interval $[-\mu_1+\sigma,\mu_1-\sigma]$ and can therefore be ignored.
Looking only at the relative difference we require that
\[
(\mu_1 - \delta)\cdot 1 \geq (\mu_1 - 2\delta)\cdot 1 + (2\mu_1 - 3\delta)\epsilon,
\]
which reduces to $\epsilon \leq \delta / (2\mu_1 - 3\delta)$, or the
sufficient condition that {\color{blue}$\epsilon \leq \delta / 2\mu_1$}.

\paragraph{Step 2 -- crossing the bundles} We analyze the crossing of
the bundles (see illustration in Figures~\ref{Fig:ProjectionPathHighDim}(b,c)) by
taking the maximum slope over the entire path and start from $-\mu_1 +
2\delta$ on the left and $\mu_1 - 2\delta$ on the right (this causes
the intersections to occur higher than they would otherwise). By
choosing $\mu_1$ large enough we can always ensure that the middle
curve remains above all intersections. Aside from this, the results in
this step are independent of $\mu_1$. We analyze the left and right
bundles in turn, starting from the left bundle. The maximum relevant
slope is directly to the left of the bundle and is equal to $(4k_1 + 1
- \epsilon) / n$, and we can use $(4k_1 + 1) / n$ for simplicity. We
want the value of $\lambda(\alpha)$ to be below the bundle at $\alpha
= -\mu_1 - \beta$ for some {\color{blue}$\beta > \sigma$}. This gives
\[
(\beta - \sigma)\cdot 2 \geq (\beta + 2\delta)\cdot (4k_1+1)/n
\]
or
\[
(2n - 4k_1 - 1)\beta \geq 2\delta(4k_1+1) + 2n\delta
\]
Choosing {\color{blue}$\sigma = \delta/2$} and using the fact that $k_1 \leq (n-3) /
2$ gives  the sufficient condition
\[
(2n - 2(n-3) - 1)\beta \geq 4\delta (n-3) + n\delta,
\]
which reduces to $\delta \leq 5/(5n-3)\beta$, which certainly holds
whenever {\color{blue}$\delta \leq \beta/n$}. Because the zero
crossing of the curves within the bundle are chosen uniformly at
random it holds with probability one that $\lambda(\alpha)$ does not
cross at any of the intersections between the lines in the bundle.

For the right bundle we find a maximum slope of $(4k_1 + 1 + \epsilon)
/ n$, and to guaranteed $\lambda(\alpha)$ to be below the bundle at
$\alpha = \mu_1+\beta$ we require
\[
(\beta - \sigma)\cdot 2k_1/k_2 \geq (\beta + 2\delta)\cdot (4k_1+1+\epsilon)/n,
\]
or, using $\sigma=\delta/2$,
\begin{equation}\label{Eq:RightBundle}
(2n k_1/k_2 - (4k_1 + 1 + \epsilon))\cdot \beta \geq 2\delta(4k_1 + 1
+ \epsilon) + n(k_1/k_2)\delta.
\end{equation}
For $k_1=k_2$ we have $k_1 = (n-3)/2$, and
\begin{equation}\label{Eq:BetaK1=K2}
(2nk_1/k_2 - (4k_1+1+\epsilon))\beta = (2n - ((2n-6)+1+\epsilon))\beta = (5 - \epsilon)\beta
\geq 4\beta,
\end{equation}
for {\color{blue}$\epsilon \leq 1$}. For the right hand side of
\eqref{Eq:RightBundle} we use $2k_1 = n-3$, and have
\begin{equation}\label{Eq:DeltaK1=K2}
((8k_1+2+2\epsilon) + n(k_1/k_2))\delta = ((4n-12) + 2 + 2\epsilon +
n)\delta = (5n-10+2\epsilon)\delta \leq 5n\delta
\end{equation}
Combining \eqref{Eq:BetaK1=K2} and \eqref{Eq:DeltaK1=K2} gives the
sufficient condition $\delta \leq 4\beta/5n$. For the case where $k_2
= k_1+1$ we can use $k_1 = (n-4)/2$, and $2n/k_2 \leq 6$.  For the
left-hand side of \eqref{Eq:RightBundle} this gives
\begin{eqnarray*}
(2nk_1/k_2 - (4k_1+1+\epsilon))\beta & = &
(2nk_1/k_2 - (2n-8 + 1 +\epsilon))\beta\\
& = &
(2nk_1/k_2 - 2n(k_1+1)/k_2 + 7-\epsilon)\beta\\
& \geq & (- 2n/k_2 + 7 -\epsilon)\beta \geq (1-\epsilon)\beta
\end{eqnarray*}
For the right-hand side of \eqref{Eq:RightBundle} we have
\begin{eqnarray*}
2\delta(4k_1 + 1
+ \epsilon) + n(k_1/k_2)\delta & \leq & (8k_1 + 2 + 2\epsilon)\delta +
6k_1\delta \\
& = & (14k_1 + 2 + 2\epsilon)\delta = (7n - 28 + 2 + 2\epsilon) \delta
\\
& \leq & 7n\delta
\end{eqnarray*}
This gives a sufficient condition of $\delta \leq \beta / 8n$ for
{\color{blue}$\epsilon \leq 1/8$}. This is the strongest condition of
the three cases, and we can therefore choose $\beta = 8n\delta $.

\paragraph{Step 3 -- crossing the outer curves}

First we have to find a minimum distance between $\mu_1$ and $\mu_2$
such that the outer curves are above $\lambda(-\mu_1 -\beta) \leq
2\beta$ and $\lambda(\mu_1+\beta) \leq 2\beta k_1/k_2 \leq 2\beta$.
Taking the smaller slope of $4-\epsilon$ it suffices to have
\[
(\mu_2 - \mu_1 - \beta)\cdot(4-\epsilon) \geq
2\beta,\quad\mbox{or}\quad
\frac{6-\epsilon}{4-\epsilon}\beta \leq \mu_2 - \mu_1.
\]
For {\color{blue}$\epsilon \leq 1$} it then suffices to have
$(6/3)\beta \leq \mu_2 - \mu_1$, which is satisfied for {\color{blue}
  $\beta = (\mu_2 - \mu_1) /2$}.

We now consider the intersection of $\lambda(\alpha)$ with the outer
curves, starting with the left (steepest) curve. Because we only need
to intersect a single curve twice, we can use the maximum slope before
the second intersection. The slope before the first crossing is less
than $2$ since we just intersected the bunch around $-\mu_1$. After
the first crossing with the outer curve at $-\mu_2$ the slope changes
to $(4k_1 + 1 + (4-\epsilon)) / (n-1)$. For $k_1 = k_2$ we have $2k_1
= n-3$ and
\[
\frac{4k_1 + 1 + (4-\epsilon)}{n-1} = \frac{2n - 6 + 5 - \epsilon}{n-1}
= \frac{2(n-1) + 1 - \epsilon}{n-1} > 2,
\]
whenever {\color{blue}$\epsilon < 1$}. This means that it can
intersect the curves from the left bundle before the second
intersection and for the maximum slope we therefore presume all of
these curves are crossed, giving a slope of
\[
\frac{2k_1 + 1 + (4-\epsilon)}{k_2+2} = \frac{2(k_1+2) + 1 -
  \epsilon}{k_1+2} \leq 2 + \frac{1-\epsilon}{k_1+2} \leq 2 + \sfrac{1}{3}.
\]
For $k_2 = k_1 + 1$ we have $2k_1 = n-4$ and
\[\frac{4k_1 + 1 + (4-\epsilon)}{n-1} = \frac{2n - 8 + 5 -
  \epsilon}{n-1}
= \frac{2(n-1) - 1 - \epsilon}{n-1} < 2.
\]
This means that the curves from the bundle will not be intersected
before the second intersection with the left-most curve.

For the right-most curve, we have a slope of $(4k_1+1+4)/(n-1)$
directly after the first intersection. For $k_1 = k_2$ we can follow
the same derivation as above (in this case with $\epsilon=0$) to find
that the maximum slope, after intersecting the curves from the right
bundle is bounded by $2+\sfrac{1}{3}$. For $k_2 = k_1+1$ we use $2k_2
= n-2$ and find
\[
\frac{4k_1 + 1 + 4}{n-1} = \frac{4k_1 + 5}{n - 1} > \frac{4k_1}{n - 2}
= 2\frac{k_1}{k_2},
\]
so again, the curves from the right bundle can be
intersected. Assuming all of them are intersected before the second
intersection with the right-most curve we have a slope bounded again
by $2 + \sfrac{1}{3}$.

We now take the maximum slope to be $3$ and compute the maximum
distance $\gamma$ from $\pm \mu_2$ for which the second intersection
with the outer curves occurs. It suffices to work with the right-most
curve. We can start at $\alpha = \mu_1+\beta$, the minimum curve of the
right bundle can be below $2\beta$, which immediately gives
$\lambda(\alpha) \leq 2\beta$. Taken together we need to find $\gamma$
such that
\[
2\beta + 3(\mu_2 - (\mu_1 + \beta) + \gamma) \leq (4-\epsilon)\gamma 
\]
Reorganizing gives
\[
3(\mu_2 - \mu_1) - \beta \leq (1-\epsilon)\gamma
\]
With the choice of $\beta = (\mu_2 - \mu_1) / 2$ this simplifies to
$5(\mu_2 - \mu_1) \leq 2(1-\epsilon)\gamma$. Since we chose $\epsilon
\leq 1/8$ we have $\epsilon \leq 1/5$ and therefore it suffices to
have $(\mu_2-\mu_1) \leq 8\gamma$, which allows us to
choose {\color{blue}$\gamma = 4\beta$}.

\paragraph{Step 4 -- Controlling $\mu_1$} We can satisfy all required
inequalities thus far by fixing $\beta > 0$ and choosing $\mu_1 >
0$. In particular we can set $\mu_2 = \mu_1 + 2\beta$, and choose
$\delta = \beta / n$. Since we already fixed $\sigma = \delta / 2$ and
$\gamma = 4\beta$, it suffices to choose $\epsilon = \min\{ \delta /
2\mu_1, 1/8\}$, which always gives $\epsilon > 0$. For the
construction to work we must make sure that the central curve is above
$\lambda(\alpha)$ for $\alpha = \pm(\mu_2 + \gamma)$. It suffices to
to be above the outer curves at this point and therefore
\[
1\cdot (\mu_2 + \gamma) \geq 4\cdot \gamma.
\]
Expanding $\gamma$ and $\mu_2$ gives $\mu_1 + 6\beta \geq 16\beta$,
which is satisfied for {\color{blue}$\mu_1 = 10\beta$}. The plot in
Figure~\ref{Fig:ProjectionPathHighDim} illustrates the construction
for $n=12$.

\section{Line search}

\subsection{First Wolfe condition}\label{AppSec:Wolfe1}

We want
\begin{equation}\label{Eq:DerivationCondOne}
f(x+\alpha d) \leq f(x) + \gamma_1 \alpha(\nabla f(x))^Td.
\end{equation}
Expanding the left-hand side gives
\begin{eqnarray*}
f(x+\alpha d) & = & 
\half\norm{Ax-b+\alpha A d}_2^2 + \frac{\mu}{2}\norm{x+\alpha d}_2^2 +
c^T(x+\alpha d) \\
& = & \half\norm{ r + \alpha A d}_2^2 + \frac{\mu}{2}\norm{x}_2^2 + \mu
\alpha x^Td + \frac{\mu}{2}\alpha^2\norm{d}_2^2 + c^Tx + \alpha c^Td
\\
& = & \half\norm{r}_2^2 + \alpha r^TAd + \half\alpha^2\norm{Ad}_2^2
+ \frac{\mu}{2}\norm{x}_2^2 + \mu
\alpha x^Td + \frac{\mu}{2}\alpha^2\norm{d}_2^2 + c^Tx + \alpha c^Td
\end{eqnarray*}
For the right-hand side we have
\[
f(x) + \gamma_1 \alpha(\nabla f(x))^Td = \half \norm{r}_2^2 +
\frac{\mu}{2}\norm{x}_2^2 + c^Tx + \gamma_1\alpha(A^Tr + \mu x + c)^Td
\]
To satisfy \eqref{Eq:DerivationCondOne} we subtract the two
\begin{eqnarray*}
  f(x+\alpha d) - f(x) - \gamma_1 \alpha(\nabla f(x))^Td &\leq & 0\\
  \alpha r^TAd + \half\alpha^2\norm{Ad}_2^2
  + \mu \alpha x^Td + \frac{\mu}{2}\alpha^2\norm{d}_2^2  + \alpha c^Td
- \gamma_1\alpha(A^Tr + \mu x + c)^Td & \leq & 0
\end{eqnarray*}
Regrouping the terms gives
\begin{eqnarray*}
(1-\gamma_1) \alpha r^TAd
+\half\alpha^2\norm{Ad}_2^2
+ \mu (1-\gamma_1) \alpha x^Td 
+ \frac{\mu}{2}\alpha^2\norm{d}_2^2
+  (1- \gamma_1)\alpha c^Td
& \leq & 0 \\
\alpha^2(\half\norm{Ad}_2^2 + \frac{\mu}{2}\norm{d}_2^2)
+ \alpha (1-\gamma_1)\cdot( r^TAd
+ \mu x^Td  +  c^Td) & \leq & 0
\end{eqnarray*}
Dividing by $\alpha$ and rearranging gives
\[
\alpha \leq -(1 - \gamma_1)\frac{r^TAd + \mu x^Td + c^Td}{\half
  \norm{Ad}_2^2 + \half\mu\norm{d}_2^2} 
= -2(1 - \gamma_1)\frac{(A^Tr + \mu x + c)^Td}{
  \norm{Ad}_2^2 + \mu\norm{d}_2^2}  = 2(1-\gamma_1)\alpha_{\mathrm{opt}}
\]
\subsection{Second Wolfe condition}\label{AppSec:Wolfe2}

We require
\[
(\nabla f(x + \alpha d))^Td \geq \gamma_w(\nabla f(x))^Td
\]
Expanding gives
\begin{eqnarray*}
[A^T(r+\alpha Ad) + \mu(x+\alpha d) + c]^Td & \geq & \gamma_2(A^Tr
+ \mu x + c)^T d \\
\alpha \norm{Ad}_2^2 + \mu\alpha \norm{d}_2^2 & \geq & (\gamma_2 - 1)(A^Tr
+ \mu x + c)^T d\\
\alpha & \geq  & -(1-\gamma_2)\frac{(A^Tr
+ \mu x + c)^T d}{\norm{Ad}_2^2 + \mu\alpha \norm{d}_2^2}\\
\alpha & \geq & (1-\gamma_2)\alpha_{\mathrm{opt}}
\end{eqnarray*}


\section*{Acknowledgments}

This work was supported by National Science Foundation Grant DMS 0906812
(American Reinvestment and Recovery Act).  The implementation of the
code uses L-BFGS update routines kindly provided by Michael
P. Friedlander.

\bibliography{articles}

\end{document}